\documentclass[a4paper,12pt]{article}

\usepackage{amsmath,amssymb,amsthm}

\newcommand{\Fg}{\mathfrak{g}}
\newcommand{\Fh}{\mathfrak{h}}

\newcommand{\BB}{\mathbb{B}}
\newcommand{\BP}{\mathbb{P}}
\newcommand{\BZ}{\mathbb{Z}}
\newcommand{\BQ}{\mathbb{Q}}

\newcommand{\CB}{\mathcal{B}}

\newcommand{\wt}{\mathop{\rm wt}\nolimits}
\newcommand{\cl}{\mathop{\rm cl}\nolimits}
\newcommand{\af}{\mathop{\rm aff}\nolimits}
\newcommand{\id}{\mathop{\rm id}\nolimits}
\newcommand{\dist}{\mathop{\rm dist}\nolimits}
\newcommand{\Hom}{\mathop{\rm Hom}\nolimits}
\newcommand{\GL}{\mathop{\rm GL}\nolimits}
\newcommand{\Bij}{\mathop{\rm Bij}\nolimits}
\newcommand{\turn}{\mathop{\rm Turn}\nolimits}
\newcommand{\supp}{\mathop{\rm supp}\nolimits}
\newcommand{\Supp}{\mathop{\rm Supp}\nolimits}
\newcommand{\Pos}{\mathop{\rm Pos}\nolimits}
\newcommand{\Neg}{\mathop{\rm Neg}\nolimits}
\newcommand{\sgn}{\mathop{\rm sgn}\nolimits}
\newcommand{\img}{\mathop{\rm Im}\nolimits}
\newcommand{\fp}{\mathop{\rm fin}\nolimits}

\newcommand{\rr}{\Delta^{\mathrm{re}}}
\newcommand{\prr}{\Delta^{\mathrm{re}}_{+}}
\newcommand{\fr}{\fin{\Delta}}
\newcommand{\pfr}{\fin{\Delta}_{+}}

\newcommand{\ve}{\varepsilon}
\newcommand{\vpi}{\varpi}
\newcommand{\vp}{\varphi}

\newcommand{\bzero}{{\bf 0}}

\newcommand{\fin}[1]{\overset{\circ}{#1}}
\newcommand{\ba}[1]{\overline{#1}}
\newcommand{\ud}[1]{\underline{#1}}
\newcommand{\ti}[1]{\widetilde{#1}}
\newcommand{\ha}[1]{\widehat{#1}}

\makeatletter
\@addtoreset{equation}{subsection}

\renewcommand\section{\@startsection{section}{1}{0pt}
{-3.5ex plus -1ex minus -.2ex}{1.0ex plus .2ex}{\large\bf}}
\renewcommand\subsection{\@startsection{subsection}{1}{0pt}
{2.5ex plus 1ex minus .2ex}{-1em}{\bf}}
\makeatother

\newcommand{\vsp}{\vspace{1.5mm}}

\theoremstyle{plain}
\newtheorem{thm}{Theorem}[subsection]
\newtheorem{lem}[thm]{Lemma}
\newtheorem{prop}[thm]{Proposition}
\newtheorem{cor}[thm]{Corollary}

\newtheorem{claim}{Claim}[thm]

\newtheorem*{aclaim}{Claim}

\theoremstyle{definition}
\newtheorem{dfn}[thm]{Definition}

\theoremstyle{remark}
\newtheorem{rem}[thm]{Remark}
\newtheorem{eg}[thm]{Example}

\begin{document}

\setlength{\baselineskip}{20pt}

\title{{\Large\bf Crystal structure of 
the set of Lakshmibai-Seshadri \\ 
paths of a level-zero shape for an affine Lie algebra}}
\author{
 Satoshi Naito \\ 
 \small Institute of Mathematics, University of Tsukuba, \\
 \small Tsukuba, Ibaraki 305-8571, Japan \ 
 (e-mail: {\tt naito@math.tsukuba.ac.jp})
 \\[2mm] and \\[2mm]
 Daisuke Sagaki \\ 
 \small Institute of Mathematics, University of Tsukuba, \\
 \small Tsukuba, Ibaraki 305-8571, Japan \ 
 (e-mail: {\tt sagaki@math.tsukuba.ac.jp})
}
\date{}
\maketitle

%
\begin{abstract} \setlength{\baselineskip}{16pt}
Let $\lambda = \sum_{i \in I_{0}} m_{i} \vpi_{i}$, 
with $m_{i} \in \BZ_{\ge 0}$ for $i \in I_{0}$, 
be a level-zero dominant integral weight 
for an affine Lie algebra $\Fg$ over $\BQ$, 
where the $\vpi_{i}$, $i \in I_{0}$, are 
the level-zero fundamental weights, and 
let $\BB(\lambda)$ be the crystal of 
all Lakshmibai-Seshadri paths of shape $\lambda$.
First, we give an explicit description of 
the decomposition of the crystal $\BB(\lambda)$ 
into a disjoint union of connected components, and 
show that all the connected components are 
pairwise ``isomorphic'' (up to a shift of weights).
Second, we ``realize'' the connected component 
$\BB_{0}(\lambda)$ of $\BB(\lambda)$ containing 
the straight line $\pi_{\lambda}$ 
as a specified subcrystal of 
the affinization $\ha{\BB(\lambda)_{\cl}}$ 
(with weight lattice $P$) of the crystal 
$\BB(\lambda)_{\cl} \simeq 
\bigotimes_{i \in I_{0}} 
\bigl(\BB(\vpi_{i})_{\cl}\bigr)^{\otimes m_{i}}$ 
(with weight lattice $P_{\cl} = P/(\BQ \delta \cap P)$, 
where $\delta$ is the null root of $\Fg$), 
which was studied in a previous paper \cite{NSz}.
\end{abstract}
%
%
%
%
\section{Introduction.}
\label{sec:intro}
%
%
Lakshmibai-Seshadri paths (LS paths for short) and 
root operators acting on them 
for a Kac-Moody algebra $\Fg$ over the field $\BQ$ were 
introduced by Littelmann in \cite{L1}, \cite{L2}, and 
it was proved independently by Kashiwara \cite{Kassim} and 
Joseph \cite{J} that 
if $\lambda$ is a dominant integral weight, then 
the set (or rather, crystal) $\BB(\lambda)$ of 
all LS paths of shape $\lambda$ is 
isomorphic as a crystal 
to the crystal base $\CB(\lambda)$ of 
the irreducible highest weight module $V(\lambda)$ of 
highest weight $\lambda$ over the quantized universal 
enveloping algebra (quantum affine algebra) $U_{q}(\Fg)$ of 
$\Fg$ over the field $\BQ(q)$.
But, the crystal structure of the set $\BB(\lambda)$ of 
all LS paths of shape $\lambda$ 
for a general integral weight $\lambda \in P$ remained unknown.
The purpose of this paper is 
to study the crystal structure of 
this $\BB(\lambda)$ in the case that $\Fg$ 
is an affine Lie algebra.

Let us describe our results more explicitly.
Let $\Fg$ be an affine Lie algebra over $\BQ$ 
with Cartan subalgebra $\Fh$, 
simple roots $\bigl\{\alpha_{j}\bigr\}_{j \in I} \subset \Fh^{\ast}$, 
simple coroots $\bigl\{h_{j}\bigr\}_{j \in I} \subset \Fh$, and 
Weyl group $W \subset \GL(\Fh^{\ast})$.
We denote by 
$\delta = \sum_{j \in I} a_{j} \alpha_{j} \in \Fh^{\ast}$ 
the null root, and by 
$c = \sum_{j \in I} a^{\vee}_{j} h_{j} \in \Fh$ 
the canonical central element.
Note that in the case of an affine Lie algebra, 
an integral weight $\lambda \in P$ is 
of positive level (i.e., $\lambda(c) > 0$), 
of negative level (i.e., $\lambda(c) < 0$), or 
of level zero (i.e., $\lambda(c) = 0$).
If $\lambda \in P$ is of positive (resp., negative) level, 
then $\lambda$ is equivalent, under the Weyl group $W$ of $\Fg$, 
to a dominant (resp., antidominant) integral weight $\Lambda \in P$, 
and hence 
the $\BB(\lambda) = \BB(\Lambda)$ 
(see Remark~\ref{rem:LS}\,(3)) is 
isomorphic as a crystal to the crystal base $\CB(\Lambda)$ of 
the irreducible highest (resp., lowest) weight module of 
highest (resp., lowest) weight $\Lambda$ 
over the quantum affine algebra $U_{q}(\Fg)$.
Thus it remains to consider the case 
in which $\lambda$ is of level zero, and 
hence (by Remark~\ref{rem:LS}\,(3) and 
Lemma~\ref{lem:d-shift4}) level-zero dominant, i.e., 
of the following form: 
$\lambda = 
 \sum_{i \in I_{0}} m_{i} \vpi_{i}$, 
with $m_{i} \in \BZ_{\ge 0}$ for $i \in I_{0}$, 
where the $\vpi_{i}$, $i \in I_{0} := I \setminus \{0\}$, are 
the level-zero fundamental weights 
for the affine Lie algebra $\Fg$.

In previous papers \cite{NSp1}, \cite{NSp2}, 
we proved that for each $i \in I_{0}$, 
the crystal $\BB(\vpi_{i})$ of 
all LS paths of shape $\vpi_{i}$ is connected, and 
also proved that this crystal $\BB(\vpi_{i})$ is 
isomorphic to the crystal base $\CB(\vpi_{i})$ of 
the extremal weight module $V(\vpi_{i})$ of 
extremal weight $\vpi_{i}$ over $U_{q}(\Fg)$, 
introduced by Kashiwara \cite{Kasmod}, \cite{Kaslz} 
in a more general setting.
Furthermore, in \cite{NSp2}, 
we obtained an explicit description of 
the decomposition of the crystal $\BB(m \vpi_{i})$ of 
all LS paths of shape $m \vpi_{i}$, 
with $m \in \BZ_{\ge 2}$, into 
a disjoint union of connected components, 
and then proved that also for $m \in \BZ_{\ge 2}$, 
this crystal $\BB(m \vpi_{i})$ is isomorphic to 
the crystal base $\CB(m \vpi_{i})$ of 
the extremal weight module $V(m \vpi_{i})$ of 
extremal weight $m \vpi_{i}$.
Therefore, it seems natural to expect that 
for an arbitrary level-zero dominant integral weight 
$\lambda = \sum_{i \in I_{0}} m_{i} \vpi_{i}$, 
with $m_{i} \in \BZ_{\ge 0}$ for $i \in I_{0}$, 
the crystal $\BB(\lambda)$ is isomorphic to 
the crystal base $\CB(\lambda)$ of 
the extremal weight module $V(\lambda)$ of 
extremal weight $\lambda$.
However, this is not true even for 
the (simple) case in which $\Fg$ is 
of type $A_2^{(1)}$ and 
$\lambda = \varpi_{1} + \varpi_{2}$, 
as can be seen from \cite[Example~4.1]{NSz} 
(and \cite[Remark~5.2]{NSp1}).
In fact, the crystals $\BB(\lambda)$ and $\CB(\lambda)$ are 
isomorphic only if $\lambda = m_{i} \vpi_{i}$ 
for some $i \in I_{0}$ and 
$m_{i} \in \BZ_{\ge 0}$ (see Appendix).

In this paper, 
we give an explicit description (Theorem~\ref{thm:comps}) of 
the decomposition of the crystal $\BB(\lambda)$ of 
all LS paths of shape $\lambda$ into 
a disjoint union of connected components, 
where $\lambda \in P$ is an arbitrary 
level-zero dominant integral weight of 
the form $\lambda = \sum_{i \in I_{0}} m_{i} \vpi_{i}$, 
with $m_{i} \in \BZ_{\ge 0}$ for $ i \in I_{0}$.
Namely, we prove that each connected component of 
the crystal $\BB(\lambda)$ contains 
exactly one (extremal) LS path having an expression of 
the form \eqref{eq:comps}, and 
that the set of all extremal LS paths 
in each connected component of 
the crystal $\BB(\lambda)$ coincides 
with the $W$-orbit of an (extremal) LS path 
having an expression of the form \eqref{eq:comps}.
A key to the proof of Theorem~\ref{thm:comps} is 
one of main results of a previous paper \cite{NSz} 
(see Theorem~\ref{thm:NSz}) that 
the crystal $\BB(\lambda)_{\cl}$ 
with weight lattice $P_{\cl}$ is ``simple'', and 
hence connected.
Here the lattice $P_{\cl}$ is 
equal to $P/(\BQ \delta \cap P)$, and 
$\BB(\lambda)_{\cl}$ is a crystal 
with weight lattice $P_{\cl}$ 
(which we call a $P_{\cl}$-crystal) obtained 
from the crystal $\BB(\lambda)$ 
with weight lattice $P$ 
(which we call a $P$-crystal) as follows: 
$\BB(\lambda)_{\cl} = 
 \bigl\{
  \cl(\pi) \mid \pi \in \BB(\lambda)
 \bigr\}$, 
where $(\cl(\pi))(t) \in P_{\cl}$ 
for $t \in [0,1]$ is defined to be 
the image $\cl(\pi(t))$ of $\pi(t)$ 
under the canonical projection 
$\cl:\BQ \otimes_{\BZ} P \twoheadrightarrow 
 \BQ \otimes_{\BZ} P_{\cl}$.

It immediately follows from 
Theorem~\ref{thm:comps}, along with Lemma~\ref{lem:d-shift2}, 
that the connected components of $\BB(\lambda)$ are 
pairwise ``isomorphic'' (up to a shift of weights), 
like the case of the crystal base $\CB(\lambda)$ of 
the extremal weight module $V(\lambda)$ of 
extremal weight $\lambda$ (see \cite[Theorem~4.15]{BN}).
Hence it suffices to study the crystal structure of 
the connected component $\BB_{0}(\lambda)$ of 
$\BB(\lambda)$ containing the straight line $\pi_{\lambda}$, 
which is obviously an extremal LS path 
having an expression of the form \eqref{eq:comps}.
For this purpose, 
generalizing \cite[Proposition~5.9]{GL} 
for the case in which $\lambda = m_{i} \vpi_{i}$, 
with $i \in I_{0}$ and $m_{i} \in \BZ_{\ge 1}$, 
to the case of the level-zero dominant integral weight 
$\lambda = \sum_{i \in I_{0}} m_{i} \vpi_{i}$, 
we introduce the ``affinization'' 
$\ha{\BB(\lambda)_{\cl}}$ (which is a $P$-crystal) of 
the $P_{\cl}$-crystal $\BB(\lambda)_{\cl}$, and 
prove that the affinization $\ha{\BB(\lambda)_{\cl}}$ is 
isomorphic as a $P$-crystal to a disjoint union of 
connected components of $\BB(\lambda)$, 
including $\BB_{0}(\lambda)$ (see Theorem~\ref{thm:isom}).
Furthermore, we give a condition 
(condition (C) of Corollary~\ref{cor:con}) 
for an element of the affinization 
$\ha{\BB(\lambda)_{\cl}}$ to 
lie in the isomorphic image of 
the connected component $\BB_{0}(\lambda)$ of 
$\BB(\lambda)$, which enables us to identify 
$\BB_{0}(\lambda)$ with a specified subcrystal of 
$\ha{\BB(\lambda)_{\cl}}$.  
Here we recall from \cite{NSp2} and \cite{NSz} 
(see Theorem~\ref{thm:NSz} and the comments just below it) 
that the $P_{\cl}$-crystal $\BB(\lambda)_{\cl}$ is 
isomorphic to the tensor product 
$\bigotimes_{i \in I_{0}} 
 \bigl(\BB(\vpi_{i})_{\cl}\bigr)^{\otimes m_{i}}$ of 
the $P_{\cl}$-crystals 
$\BB(\vpi_{i})_{\cl}$, $i \in I_{0}$, and 
that for each $i \in I_{0}$, 
the $P_{\cl}$-crystal $\BB(\vpi_{i})_{\cl}$ is 
isomorphic to the crystal base of 
the level-zero fundamental module $W(\vpi_{i})$ 
over the quantized universal enveloping algebra 
$U_{q}^{\prime}(\Fg)$ of $\Fg$ 
with weight lattice $P_{\cl}$, 
which is introduced by Kashiwara in \cite{Kaslz}.
Thus, our results in this paper, 
together with this description of 
the $P_{\cl}$-crystal structure of $\BB(\lambda)_{\cl}$, 
completely determine the $P$-crystal structure of 
$\BB(\lambda)$ for an arbitrary level-zero 
dominant integral weight $\lambda \in P$.  

This paper is organized as follows. 
In Section~\ref{sec:pre}, 
we fix our notation, and 
recall some basic facts concerning 
affine Lie algebras and crystals of all LS paths.
In addition, 
we briefly review our results about
the crystal structure of 
the $P_{\cl}$-crystal $\BB(\lambda)_{\cl}$ 
for a level-zero dominant integral weight 
$\lambda = \sum_{i \in I_{0}} m_{i} \vpi_{i}$, 
with $m_{i} \in \BZ_{\ge 0}$ for $i \in I_{0}$.
In Section~\ref{sec:con}, 
we give an explicit description (Theorem~\ref{thm:comps}) of 
the decomposition of 
the $P$-crystal $\BB(\lambda)$ into 
connected components for 
the $\lambda = \sum_{i \in I_{0}} m_{i} \vpi_{i}$ above.
In Section~\ref{sec:aff}, 
we study the decomposition of 
the affinization $\ha{\BB(\lambda)_{\cl}}$ of 
$\BB(\lambda)_{\cl}$ into 
connected components (Theorem~\ref{thm:isom}), and 
then obtain a condition (Corollary~\ref{cor:con}) 
for an element of $\ha{\BB(\lambda)_{\cl}}$ to 
lie in (the isomorphic image of) 
the connected component $\BB_{0}(\lambda)$ 
of $\BB(\lambda)$.
%
%
%
%
\section{Preliminaries.}
\label{sec:pre}

%
\subsection{Affine Lie algebras and quantum affine algebras.}
\label{subsec:affine}

Let $A=(a_{ij})_{i,j \in I}$ 
be a generalized Cartan matrix of affine type. 
Throughout this paper, we assume that 
the elements of the index set $I$ are numbered 
as in \cite[Section~4.8, Tables Aff~1\,--\,Aff~3]{Kac}. 
Take a special vertex $0 \in I$ as in these tables, 
and set $I_{0}:=I \setminus \{0\}$. 
Let $\Fg=\Fg(A)$ be the affine Lie algebra 
associated to the Cartan matrix 
$A=(a_{ij})_{i,j \in I}$ of affine type 
over the field $\BQ$ of rational numbers, and 
let $\Fh$ be its Cartan subalgebra. 
If we denote by 
$\Pi^{\vee}:=\bigl\{h_{j}\bigr\}_{j \in I} \subset \Fh$ 
the set of the simple coroots, and 
by $d \in \Fh$ the scaling element, then 
we have
$\Fh=\bigl(\bigoplus_{j \in I} \BQ h_{j}\bigr) \oplus \BQ d$.
Also, we denote by 
$\Pi:=\bigl\{\alpha_{j}\bigr\}_{j \in I} \subset 
\Fh^{\ast}:=\Hom_{\BQ}(\Fh,\BQ)$ 
the set of the simple roots, and by 
$\Lambda_{j} \in \Fh^{\ast}$, $j \in I$, 
the fundamental weights; 
note that $\alpha_{j}(d)=\delta_{j,0}$ and 
$\Lambda_{j}(d)=0$ for all $j \in I$ 
(see \cite[Section~6.1]{Kac} and \cite[Section~10.1]{HK}). 
Let 
\begin{equation}
\delta=\sum_{j \in I} a_{j}\alpha_{j} \in \Fh^{\ast}
\qquad \text{ and } \qquad 
c=\sum_{j \in I} a^{\vee}_{j} h_{j} \in \Fh
\end{equation}
be the null root and the canonical central element 
of $\Fg$, respectively. 
Here we should note that 
\begin{equation}
a_{0}=
 \begin{cases}
 2 & \text{if $\Fg$ is of type $A^{(2)}_{2\ell}$}, 
 \\[1mm]
 1 & \text{otherwise}.
 \end{cases}
\end{equation}
We define the Weyl group $W$ of $\Fg$ by: 
$W=\langle r_{j} \mid j \in I\rangle 
\subset \GL(\Fh^{\ast})$, where 
$r_{j} \in \GL(\Fh^{\ast})$ is the simple reflection 
associated to $\alpha_{j}$ for $j \in I$, and then 
define the set $\rr$ of real roots by: $\rr=W\Pi$. 
The set of positive real roots is denoted by 
$\prr \subset \rr$. Also, let us denote by 
$(\cdot\,,\,\cdot)$ the (standard) 
bilinear form on $\Fh^{\ast}$ 
normalized as in \cite[Section~6.2]{Kac};
we have 
$(\alpha_{j},\alpha_{j})=
2a_{j}^{\vee}a_{j}^{-1}$ for all $j \in I$.

We take a dual weight lattice $P^{\vee}$ 
and a weight lattice $P$ as follows:
%
%
\begin{equation} \label{eq:lattices}
P^{\vee}=
\left(\bigoplus_{j \in I} \BZ h_{j}\right) \oplus \BZ d \, 
\subset \Fh
\quad \text{and} \quad 
P= 
\left(\bigoplus_{j \in I} \BZ \Lambda_{j}\right) \oplus 
   \BZ a_{0}^{-1}\delta
   \subset \Fh^{\ast}.
\end{equation}
It is easily seen that 
$P \cong \Hom_{\BZ}(P^{\vee},\BZ)$. 
The quintuplet
$(A,P,P^{\vee},\Pi,\Pi^{\vee})$ is called 
a Cartan datum for the generalized Cartan matrix 
$A=(a_{ij})_{i,j \in I}$ of affine type 
(see \cite[Definition~2.1]{HK}).

Let 
$\cl:\Fh^{\ast} \twoheadrightarrow \Fh^{\ast}/\BQ\delta$ be 
the canonical projection from $\Fh^{\ast}$ onto 
$\Fh^{\ast}/\BQ\delta$, and define the classical weight lattice 
$P_{\cl}$ and the dual weight lattice $P_{\cl}^{\vee}$ by:
\begin{equation}
P_{\cl} = \cl(P) = 
 \bigoplus_{j \in I} \BZ \cl(\Lambda_{j})
\quad \text{and} \quad 
P_{\cl}^{\vee} = 
 \bigoplus_{j \in I} \BZ h_{j} 
 \subset P^{\vee}.
\end{equation}
Note that $P_{\cl} \simeq P/(\BQ\delta \cap P)$, and that 
$P_{\cl}$ can be identified with 
$\Hom_{\BZ}(P_{\cl}^{\vee},\BZ)$ as a $\BZ$-module by: 
$(\cl(\lambda))(h)=\lambda(h)$ 
for $\lambda \in P$ and $h \in P_{\cl}^{\vee}$.  
The quintuple $(A, \cl(\Pi), \Pi^{\vee}, 
P_{\cl}, P_{\cl}^{\vee})$ is 
called a classical Cartan datum 
(see \cite[Section~10.1]{HK}). 

We set 
\begin{align}
& 
\fr:=\rr \, \cap \, 
     \sum_{j \in I_{0}}\BZ \alpha_{j}, 
\qquad 
\pfr:=\prr \, \cap \, 
      \sum_{j \in I_{0}}\BZ \alpha_{j}, \\[1.5mm]
& 
\fin{\Pi}:=
  \bigl\{\alpha_{j}\bigr\}_{j \in I_{0}} \subset \Pi, 
\qquad
\fin{W}:=\langle r_{j} \mid j \in I_{0}\rangle \subset W.
\end{align}
Note that $\fr=\fin{W}\fin{\Pi}$. In this paper, 
we call a root $\beta \in \fr$ a long (resp., short) 
root in $\fr$ if $(\beta,\beta) \ge (\gamma,\gamma)$ 
for all $\gamma \in \fr$ 
(resp., otherwise).
For each $\beta \in \fr$, we set 
%
%
\begin{equation} \label{eq:cbeta}
c_{\beta}:=\max\left\{1, \, 
 \frac{(\beta,\beta)}{2}\right\}.
\end{equation}
Note that $c_{w\beta}=c_{\beta}$ 
for all $\beta \in \fr$ and 
$w \in \fin{W}$, since 
the bilinear form $(\cdot\,,\,\cdot)$ on 
$\Fh^{\ast}$ is $W$-invariant. 
%
%
\begin{rem} \label{rem:cbeta}
By direct computation, we obtain the following:

\vsp

\noindent 
(1) If $\Fg$ is nontwisted, then 
$(\alpha_{j},\alpha_{j})/2=1$ 
for every long simple root 
$\alpha_{j} \in \fin{\Pi}$, 
and $(\alpha_{j},\alpha_{j})/2$ 
is equal to $1/3$, 
$1/2$, or $1$ for each $j \in I_{0}$. 
Therefore, we have 
$c_{\alpha_{j}}=1$ for all $j \in I_{0}$, 
and hence 
$c_{\beta}=1$ for all $\beta \in \fr$.

\vsp

\noindent 
(2) If $\Fg$ is twisted, then 
$(\alpha_{j},\alpha_{j})/2=1$ 
for every short simple root 
$\alpha_{j} \in \fin{\Pi}$, 
and $(\alpha_{j},\alpha_{j})/2$ is 
equal to $1$, $2$, or $3$ 
for each $j \in I_{0}$. 
Therefore, we have $c_{\alpha_{j}}=
(\alpha_{j},\alpha_{j})/2$ for all $j \in I_{0}$, 
and hence for each $\beta \in \fr$, 
\begin{equation*}
c_{\beta}=
 \frac{(\beta,\beta)}{2}=
\begin{cases}
1 & \text{if $\beta$ is a short root in $\fr$}, \\[1mm]
\text{$2$ or $3$} 
  & \text{if $\beta$ is a long root in $\fr$}. 
\end{cases}
\end{equation*}
\end{rem}

For $j \in I_{0}$, we set 
\begin{equation}
\ti{\alpha}_{j}:=
\begin{cases}
\alpha_{j}
 & \text{if $\Fg$ is twisted, 
             and not of type $A^{(2)}_{2\ell}$},
\\[3mm]
\dfrac{2\alpha_{j}}{(\alpha_{j},\alpha_{j})}
 & \text{otherwise}, 
\end{cases}
\end{equation}
and $\Gamma:=\sum_{j \in I_{0}} 
\BZ \ti{\alpha}_{j}$ (see \cite[(6.5.8)]{Kac}). 
For $\beta \in \Gamma$, we denote by 
$t_{\beta} \in \GL(\Fh^{\ast})$ 
the translation of $\Fh^{\ast}$ with respect to $\beta$ 
(see \cite[(6.5.2)]{Kac}). 
We know from \cite[(6.5.6)]{Kac} that 
$T:=\bigl\{t_{\beta} \mid \beta \in \Gamma\bigr\}$ is 
an abelian group, and from \cite[Proposition~6.5]{Kac} 
that the Weyl group $W$ of $\Fg$ decomposes into 
the semidirect product $\fin{W} \ltimes T$ 
of the Weyl group $\fin{W}$ (of finite type) and 
the group $T$ of translations. 
%
%
%
%
\subsection{Real roots of an affine Lie algebra and 
their finite parts.}
\label{subsec:fp}
The following proposition immediately follows 
from \cite[Proposition~6.3]{Kac} and 
Remark~\ref{rem:cbeta}.
%
%
\begin{prop} \label{prop:real}
{\rm (1) }
If $\Fg$ is not of type 
$A^{(2)}_{2\ell}$, then 
$\prr=\bigl\{
 \beta+nc_{\beta}\delta \mid 
 \beta \in \pfr,\, n \in \BZ_{\ge 0}
 \bigr\} \cup \bigl\{
 -\beta+nc_{\beta}\delta \mid 
 \beta \in \pfr,\, n \in \BZ_{\ge 1}
 \bigr\}$.

\vsp

\noindent {\rm (2) }
If $\Fg$ is of type $A^{(2)}_{2\ell}$, 
then 
\begin{equation*}
\begin{array}{ll}
\prr = & 
 \bigl\{
 \beta+nc_{\beta}\delta \mid 
 \beta \in \pfr,\, n \in \BZ_{\ge 0} \bigr\} \cup 
 \bigl\{
 -\beta+nc_{\beta}\delta \mid 
 \beta \in \pfr,\, n \in \BZ_{\ge 1} \bigr\} \, 
 \cup \\[1.5mm]
&
 \bigl\{
 \tfrac{1}{2}\bigl(\beta+(2n-1)\delta\bigr) \mid 
 \text{\rm $\beta$ is a long root in $\fr$}, \ 
 n \in \BZ_{\ge 1}
 \bigr\}.
\end{array}
\end{equation*}
\end{prop}
%
%
\begin{dfn} \label{dfn:fp-rr}
For $\xi \in \prr$, we define 
the finite part $\ba{\xi} \in \pfr$ of 
$\xi$ as follows: 

\vsp

\noindent 
(1) If $\xi$ is of the form 
$\xi=\beta+nc_{\beta}\delta$ with
$\beta \in \fr$ and $n \in \BZ_{\ge 0}$, then 
we set $\ba{\xi}:=\sgn(\beta)\beta$, where 
$\sgn(\beta):=1$ if $\beta \in \pfr$, and 
$\sgn(\beta):=-1$ if $\beta \in -\pfr$. 

\vsp

\noindent 
(2) If $\xi$ is of the form 
$\xi=\frac{1}{2}(\beta+(2n-1)\delta)$ with
$\beta \in \fr$ and $n \in \BZ_{\ge 1}$, then 
we set $\ba{\xi}:=\sgn(\beta)\beta$. 
\end{dfn}
%
%
%
%
\subsection{Level-zero integral weights.}
\label{subsec:lvzero}
%
%
\begin{dfn} \label{dfn:lvzero}
An integral weight $\lambda \in P$ is said to 
be level zero if $\lambda(c)=0$. 
An integral weight $\lambda \in P$ of 
level zero is said to be level-zero dominant 
if $\lambda(h_{j}) \ge 0$ for all $j \in I_{0}$. 
\end{dfn}
For each $i \in I_{0}$, 
we define a level-zero fundamental weight 
$\vpi_{i} \in P$ by: 
$\vpi_{i}=\Lambda_{i}-a_{i}^{\vee}\Lambda_{0}$, 
and $d_{i} \in \BZ_{> 0}$ by: 
%
%
\begin{equation} \label{eq:di}
d_{i}=
\begin{cases}
1 & 
 \text{if $\Fg$ is of 
 type $A^{(2)}_{2\ell}$ and $i=\ell$}, 
\\[1.5mm]
c_{\alpha_{i}} & \text{otherwise}. 
\end{cases}
\end{equation}

%
\begin{lem} \label{lem:dlam}
Let $\lambda=\sum_{i \in I_{0}} m_{i}\vpi_{i}$, 
with $m_{i} \in \BZ_{\ge 0}$ for $i \in I_{0}$, be 
a level-zero dominant integral weight. 
Then, 
$T\lambda=\lambda+
\left( \sum_{i \in I_{0}} m_{i}d_{i}\BZ \right) \delta$.
\end{lem}

\begin{proof}
We deduce from Remark~\ref{rem:cbeta} and 
the definition \eqref{eq:di} of $d_{i}$ that 
$(\vpi_{i},\ti{\alpha}_{j})=\delta_{ij}d_{i}$ 
for all $i,\,j \in I_{0}$. 
Since $\lambda$ is level-zero, 
it follows from \cite[(6.5.5)]{Kac} that 
$t_{\beta} \lambda = 
\lambda-(\lambda,\beta)\delta$
for all $\beta \in \Gamma$. 
Therefore, we have
\begin{align*}
T\lambda
 & =\lambda + (\lambda,\Gamma)\delta
   =\lambda + \left(\sum_{j \in I_{0}}
    (\lambda, \ti{\alpha}_{j}) \BZ \right)
    \delta
\\[3mm]
&  =\lambda + \left(\sum_{i,j \in I_{0}}
    m_{i}(\vpi_{i}, \ti{\alpha}_{j}) \BZ \right)
    \delta
   =\lambda-\left(\sum_{i \in I_{0}}
    m_{i}d_{i}\BZ \right)
    \delta.
\end{align*}
This proves the lemma. 
\end{proof}
%
%
%
%
\begin{lem} \label{lem:Worbit}
Let $\lambda=\sum_{i \in I_{0}}m_{i}\vpi_{i}$, 
with $m_{i} \in \BZ_{\ge 0}$ for $i \in I_{0}$, 
be a level-zero dominant integral weight. 
Then we have $W\lambda \subset \lambda-\fin{Q}_{+}+\BZ\delta$, 
where $\fin{Q}_{+}:=\sum_{j \in I_{0}}\BZ_{\ge 0}\alpha_{j}$. 
\end{lem}

\begin{proof}
We have 
\begin{align*}
W\lambda & = \fin{W}T\lambda 
 \subset \fin{W}(\lambda+\BZ\delta)
\quad \text{by Lemma~\ref{lem:dlam}}
\\[1.5mm]
& = \fin{W}\lambda+\BZ\delta 
  \subset \lambda-\fin{Q}_{+}+\BZ\delta 
\quad \text{since $\lambda$ is level-zero dominant}.
\end{align*}
This proves the lemma. 
\end{proof}
For a real root $\beta \in \rr$, 
let $\beta^{\vee} \in \Fh$ denote its dual root. 
The proof of the following lemma is straightforward. 
%
%
%
%
\begin{lem} \label{lem:pair0}
{\rm (1) }
Let $\xi \in \prr$ be a positive real root of 
the form $\xi=\beta+nc_{\beta}\delta$ with 
$\beta \in \fr$ and $n \in \BZ$, and 
let $\nu \in P$ be 
a level-zero integral weight. 
Then, $\nu(\xi^{\vee})=\nu(\beta^{\vee})$, and 
$r_{\xi}(\nu)=
 \nu-\nu(\beta^{\vee})\beta-
 \nu(\beta^{\vee})nc_{\beta}\delta$.

\noindent {\rm (2) }
Let $\xi \in \prr$ be a positive real root of 
the form $\xi=\frac{1}{2}(\beta+(2n-1)\delta)$ 
with $\beta \in \fr$ and $n \in \BZ$, and 
let $\nu \in P$ be 
a level-zero integral weight. 
Then, $\nu(\xi^{\vee})=
2\nu(\beta^{\vee})$, and 
$r_{\xi}(\nu)=
 \nu-\nu(\beta^{\vee})\beta-
 \nu(\beta^{\vee})(2n-1)\delta$.
\end{lem}
%
%
%
\subsection{Definition of $\sigma$-chains.}
\label{subsec:achain}
%
%
%
%
\begin{dfn} \label{dfn:bruhat}
Let $\mu,\,\nu \in P$ be integral weights. 
We write $\mu > \nu$ 
if there exist a sequence 
$\mu=\nu_{0},\,\nu_{1},\,
\dots,\,\nu_{k}=\nu$ of integral weights and 
a sequence $\xi_{1},\,\dots,\,\xi_{k}$ of 
positive real roots such that
$\nu_{l}=r_{\xi_{l}}(\nu_{l-1})$ and 
$\nu_{l-1}(\xi_{l}^{\vee}) \in \BZ_{< 0}$ 
for $1 \le l \le k$ 
(we write $\mu \ge \nu$ if 
$\mu > \nu$ or $\mu = \nu$). 
In this case, the sequence 
$\nu_{0},\,\nu_{1},\,\dots,\,\nu_{k}$ above is 
called a chain for $(\mu,\nu)$. 
If $\mu > \nu$, then 
we define $\dist(\mu,\nu)$ to 
be the maximum length $k$ of 
all possible chains for $(\mu,\nu)$. 
\end{dfn}
%
%
\begin{rem} \label{rem:dist=1}
Let $\mu,\,\nu \in P$ be integral weights such that $\mu > \nu$. 
Let $\mu=\nu_{0},\,\nu_{1},\,
\dots,\,\nu_{k}=\nu$ be a chain for $(\mu,\nu)$, and 
assume that both $\xi_{1},\,\dots,\,\xi_{k}$ and 
$\zeta_{1},\,\dots,\,\zeta_{k}$ are sequences of 
positive real roots corresponding to the chain 
$\mu=\nu_{0},\,\nu_{1},\,\dots,\,\nu_{k}=\nu$.
Then it immediately follows from \cite[Lemma~2.1.4]{NSp2}
that $\xi_{l}=\zeta_{l}$ for all $1 \le l \le k$. 
Thus, a sequence of positive real roots corresponding to 
each chain is uniquely determined. 
\end{rem}
%
%
\begin{rem} \label{rem:Q+}
(1) Let $\mu,\,\nu \in P$ be integral weights such that
$\mu > \nu$. 
Assume that $\mu=\nu_{0},\,\nu_{1},\,
\dots,\,\nu_{k}=\nu$ is a chain for $(\mu,\nu)$, 
with $\xi_{1},\,\dots,\,\xi_{k}$ 
the corresponding positive real roots.
Then it follows that 
$\nu-\mu \in \sum_{l=1}^{k}\BZ_{> 0} \xi_{l}$. 
In particular, $\nu-\mu \in 
Q_{+} \setminus \{0\}$, 
where $Q_{+}:=\sum_{j \in I} \BZ_{\ge 0}\alpha_{j}$. 

\vsp

\noindent
(2) Let $\mu \in P$, and assume 
that $\mu(h_{j}) \in \BZ_{< 0}$ for some $j \in I$. 
Then we have $\mu > r_{j}(\mu)$. 
Furthermore, 
it immediately follows from part (1) that 
$\dist(\mu,\,r_{j}(\mu))=1$. 
Indeed, assume that 
$\dist(\mu,\,r_{j}(\mu))=k$, and 
let $\xi_{1},\,\dots,\,\xi_{k}$ be 
the positive real roots corresponding to a chain 
of the maximum length. 
Then, by part (1), we have 
$-\mu(h_{j})\alpha_{j} = 
 r_{j}(\mu)-\mu 
 \in -\mu(\xi_{1}^{\vee})\xi_{1}+
 \sum_{l=2}^{k}\BZ_{> 0} \xi_{l}$. 
Therefore, 
from the linearly independence of the simple roots, 
we obtain that 
$\xi_{l}=\alpha_{j}$ for all $1 \le l \le k$. 
Hence we have $k=1$, 
since $\xi_{1}=\alpha_{j}$ implies $\xi_{1}^{\vee}=h_{j}$. 
\end{rem}
%
%
%
%
\begin{lem} \label{lem:dist=1}
Let $\mu$, $\nu \in P$ be level-zero integral weights. 
Assume that $\mu > \nu$ and $\dist(\mu,\nu)=1$, and let 
$\xi \in \prr$ be the corresponding positive real root. 
Then, $\xi$ is one of the following forms\,{\rm:}
$\xi=\beta$ with $\beta \in \pfr$, 
$\xi=-\beta+c_{\beta}\delta$ 
with $\beta \in \pfr$, or 
$\xi=\frac{1}{2}(-\beta+\delta)$ 
with $\beta \in \pfr$. 
\end{lem}

\begin{proof}
First, suppose that $\xi=\beta+nc_{\beta}\delta$ 
for some $\beta \in \pfr$ and $n \in \BZ_{\ge 1}$. 
Then we see from Proposition~\ref{prop:real} that
$\xi^{\prime}:=-\beta+nc_{\beta}\delta$ is 
a positive real root. 
By simple computations 
(using Lemma~\ref{lem:pair0}), 
we can easily verify that 
\begin{equation*}
\nu_{0}:=\mu, \quad 
\nu_{1}:=r_{\beta}(\mu), \quad 
\nu_{2}:=r_{\xi^{\prime}}r_{\beta}(\mu), \quad 
\nu_{3}:=r_{\beta}r_{\xi^{\prime}}r_{\beta}(\mu)=\nu
\end{equation*}
is a chain for $(\mu,\,\nu)$. Hence we have 
$\dist(\mu,\nu) \ge 3$, which is a contradiction. 

Now, suppose that $\xi=-\beta+nc_{\beta}\delta$ 
for some $\beta \in \pfr$ and $n \in \BZ_{\ge 2}$. 
Then we see from Proposition~\ref{prop:real} that
$\xi^{\prime}:=\beta+(n-2)c_{\beta}\delta$ and 
$\xi^{\prime\prime}:=-\beta+c_{\beta}\delta$ are 
positive real roots. 
By simple computations (using Lemma~\ref{lem:pair0}), 
we can easily verify that 
\begin{equation*}
\nu_{0}:=\mu, \quad 
\nu_{1}:=r_{\xi^{\prime\prime}}(\mu), \quad
\nu_{2}:=
  r_{\xi^{\prime}}
  r_{\xi^{\prime\prime}}(\mu), \quad 
  \nu_{3}:=
  r_{\xi^{\prime\prime}}
  r_{\xi^{\prime}}
  r_{\xi^{\prime\prime}}(\mu)=\nu
\end{equation*}
is a chain for $(\mu,\,\nu)$. Hence we have
$\dist(\mu,\nu) \ge 3$, which is a contradiction. 

Now, suppose that
$\xi=\frac{1}{2}\bigl(\beta+(2n-1)\delta\bigr)$ 
for some $\beta \in \pfr$ and $n \in \BZ_{\ge 1}$; 
note that $\Fg$ is of type $A_{2\ell}^{(2)}$. 
Then we see 
from Proposition~\ref{prop:real}\,(2) that
$\xi^{\prime}:=
\frac{1}{2}\bigl(-\beta+(2n-1)\delta\bigr)$ is 
a positive real root. 
By simple computations (using Lemma~\ref{lem:pair0}), 
we can easily verify that 
\begin{equation*}
\nu_{0}:=\mu, \quad 
\nu_{1}:=r_{\beta}(\mu), \quad 
\nu_{2}:=r_{\xi^{\prime}}r_{\beta}(\mu), \quad 
\nu_{3}:=r_{\beta}r_{\xi^{\prime}}r_{\beta}(\mu)=\nu
\end{equation*}
is a chain for $(\mu,\,\nu)$. Hence we have
$\dist(\mu,\nu) \ge 3$, which is a contradiction. 

Finally, suppose that
$\xi=\frac{1}{2}\bigl(-\beta+(2n-1)\delta\bigr)$ 
for some $\beta \in \pfr$ and $n \in \BZ_{\ge 2}$; 
note that $\Fg$ is of type $A_{2\ell}^{(2)}$.
Then we see from Proposition~\ref{prop:real}\,(2) 
that $\xi^{\prime}:=
\frac{1}{2}\bigl(\beta+(2n-3)\delta\bigr)$ and 
$\xi^{\prime\prime}:=
\frac{1}{2}\bigl(-\beta+\delta\bigr)$ are
positive real roots. 
By simple computations (using Lemma~\ref{lem:pair0}), 
we can easily verify that
\begin{equation*}
\nu_{0}:=\mu, \quad 
\nu_{1}:=r_{\xi^{\prime\prime}}(\mu), \quad 
\nu_{2}:=r_{\xi^{\prime}}r_{\xi^{\prime\prime}}(\mu), \quad 
\nu_{3}:=r_{\xi^{\prime\prime}} r_{\xi^{\prime}} 
 r_{\xi^{\prime\prime}}(\mu)=\nu
\end{equation*}
is a chain for $(\mu,\,\nu)$. Hence we have
$\dist(\mu,\nu) \ge 3$, which is a contradiction. 
This proves the lemma. 
\end{proof}
%
%
%
%
\begin{dfn} \label{dfn:achain}
Let $\mu,\,\nu \in P$ be 
integral weights such that $\mu \ge \nu$, 
and let $0 < \sigma < 1$ be a rational number. 
A $\sigma$-chain for $(\mu,\nu)$ is, by definition, 
a sequence $\mu=\nu_{0} > \nu_{1} > 
\dots > \nu_{k}=\nu$ of integral weights 
satisfying either of the following conditions:

\vsp

\noindent
(1) $k=0$ and $\mu=\nu_{0}=\nu$, or 

\vsp

\noindent
(2) $k \ge 1$, and 
$\dist(\nu_{l-1},\nu_{l})=1$ and 
$\nu_{l-1}(\xi_{l}^{\vee}) \in \sigma^{-1}\BZ_{< 0}$ 
for every $1 \le l \le k$, 
where $\xi_{l}$ is the positive real root 
corresponding to the chain for 
$(\nu_{l-1},\nu_{l})$, $1 \le l \le k$. 
\end{dfn}
Using the fact that $\delta(\beta^{\vee})=0$ 
for all $\beta \in \prr$, 
we can easily derive the following lemma.
%
%
\begin{lem} \label{lem:d-shift}
Let $\mu,\,\nu \in P$ be integral weights 
such that $\mu \ge \nu$, and 
let $0 < \sigma < 1$ be a rational number. If there exists 
a $\sigma$-chain for $(\mu,\nu)$, then there also exists 
a $\sigma$-chain for $(\mu+n\delta,\ \nu+n\delta)$ 
for each $n \in a_{0}^{-1}\BZ$.
\end{lem}
%
%
%
\subsection{Path crystals with weight lattice $P$, or $P_{\cl}$.}
\label{subsec:path}
A path with weight in $P$ is, by definition, 
a piecewise linear, continuous map 
$\pi:[0,1] \rightarrow \BQ \otimes_{\BZ} P$ 
from $[0,1]:=\bigl\{t \in \BQ \mid 0 \le t \le 1\bigr\}$ 
to $\BQ \otimes_{\BZ} P$ such that 
$\pi(0)=0$ and $\pi(1) \in P$. 
We denote by $\BP$ the set of all paths 
$\pi:[0,1] \rightarrow \BQ \otimes_{\BZ} P$ 
with weight in $P$.
For two paths $\pi_{1},\,\pi_{2} \in \BP$, we define 
$\pi_{1} \pm \pi_{2} \in \BP$ by: 
$(\pi_{1} \pm \pi_{2})(t)=
\pi_{1}(t) \pm \pi_{2}(t)$ for $t \in [0,1]$.
%
%
%
\begin{dfn} \label{dfn:expr}
Let $\pi \in \BP$ be a path with weight in $P$. 
A pair $(\ud{\nu}\,;\,\ud{\sigma})$ of a sequence 
$\ud{\nu}:\nu_{1},\,\nu_{2},\,\dots,\,\nu_{s}$ of 
integral weights and a sequence 
$\ud{\sigma}:
 0=\sigma_{0} < \sigma_{1} < \cdots < \sigma_{s}=1$ of 
rational numbers is called 
an expression of $\pi$ if the following equation
%
%
\begin{equation} \label{eq:path}
\pi(t)=\sum_{u^{\prime}=1}^{u-1}
(\sigma_{u^{\prime}}-\sigma_{u^{\prime}-1})\nu_{u^{\prime}}+
(t-\sigma_{u-1})\nu_{u} 
\quad
\text{for all \, }
   \sigma_{u-1} \le t \le \sigma_{u}, \ 
   1 \le u \le s, 
\end{equation}
holds. If $(\ud{\nu}\,;\,\ud{\sigma})$ is an expression of 
a path $\pi \in \BP$, then 
we write $\pi=(\ud{\nu}\,;\,\ud{\sigma})$. 
\end{dfn}
%
%
\begin{rem} \label{rem:expr}
An expression for a path $\pi \in \BP$ is uniquely determined, 
up to the following two operations: 

\vsp 

\noindent (1) 
Let $(\nu_{1},\,\nu_{2},\,\dots,\,\nu_{s}\,;\,
\sigma_{0},\,\sigma_{1},\,\dots,\,\sigma_{s})$ be 
an expression of a path $\pi \in \BP$.
If $\nu_{u}=\nu_{u+1}$ for some $0 \le u \le s-1$, 
then we can ``omit'' $\nu_{u}$ and $\sigma_{u}$ from 
this expression of $\pi$. 
For example, if $\pi=
(\nu_{1},\,\nu_{2},\,\nu_{3}\,;\,
\sigma_{0},\,\sigma_{1},\,\sigma_{2},\,\sigma_{3})$ and 
$\nu_{2}=\nu_{3}$, 
then $(\nu_{1},\,\nu_{3}\,;\,
\sigma_{0},\,\sigma_{1},\,\sigma_{3})$ is also 
an expression of $\pi$, i.e., 
$\pi=(\nu_{1},\,\nu_{3}\,;\,
\sigma_{0},\,\sigma_{1},\,\sigma_{3})$. 

\vsp

\noindent (2) 
Let $(\nu_{1},\,\nu_{2},\,\dots,\,\nu_{s}\,;\,
\sigma_{0},\,\sigma_{1},\,\dots,\,\sigma_{s})$ be 
an expression of a path $\pi \in \BP$.
If $\sigma_{u} < \tau < \sigma_{u+1}$ 
for some $0 \le u \le s-1$, 
then we can ``insert'' $\tau$ 
into this expression of $\pi$. 
For example, if $\pi=
(\nu_{1},\,\nu_{2}\,;\,\sigma_{0},\,\sigma_{1},\,\sigma_{2})$ and 
$\sigma_{0} < \tau < \sigma_{1}$, then 
$(\nu_{1},\,\nu_{1},\,\nu_{2}\,;\,
 \sigma_{0},\,\tau,\,\sigma_{1},\,\sigma_{2})$ is 
also an expression of $\pi$, i.e., 
$\pi=(\nu_{1},\,\nu_{1},\,\nu_{2}\,;\,
 \sigma_{0},\,\tau,\,\sigma_{1},\,\sigma_{2})$. 
\end{rem}

For an integral weight $\nu \in P$, 
let us denote by $\pi_{\nu}$ 
the straight line connecting $0 \in P$ with $\nu \in P$, i.e., 
$\pi_{\nu}(t):=t\nu$ for $t \in [0,1]$. Note that 
$\pi_{\nu}$ has an expression of the form:
%
%
\begin{equation} \label{eq:straight}
\pi_{\nu}=
(\underbrace{\nu,\,\nu,\,\dots,\,\nu}_{
 \text{$s$ times}} \,;\, 
 \sigma_{0},\,\sigma_{1},\,\dots,\,\sigma_{s})
\end{equation}
for each 
$0=\sigma_{0} < \sigma_{1} < \cdots < \sigma_{s}=1$.

For each $j \in I$, we denote by 
$e_{j},\,f_{j} : \BP \cup \{\bzero\} 
\rightarrow \BP \cup \{\bzero\}$ 
the root operators corresponding to 
the simple root $\alpha_{j}$, 
which are introduced in \cite[Section~1]{L2} 
(see also \cite[Section~1.2]{NSp2} for 
an explicit description of $e_{j}$ and $f_{j}$). 
Here the $\bzero$ is an additional element 
corresponding to the $0$ 
in the theory of crystals. 
For each $j \in I$, we define 
$S_{j}:\BP \rightarrow \BP$ by:
%
%
\begin{equation} \label{eq:sj}
S_{j}\pi=
\begin{cases}
f_{j}^{n}\pi 
  & \text{if \,} n=(\pi(1))(h_{j}) \ge 0, \\[3mm]
e_{j}^{-n}\pi 
  & \text{if \,} n=(\pi(1))(h_{j}) < 0,
\end{cases}
\end{equation}
for $\pi \in \BP$.
We know from \cite[Theorem 8.1]{L2} that 
there exists a unique action 
$S:W \rightarrow \Bij(\BP)$, 
$w \mapsto S_{w}$, of the Weyl group $W$ 
on the set $\BP$ of all paths such that 
$S_{r_{j}}=S_{j}$ for all $j \in I$, 
where $\Bij(\BP)$ denotes
the group of all bijections 
from the set $\BP$ to itself; 
if $w \in W$ equals 
$r_{i_{1}}r_{i_{2}} \cdots r_{j_{l}}$ 
for some $j_{1},\,j_{2},\,\dots,\,j_{l} \in I$, then 
$S_{w}=S_{j_{1}}S_{j_{2}} \cdots S_{j_{l}}$. 
Note that we have $(S_{w}\pi)(1)=w(\pi(1))$ 
for all $w \in W$ and $\pi \in \BP$. 
Moreover, 
we know from \cite[Lemma~5.2]{GL} that 
%
%
\begin{equation} \label{eq:weyl}
S_{w}\pi_{\nu}=\pi_{w\nu} \quad
\text{for $\nu \in P$ and $w \in W$}. 
\end{equation}
\begin{dfn}
A path $\pi \in \BP$ is said to be extremal 
if for every $w \in W$, either 
$e_{j}S_{w}\pi=\bzero$ or 
$f_{j}S_{w}\pi=\bzero$ holds for each $j \in I$. 
\end{dfn}
%
%
\begin{rem} \label{rem:ext}
We can easily deduce from \eqref{eq:weyl} and 
the definition of an extremal path that 
for every $\nu \in P$, the straight line 
$\pi_{\nu} \in \BP$ is an extremal path. 
\end{rem}

A path with weight in $P_{\cl}$ is, by definition, 
a piecewise linear, continuous map 
$\eta:[0,1] \rightarrow \BQ \otimes_{\BZ} P_{\cl}$ 
from $[0,1]=\bigl\{t \in \BQ \mid 0 \le t \le 1\bigr\}$ 
to $\BQ \otimes_{\BZ} P_{\cl}$ such that 
$\eta(0)=0$ and $\eta(1) \in P_{\cl}$. 
We denote by $\BP_{\cl}$ the set of all paths 
$\eta:[0,1] \rightarrow \BQ \otimes_{\BZ} P_{\cl}$ 
with weight in $P_{\cl}$. 
For $\pi \in \BP$, we define a path 
$\cl(\pi) \in \BP_{\cl}$ by: 
$(\cl(\pi))(t)=\cl(\pi(t))$ 
for $t \in [0,1]$.
%
%
\begin{rem} \label{rem:exp-cl}
Let $\pi$ be a path with weight in $P$, and $\pi=
(\nu_{1},\,\nu_{2},\,\dots,\,\nu_{s} \,;\,
\sigma_{0},\,\sigma_{1},\,\dots,\,\sigma_{s})$
an expression of $\pi$.
Then the path $\cl(\pi)$ with weight in $P_{\cl}$ 
has an expression of the form: 
$\cl(\pi)=
(\cl(\nu_{1}),\,\cl(\nu_{2}),\,\dots,\,\cl(\nu_{s}) \,;\,
\sigma_{0},\,\sigma_{1},\,\dots,\,\sigma_{s})$.
Here an expression of a path with weight in $P_{\cl}$ is 
defined similarly to that of a path with weight in $P$ 
(see Definition~\ref{dfn:expr}).
\end{rem}
Root operators on $\BP_{\cl} \cup \{\bzero\}$ 
are defined similarly to those on $\BP \cup \{\bzero\}$; 
we denote them also by $e_{j}$ and $f_{j}$ 
for each $j \in I$. 
We know from \cite[Lemma~2.7\,(1)]{NSz} that 
%
%
\begin{equation} \label{eq:cl-ro}
\cl(e_{j}\pi)=e_{j}\cl(\pi), \quad 
\cl(f_{j}\pi)=f_{j}\cl(\pi) \quad 
\text{for \,} \pi \in \BP 
\text{\, and \,} j \in I,
\end{equation}
where $\cl(\bzero)$ is understood to be $\bzero$. 
Note that we have a natural action of the Weyl group $W$ 
on $\BQ \otimes_{\BZ} P_{\cl} \simeq \Fh^{\ast}/\BQ\delta$, 
since $W\delta=\delta$.
In the same manner as we defined 
$S_{w}$, $w \in W$, on $\BP$, this time 
using root operators on $\BP_{\cl}$, 
we can define an action, 
also denoted by $S_{w}$, $w \in W$, 
of the Weyl group $W$ 
on the set $\BP_{\cl}$ such that $(S_{w}\eta)(1)=w(\eta(1))$ 
for every $w \in W$ and $\eta \in \BP_{\cl}$. 
In fact, we have 
%
%
\begin{equation} \label{eq:cl-sw}
\cl(S_{w}\pi)=S_{w}\cl(\pi) \quad 
\text{for every $w \in W$ and $\pi \in \BP$}.
\end{equation}
A path $\eta \in \BP_{\cl}$ is 
said to be extremal if for every $w \in W$, 
either $e_{j}S_{w}\eta=\bzero$ or 
$f_{j}S_{w}\eta=\bzero$ holds for each $j \in I$. 
%
%
\begin{rem} \label{rem:cl-ext}
It immedaitely follows from \eqref{eq:cl-sw} that 
a path $\pi \in \BP$ is extremal if and only if 
$\cl(\pi) \in \BP_{\cl}$ is extremal. 
In particular, it follows from Remark~\ref{rem:ext} that 
the straight line $\cl(\pi_{\nu}) \in \BP_{\cl}$ is 
extremal for every $\nu \in P$.
\end{rem}
%
%
%
%
\subsection{Lakshmibai-Seshadri paths.}
\label{subsec:LS}
Let $\lambda \in P$ be an integral weight. 
Now we give the definition of 
a Lakshmibai-Seshadri path of shape $\lambda$.
%
%
%
%
\begin{dfn}[{cf. \cite[Section~4]{L2}}] \label{dfn:LS}
A Lakshmibai-Seshadri path (LS path for short) of 
shape $\lambda$ is a path $\pi \in \BP$ having an 
expression of the form: 
$(\nu_{1},\,\nu_{2},\,\dots,\,\nu_{s}\,;\,
\sigma_{0},\,\sigma_{1},\,\dots,\,\sigma_{s})$, where 
$\nu_{1},\,\nu_{2},\,\dots,\,
\nu_{s} \in W\lambda$, and where 
for each $1 \le u \le s-1$, 
there exists a $\sigma_{u}$-chain 
for $(\nu_{u},\,\nu_{u+1})$.
We denote by $\BB(\lambda)$ 
the set of all LS paths of shape $\lambda$. 
\end{dfn}
%
%
%
\begin{rem} \label{rem:LS}
(1) Let $\pi$ be 
an LS path of shape $\lambda$, and 
$(\nu_{1},\,\nu_{2},\,\dots,\,\nu_{s}\,;\,
\sigma_{0},\,\sigma_{1},\,\dots,\,\sigma_{s})$ 
an arbitrary expression of $\pi$. 
Then, we see from Remark~\ref{rem:expr} that 
$\nu_{1},\,\nu_{2},\,\dots,\,
\nu_{s} \in W\lambda$, and that 
for each $1 \le u \le s-1$, 
there exists an $\sigma_{u}$-chain 
for $(\nu_{u},\,\nu_{u+1})$.

\vsp

\noindent (2) It can easily be seen 
from part (1) and \eqref{eq:straight} that 
$\pi_{\nu} \in \BB(\lambda)$
if and only if $\nu \in W\lambda$.

\vsp

\noindent (3) 
It follows from the definitions that 
$\BB(w\lambda)=\BB(\lambda)$ for every $w \in W$. 
\end{rem}
%
%
%
%
\begin{thm}[{\cite[Sections~2 and 4]{L2}}] \label{thm:LSc}
Let $\lambda \in P$ be an integral weight. 

\vsp

\noindent {\rm (1) } 
We have $\pi(1) \in \lambda + Q$ for all $\pi \in \BB(\lambda)$, 
where $Q:=\sum_{j \in I}\BZ \alpha_{j}$. 

\vsp

\noindent {\rm (2) }
The set $\BB(\lambda) \cup \{\bzero\}$ is stable 
under the root operators $e_{j}$ and $f_{j}$ 
for all $j \in I$. 

\vsp

\noindent {\rm (3) }
We set 
%
%
\begin{equation} \label{eq:LS-c1}
\begin{cases}
\wt(\pi):=\pi(1) 
  & \text{\rm for \,} \pi \in \BB(\lambda), \\[1.5mm]
\ve_{j}(\pi):=\max\bigl\{l \in \BZ_{\ge 0} \mid 
  e_{j}^{l}\pi \ne \bzero \bigr\} 
  & \text{\rm for \,} \pi \in \BB(\lambda)
    \text{\rm \, and \,} j \in I, \\[1.5mm]
\vp_{j}(\pi):=\max\bigl\{l \in \BZ_{\ge 0} \mid 
  f_{j}^{l}\pi \ne \bzero \bigr\} 
  & \text{\rm for \,} \pi \in \BB(\lambda)
    \text{\rm \, and \,} j \in I.
\end{cases}
\end{equation}
Then the set $\BB(\lambda)$, together 
with the root operators $e_{j}$, $f_{j}$ for $j \in I$ 
and the maps $\wt$, $\ve_{j}$, $\vp_{j}$, $j \in I$, 
becomes a $P$-crystal, i.e.,
a crystal associated to 
the Cartan datum $(A,\Pi,\Pi^{\vee},P,P^{\vee})$ 
{\rm(}see {\rm \cite[Section~4.5 and Section~10.2]{HK}}{\rm)}. 
\end{thm}
%
%
%
%
\begin{lem} \label{lem:wt}
Let $\lambda=\sum_{i \in I_{0}}m_{i}\vpi_{i}$, 
with $m_{i} \in \BZ_{\ge 0}$ for $i \in I_{0}$, 
be a level-zero dominant integral weight, 
and let $\pi \in \BB(\lambda)$. Then, 
$\pi(1) \in \lambda-
a_{0}^{-1}\fin{Q}_{+}+a_{0}^{-1}\BZ\delta$. 
\end{lem}

\begin{proof}
It follows from 
Theorem~\ref{thm:LSc}\,(1) that 
$\pi(1) \in \lambda+Q$. 
We set $\theta := \delta-a_{0}\alpha_{0}$ 
(see \cite[Proposition~6.4]{Kac}).
Since $\theta \in \pfr \subset \fin{Q}_{+}$ and 
$\alpha_{0}=a_{0}^{-1}(\delta-\theta)$, we see that
$Q \subset a_{0}^{-1}\fin{Q}+a_{0}^{-1}\BZ\delta$, 
where $\fin{Q}:=\sum_{j \in I_{0}}\BZ \alpha_{j}$. 
Hence we have 
$\pi(1) \in \lambda+
a_{0}^{-1}\fin{Q}+a_{0}^{-1}\BZ\delta$.

Let $(\nu_{1},\,\nu_{2},\,\dots,\,\nu_{s} \, ; \, 
\sigma_{0},\,\sigma_{1},\,\dots,\,\sigma_{s})$ be an 
expression of the LS path 
$\pi \in \BB(\lambda)$ 
as in Definition~\ref{dfn:LS}. 
Then it follows from Definition~\ref{dfn:expr} that 
$\pi(1) = \sum_{u=1}^{s}(\sigma_{u}-\sigma_{u-1})\nu_{u}$, 
where $\nu_{1},\,\nu_{2},\,\dots,\,\nu_{s} \in W\lambda$ 
by the definition of an LS path.
Hence we see from Lemma~\ref{lem:Worbit} that 
$\pi(1) \in 
 \lambda - 
 \sum_{j \in I_{0}}\BQ_{\ge 0} \alpha_{j} + 
 \BQ\delta$.
Noting that 
$\alpha_{j}$, $j \in I_{0}$, and $\delta$ 
are linearly independent, we deduce from the above that 
$\pi(1) \in 
 \lambda-a_{0}^{-1}\fin{Q}_{+}+a_{0}^{-1}\BZ\delta$.
This proves the lemma. 
\end{proof}

Let us set
$\BB(\lambda)_{\cl}
  :=\cl(\BB(\lambda))
  =\bigl\{\cl(\pi) \mid \pi \in \BB(\lambda)\bigr\}
  \subset \BP_{\cl}$.
We know the following theorem 
from \cite[Theorem~2.4 and Section~3.1]{NSz}; 
%
%
\begin{thm} \label{thm:LScl}
Let $\lambda \in P$ be an integral weight. 

\vsp

\noindent {\rm (1) }
The set $\BB(\lambda)_{\cl} \cup \{\bzero\}$ is 
stable under the root operators 
$e_{j}$ and $f_{j}$ for all $j \in I$. 

\vsp 

\noindent {\rm (2) } 
We set $\wt(\eta):=\eta(1) \in P_{\cl}$ 
for $\eta \in \BB(\lambda)_{\cl}$, and 
define $\ve_{j}(\eta)$ and $\vp_{j}(\eta)$ 
for $\eta \in \BB(\lambda)_{\cl}$ and 
$j \in I$ as in \eqref{eq:LS-c1}.
Then, the set $\BB(\lambda)_{\cl}$, 
together with 
the root operators $e_{j}$, $f_{j}$ for $j \in I$ 
and the maps $\wt$, $\ve_{j}$, $\vp_{j}$ for $j \in I$, 
becomes a $P_{\cl}$-crystal, i.e., 
a crystal associated to the classical Cartan datum 
$(A,\cl(\Pi),\Pi^{\vee},P_{\cl},P_{\cl}^{\vee})$
{\rm(}see {\rm \cite[Section~10.2]{HK}}{\rm)}. 
\end{thm}
Here we summarize main results of \cite{NSz} 
in the following theorem.
%
%
\begin{thm}[{\cite[Section~3]{NSz}}] \label{thm:NSz}
Let $\lambda=\sum_{i \in I_{0}} m_{i}\vpi_{i}$, 
with $m_{i} \in \BZ_{\ge 0}$ for $i \in I_{0}$, 
be a level-zero dominant integral weight. 
Then\,{\rm:}

\vsp

\noindent
{\rm (1)} 
The $P_{\cl}$-crystal $\BB(\lambda)_{\cl}$ is 
connected. 

\vsp

\noindent
{\rm (2)}
There exists a unique
isomorphism of $P_{\cl}$-crystals 
$\BB(\lambda)_{\cl} \stackrel{\sim}{\rightarrow} 
\bigotimes_{i \in I_{0}} 
 \bigl(\BB(\vpi_{i})_{\cl}\bigr)^{\otimes m_{i}}$ 
that maps $\cl(\pi_{\lambda})$ to 
$\bigotimes_{i \in I_{0}}
 \bigl(\cl(\pi_{\vpi_{i}})\bigr)^{\otimes m_{i}}$. 

\vsp

\noindent 
{\rm (3)} 
The $P_{\cl}$-crystal $\BB(\lambda)_{\cl}$ is simple 
{\rm(}in the sense of {\rm \cite[Definition~3.17]{NSz}}{\rm)}.
Namely, the set of all extremal paths 
in $\BB(\lambda)_{\cl}$ coincides 
with the $W$-orbit of the straight line 
$\cl(\pi_{\lambda})$, and the set of 
all paths of weight $\cl(\lambda)$ in 
$\BB(\lambda)_{\cl}$ consists only of 
the straight line $\cl(\pi_{\lambda})$. 
\end{thm}

We proved in 
\cite[Propositions~3.4.1 and 3.4.2]{NSp2} that 
the $P_{\cl}$-crystal $\BB(\vpi_{i})_{\cl}$ is 
isomorphic as a $P_{\cl}$-crystal 
to the crystal base of 
the (finite-dimensional, irreducible) 
level-zero fundamental module $W(\vpi_{i})$, 
introduced in \cite[Section~5.2]{Kaslz}, 
over the quantized universal enveloping algebra 
$U_{q}^{\prime}(\Fg)$ of $\Fg$ 
(or rather, of the derived subalgebra $[\Fg,\Fg]$ of $\Fg$) 
over the field $\BQ(q)$ of rational functions 
with weight lattice $P_{\cl}$. 
%
%
%
\subsection{$\delta$-shifts of paths.}
\label{subsec:d-shift}
We know the following lemma 
from \cite[Lemma~3.3.1]{NSp2}.
%
%
\begin{lem} \label{lem:d-shift2}
Let $\pi,\,\pi^{\prime} \in \BP$, and let 
$F : [0,1] \rightarrow \BQ$ be a piecewise linear, 
continuous map such that 
$F(0)=0$ and $F(1) \in a_{0}^{-1}\BZ$. 
Assume that $\pi(t)=\pi^{\prime}(t)+F(t)\delta$ 
for all $t \in [0,1]$. Then\,{\rm:}

\vsp 

\noindent {\rm (1) } 
For each $j \in I$, 
$e_{j}\pi=\bzero$ if and only if 
$e_{j}\pi^{\prime}=\bzero$. Moreover, 
if $e_{j}\pi^{\prime} \ne \bzero$, then
$(e_{j}\pi)(t)=
(e_{j}\pi^{\prime})(t)+F(t)\delta$
for all $t \in [0,1]$.

\vsp 

\noindent {\rm (2) } 
For each $j \in I$, 
$f_{j}\pi=\bzero$ if and only 
if $f_{j}\pi^{\prime}=\bzero$. 
Moreover, if $f_{j}\pi^{\prime} \ne \bzero$, 
then 
$(f_{j}\pi)(t)=
(f_{j}\pi^{\prime})(t)+F(t)\delta$
for all $t \in [0,1]$.
\end{lem}

%
\begin{lem} \label{lem:d-shift3}
Let $\pi,\,\pi^{\prime} \in \BP$, and let 
$F : [0,1] \rightarrow \BQ$ be a piecewise linear, 
continuous map such that $F(0)=0$ and $F(1) \in a_{0}^{-1}\BZ$. 
Assume that $\pi(t)=\pi^{\prime}(t)+F(t)\delta$ 
for all $t \in [0,1]$. Then, for each $w \in W$, 
$(S_{w}\pi)(t)=(S_{w}\pi^{\prime})(t)+F(t)\delta$ 
for all $t \in [0,1]$. 
\end{lem}

\begin{proof}
It immediately follows from 
Lemma~\ref{lem:d-shift2} and \eqref{eq:sj} that 
for each $j \in I$, 
$(S_{j}\pi)(t)=(S_{j}\pi^{\prime})(t)+F(t)\delta$
for all $t \in [0,1]$.
Using this equation, 
we can easily prove the lemma 
by induction on 
the length $\ell(w)$ of $w \in W$.
\end{proof}
From Remark~\ref{rem:ext}, using 
Lemmas~\ref{lem:d-shift2} and \ref{lem:d-shift3}, 
we obtain the following lemma.
%
%
%
\begin{lem} \label{lem:d-shift-ext}
Let $\nu \in P$ be an integral weight, and let 
$F : [0,1] \rightarrow \BQ$ be 
a piecewise linear, continuous map 
such that $F(0)=0$ and $F(1) \in a_{0}^{-1}\BZ$. 
Then, the path $\pi_{\nu}(t)+F(t)\delta$, 
$t \in [0,1]$, is extremal. 
\end{lem}
The next lemma immediately follows from 
Definition~\ref{dfn:LS}, by using 
Lemma~\ref{lem:d-shift}.
%
%
\begin{lem} \label{lem:d-shift4}
Let $\lambda \in P$ be 
an integral weight, and let 
$n \in a_{0}^{-1}\BZ$. Then we have
$\BB(\lambda+n\delta)=
 \BB(\lambda)+\pi_{n\delta}$.
\end{lem}
%
%

\section{Connected components of 
 the $P$-crystal $\BB(\lambda)$.}
\label{sec:con}
%
%
%
\subsection{Main result.}
\label{subsec:thm}

Let $\lambda \in P$ be an integral weight 
of level zero. 
Then, $\lambda$ is 
equivalent, under the Weyl group 
$\fin{W} (\subset W)$ of finite type, to 
a unique level-zero dominant integral weight.
Since we are interested in 
the crystal structure of the set $\BB(\lambda)$ 
of all LS paths of shape $\lambda$, 
we may (and do) assume, 
in view of Remark~\ref{rem:LS}\,(3), 
that $\lambda \in P$ is a level-zero dominant 
integral weight. 
Hence, we can write the $\lambda \in P$ 
in the form $\sum_{i \in I_{0}} m_{i} \vpi_{i} + 
n \delta$, with $m_{i} \in \BZ_{\ge 0}$ 
for $i \in I_{0}$ and $n \in a_{0}^{-1}\BZ$. 
Furthermore, it follows from 
Lemmas~\ref{lem:d-shift2} and 
\ref{lem:d-shift4} that 
$\BB \bigl(
 \sum_{i \in I_{0}} m_{i} \vpi_{i} + n \delta
 \bigr)
=\BB\bigl(
 \sum_{i \in I_{0}} m_{i} \vpi_{i}
 \bigr) + \pi_{n\delta}$, 
and that 
the $P$-crystal 
$\BB \bigl(
 {\textstyle\sum_{i \in I_{0}} m_{i} \vpi_{i} + n \delta}
 \bigr)$ is ``isomorphic'' 
(up to a shift of weights by $n\delta$) 
to the $P$-crystal 
$\BB \bigl(
 {\textstyle\sum_{i \in I_{0}} m_{i} \vpi_{i}}
 \bigr)$. 
Thus, the problem of determining the $P$-crystal 
structure of $\BB(\lambda)$ 
for a level-zero integral weight 
$\lambda \in P$ is reduced to that for the case 
in which the $\lambda \in P$ is of the form: 
%
%
\begin{equation} \label{eq:lambda}
\lambda=\sum_{i \in I_{0}} m_{i} \vpi_{i}, \quad 
 \text{with $m_{i} \in \BZ_{\ge 0}$ for $i \in I_{0}$}.
\end{equation}
Therefore, in the remainder of this paper, 
we take (and fix) a level-zero dominant integral 
weight $\lambda$ of the form \eqref{eq:lambda}.

Now we set $\Supp_{\ge 2}(\lambda):=
\bigl\{i \in I_{0} \mid m_{i} \ge 2\bigr\}$, 
%
%
\begin{equation} \label{eq:turn}
\turn(\lambda):=\bigcup_{i \in \Supp_{\ge 2}(\lambda)} 
\bigl\{q/m_{i} \mid 1 \le q \le m_{i}-1\bigr\}, 
\end{equation}
and then denote by $s$ the number of elements of 
$\turn(\lambda)$ plus $1$.
When the set $\turn(\lambda)$ is nonempty, 
we enumerate the elements of $\turn(\lambda)$ in 
increasing order as follows: 
\begin{equation}
\turn(\lambda)=
\bigl\{
  \tau_{1} < \tau_{2} < \cdots < \tau_{s-1}
\bigr\}, 
\end{equation}
where $\tau_{u}$ denotes the $u$-th smallest 
element in $\turn(\lambda)$ 
for $1 \le u \le s-1$; 
by convention, we set $\tau_{0}:=0$ and $\tau_{s}:=1$. 
For each $1 \le u \le s-1$, 
we write $\tau_{u} \in \turn(\lambda)$ as 
$\tau_{u}=q_{u}/p_{u}$, 
where $p_{u},\,q_{u} \in \BZ_{> 0}$ 
are coprime, and then set $I_{0}(\lambda,p_{u}):=
\bigl\{i \in I_{0} \mid m_{i} \in p_{u}\BZ\bigr\}$. 

The aim in this section is to prove
the following theorem, which gives 
a parametrization of the connected components 
of the $P$-crystal $\BB(\lambda)$. 
%
%
\begin{thm} \label{thm:comps}
{\rm (1) } 
Let $\pi \in \BP$ be a path 
that has the following expression\,{\rm:}
%
%
\begin{equation} \label{eq:comps}
\begin{array}{l}
(\lambda-N_{1}\delta,\,\dots,\,
 \lambda-N_{s-1}\delta,\ \lambda \ ; \ 
 \tau_{0},\,\tau_{1},\,\dots,\,\tau_{s-1},\,\tau_{s}), 
   \\[3mm]
\hspace{10mm} \text{\rm with \, }
  N_{u}-N_{u+1} \in 
   {\displaystyle\sum_{j \in I_{0}(\lambda,p_{u})}}
   m_{j}d_{j} \BZ_{\ge 0}
  \quad \text{\rm for \, } 1 \le u \le s-1, 
\end{array}
\end{equation}
where $N_{s}:=0$. 
Then, $\pi \in \BP$ is 
an extremal LS path of shape $\lambda$. 

\vsp

\noindent {\rm (2) } 
Each connected component of the $P$-crystal 
$\BB(\lambda)$ contains exactly one 
{\rm(}extremal\,{\rm)} LS path 
having an expression of the form \eqref{eq:comps}.

\vsp

\noindent {\rm (3) }
The set of all extremal LS paths 
in each connected component of 
the crystal $\BB(\lambda)$ coincides 
with the $W$-orbit of an LS path 
having an expression of the form \eqref{eq:comps}.
\end{thm}

Let $\BB_{0}(\lambda)$ denote the connected component of 
$\BB(\lambda)$ containing the straight line $\pi_{\lambda}$; 
note that $\pi_{\lambda}$ has an expression of 
the form \eqref{eq:comps} for $N_{1}=\cdots=N_{s-1}=0$. 
We deduce from Theorem~\ref{thm:comps}, 
along with Lemma~\ref{lem:d-shift2}, that 
each connected component of $\BB(\lambda)$ is ``isomorphic'' 
(up to a shift of weights by a constant multiple of $\delta$) 
to the connected component $\BB_{0}(\lambda)$. 
In other words, 
the crystal graph of $\BB(\lambda)$ consists of 
``copies'' of that of $\BB_{0}(\lambda)$. 
Furthermore, it follows from Theorem~\ref{thm:comps} that 
these copies are parametrized by the sequences 
$(N_{1},\,N_{2},\,\dots,\,N_{s-1})$ of 
nonnegative integers 
satisfying the condition that 
$N_{u}-N_{u+1} \in 
   \sum_{j \in I_{0}(\lambda,p_{u})} 
   m_{j}d_{j} \BZ_{\ge 0}$ 
for $1 \le u \le s-1$, where $N_{s}:=0$. 
Later in Section~\ref{sec:aff}, 
we will study the $P$-crystal structure of 
the connected component $\BB_{0}(\lambda)$. 

%
\subsection{Some preliminary lemmas.}
\label{subsec:prelem}
The following lemma immediately follows from 
Lemma~\ref{lem:pair0}.
%
%
\begin{lem} \label{lem:pair1}
{\rm (1) }
Let $\xi \in \prr$ be a positive real root of 
the form $\xi=-\beta+c_{\beta}\delta$ 
with $\beta \in \pfr$, and 
let $\nu \in P$ be 
a level-zero integral weight. 
Then, $\nu(\xi^{\vee})=-\nu(\beta^{\vee})$, and 
$r_{\xi}(\nu)=
 \nu-\nu(\beta^{\vee})\beta+
 \nu(\beta^{\vee})c_{\beta}\delta$.

\vsp

\noindent {\rm (2) }
Let $\xi \in \prr$ be a positive real root of 
the form $\xi=\frac{1}{2}(-\beta+\delta)$ 
with $\beta \in \pfr$, and 
let $\nu \in P$ be a level-zero integral weight. 
Then, $\nu(\xi^{\vee})=-2\nu(\beta^{\vee})$, and 
$r_{\xi}(\nu)=
 \nu-\nu(\beta^{\vee})\beta+
 \nu(\beta^{\vee})\delta$.
\end{lem}
We derive the next lemma 
from Lemma~\ref{lem:pair1}.
%
%
%
\begin{lem} \label{lem:refl}
Let $\xi \in \prr$ be a positive real root of one 
of the following forms\,{\rm:} 
$\xi=\beta$ with $\beta \in \pfr$, 
$\xi=-\beta+c_{\beta}\delta$ 
with $\beta \in \pfr$, or 
$\xi=\frac{1}{2}(-\beta+\delta)$ 
with $\beta \in \pfr$. 
Let $\nu \in P$ be a level-zero 
integral weight. Then, we have
$r_{\xi}(\nu)=r_{\beta}(\nu)+n\delta$
 for some $n \in \BZ$.
\end{lem}
%
%
%
%
\begin{dfn}[see Lemma~\ref{lem:Worbit}] \label{dfn:fp}
Let $\lambda$ be the (fixed) level-zero dominant 
integral weight of the form \eqref{eq:lambda}, 
and write each $\nu \in W\lambda$ as 
$\nu=\lambda-\alpha+n\delta$ with
$\alpha \in \fin{Q}_{+}$ and $n \in \BZ$. 
We set 
\begin{equation}
\fp(\nu):=\alpha, \qquad D(\nu):=n, 
\end{equation}
and call the $\fp(\nu)=\alpha$ above 
the finite part of $\nu$.
\end{dfn}
For 
$\alpha=\sum_{j \in I_{0}}u_{j}\alpha_{j} 
 \in \fin{Q}_{+}$, we define the support 
$\supp(\alpha) \subset I_{0}$ of $\alpha$ by:
\begin{equation}
\supp(\alpha)=
 \bigl\{j \in I_{0} \mid u_{j} \ne 0\bigr\}.
\end{equation}
%
%
%
%
\begin{lem} \label{lem:fp}
{\rm (1) } 
Let $\xi \in \prr$ be 
a positive real root of the form\,{\rm:}
$\xi=\beta \in \pfr$, and assume that 
$\nu \in W\lambda$ satisfies 
$\nu(\xi^{\vee}) \in \BZ_{< 0}$. Then we have 
%
%
\begin{equation} \label{eq:fp-1}
\begin{cases}
\supp(\fp(\nu))=
\supp\bigl(\fp(r_{\xi}(\nu))\bigr) \cup 
\supp(\beta), & \\[1.5mm]
D(r_{\xi}(\nu))=D(\nu). & 
\end{cases}
\end{equation}

\vsp

\noindent {\rm (2) } 
Let $\xi \in \prr$ be 
a positive real root of the form\,{\rm:}
$\xi=-\beta+c_{\beta}\delta$ 
with $\beta \in \pfr$, and assume that 
$\nu \in W\lambda$ satisfies 
$\nu(\xi^{\vee}) \in \BZ_{< 0}$. Then we have 
\begin{equation} \label{eq:fp-2}
\begin{cases}
\supp\bigl(\fp(r_{\xi}(\nu))\bigr)=
\supp(\fp(\nu)) \cup \supp(\beta), & \\[1.5mm]
D(r_{\xi}(\nu))=
D(\nu)+\nu(\beta^{\vee})c_{\beta}. & 
\end{cases}
\end{equation}

\vsp

\noindent {\rm (3) } 
Let $\xi \in \prr$ be a positive real root of the form\,{\rm:}
$\xi=\frac{1}{2}(-\beta+\delta)$ 
with $\beta \in \pfr$, and assume that 
$\nu \in W\lambda$ satisfies 
$\nu(\xi^{\vee}) \in \BZ_{< 0}$. Then we have 
%
%
\begin{equation} \label{eq:fp-3}
\begin{cases}
\supp\bigl(\fp(r_{\xi}(\nu))\bigr)=
\supp(\fp(\nu)) \cup \supp(\beta), & \\[1.5mm]
D(r_{\xi}(\nu))=
D(\nu)+\nu(\beta^{\vee}). & 
\end{cases}
\end{equation}
\end{lem}

\begin{proof}
(1) \, 
First, note that $\nu=\lambda-\fp(\nu)+D(\nu)\delta$ 
by the definition of $\fp(\nu)$ and $D(\nu)$ 
(see Definition~\ref{dfn:fp}). 
Therefore, we have 
\begin{equation*}
r_{\xi}(\nu)
 =\nu-\nu(\xi^{\vee})\xi
 =\nu-\nu(\beta^{\vee})\beta
 =\lambda-\bigl\{
   \fp(\nu)+\nu(\beta^{\vee}) \beta\bigr\}
   +D(\nu)\delta.
\end{equation*}
Hence we obtain that 
$\fp(r_{\xi}(\nu))=\fp(\nu)+
\nu(\beta^{\vee})\beta$ and 
$D(r_{\xi}(\nu))=D(\nu)$. 
Since $\nu(\beta^{\vee})=
\nu(\xi^{\vee}) \in \BZ_{< 0}$ by assumption, 
we deduce that $\supp(\fp(\nu))=
\supp\bigl(\fp(r_{\xi}(\nu))\bigr) \cup \supp(\beta)$.
This proves part~(1). 

\vsp

\noindent (2) \, 
Because 
$r_{\xi}(\nu) = \nu-\nu(\beta^{\vee})\beta+
   \nu(\beta^{\vee})c_{\beta}\delta$ 
by Lemma~\ref{lem:pair1}\,(1), and 
$\nu=\lambda-\fp(\nu)+D(\nu)\delta$ by definition, 
we have
\begin{equation*}
r_{\xi}(\nu)=\lambda-\bigl\{
   \fp(\nu)+\nu(\beta^{\vee}) \beta\bigr\}
   +\bigl\{D(\nu)+
   \nu(\beta^{\vee})c_{\beta}\bigr\}\delta.
\end{equation*}
Hence we obtain that $\fp(r_{\xi}(\nu))=\fp(\nu)+
\nu(\beta^{\vee})\beta$ and 
$D(r_{\xi}(\nu))=
D(\nu)+\nu(\beta^{\vee})c_{\beta}$. 
Since $\nu(\beta^{\vee}) \in \BZ_{> 0}$ 
by Lemma~\ref{lem:pair1}\,(1), 
we deduce that 
$\supp\bigl(\fp(r_{\xi}(\nu))\bigr)=
\supp(\fp(\nu)) \cup \supp(\beta)$. 
This proves part~(2). 

\vsp

\noindent (3) \, 
Because 
$r_{\xi}(\nu)=
 \nu-\nu(\beta^{\vee})\beta+\nu(\beta^{\vee})\delta$ 
by Lemma~\ref{lem:pair1}\,(2), and 
$\nu=\lambda-\fp(\nu)+D(\nu)\delta$ by definition, 
we have
\begin{equation*}
r_{\xi}(\nu)
  =\lambda-\bigl\{
   \fp(\nu)+\nu(\beta^{\vee}) \beta\bigr\}
   +\bigl\{D(\nu)+
   \nu(\beta^{\vee})\bigr\}\delta.
\end{equation*}
Hence we obtain that 
$\fp(r_{\xi}(\nu))=\fp(\nu)+
\nu(\beta^{\vee})\beta$ and 
$D(r_{\xi}(\nu))=
D(\nu)+\nu(\beta^{\vee})$. 
Because $\nu(\beta^{\vee}) \in \BZ_{> 0}$ 
by Lemma~\ref{lem:pair1}\,(2), 
we deduce that 
$\supp\bigl(\fp(r_{\xi}\nu)\bigr)=
\supp(\fp(\nu)) \cup \supp(\beta)$. 
This proves part~(3). 
\end{proof}
%
%
\begin{lem} \label{lem:f-pair}
Let $\beta=\sum_{j \in I_{0}} u_{j}\alpha_{j}$ 
be an element of $\pfr$, 
with $u_{j} \in \BZ_{\ge 0}$ for $j \in I_{0}$. Then,
$u_{j}(\alpha_{j},\alpha_{j})/(\beta,\beta)
\in \BZ_{\ge 0}$ for all $j \in I_{0}$, and 
$\lambda(\beta^{\vee})
 =\sum_{j \in I_{0}} m_{j} 
  u_{j}(\alpha_{j},\alpha_{j})/(\beta,\beta)$.
\end{lem}

\begin{proof}
Since $u_{j} \in \BZ_{\ge 0}$ for all $j \in I_{0}$, 
it follows from \cite[Proposition~5.1\,(d)]{Kac} that 
$u_{j}(\alpha_{j},\alpha_{j})/(\beta,\beta) 
\in \BZ_{\ge 0}$ for all $j \in I_{0}$. 
Furthermore, we deduce that
\begin{align*}
\lambda(\beta^{\vee})
 & 
 =\frac{2(\lambda,\beta)}{(\beta,\beta)}
 =\sum_{j \in I_{0}} 
  \frac{2u_{j}(\lambda,\alpha_{j})}{(\beta,\beta)}
 =\sum_{j \in I_{0}} 
  \frac{2(\lambda,\alpha_{j})}{(\alpha_{j},\alpha_{j})}
  \cdot 
  \frac{u_{j}(\alpha_{j},\alpha_{j})}{(\beta,\beta)}
\\[3mm]
 &
 = \sum_{j \in I_{0}} 
  \lambda(h_{j})
  \frac{u_{j}(\alpha_{j},\alpha_{j})}{(\beta,\beta)}
 = \sum_{j \in I_{0}} m_{j}
  \frac{u_{j}(\alpha_{j},\alpha_{j})}{(\beta,\beta)}.
\end{align*}
This proves the lemma.
\end{proof}
%
%
%
%
\subsection{Necessary condition for the existence of 
   a $(q/p)$-chain.}
\label{subsec:nec}

Let $\lambda$ be the (fixed) level-zero dominant 
integral weight of the form \eqref{eq:lambda}.
This subsection is devoted to proving 
the following proposition.
%
%
%
%
\begin{prop} \label{prop:nec}
Let $1 \le q < p$ be coprime integers, 
and let $N \in \BZ$. 
Assume that there exists a $(q/p)$-chain for 
$(\lambda, \lambda+N\delta)$. Then, 
$N \in \sum_{
 j \in I_{0}(\lambda,p)} 
 m_{j}d_{j}\BZ_{\ge 0}$.
\end{prop}
From Proposition~\ref{prop:nec} 
along with Lemma~\ref{lem:d-shift}, 
we obtain the following corollary.
%
%
\begin{cor} \label{cor:nec}
Let $1 \le q < p$ be coprime integers, 
and let $N^{\prime},\,N^{\prime\prime} \in \BZ$. 
Assume that there exists a $(q/p)$-chain for 
$(\lambda-N^{\prime}\delta, \, 
  \lambda-N^{\prime\prime}\delta)$. 
Then, 
$N^{\prime}-N^{\prime\prime} 
\in \sum_{
 j \in I_{0}(\lambda,p)} 
 m_{j}d_{j}\BZ_{\ge 0}$.
\end{cor}

Now, let $\lambda=\nu_{0} > \nu_{1} > \cdots > 
\nu_{k}=\lambda+N\delta$ be 
a $(q/p)$-chain for $(\lambda,\lambda+N\delta)$, 
with $\xi_{1},\,\xi_{2},\,\dots,\,\xi_{k}$ 
the corresponding positive real roots; 
recall from Definition~\ref{dfn:achain} that 
$\nu_{l-1} > \nu_{l}=r_{\xi_{l}}(\nu_{l-1})$, with
$\dist (\nu_{l-1},\nu_{l})=1$, and 
$\nu_{l-1}(\xi_{l}^{\vee}) \in p\BZ_{< 0}$ 
for all $1 \le l \le k$. 
For each $1 \le l \le k$, 
let us denote by $\beta_{l} \in \pfr$ 
the finite part $\ba{\xi_{l}}$ of $\xi_{l} \in \prr$ 
(see Definition~\ref{dfn:fp-rr}).
Since $\dist (\nu_{l-1},\nu_{l})=1$, we know from 
Lemma~\ref{lem:dist=1} that $\xi_{l}$ is of one of 
the following forms: $\xi_{l}=\beta_{l}$, 
$\xi_{l}=-\beta_{l}+c_{\beta_{l}}\delta$, or 
$\xi_{l}=\frac{1}{2}(-\beta_{l}+\delta)$.
We define subsets $\Pos$, $\Neg_{1}$, $\Neg_{2}$, 
and $\Neg$ of $\bigl\{1,\,2,\,\dots,\,k\bigr\}$ as follows:
%
%
\begin{align} 
& 
\Pos=
\bigl\{1 \le l \le k \mid 
\xi_{l}=\beta_{l}\bigr\}, \label{eq:pos} \\[1mm]
& 
\Neg_{1}=
\bigl\{1 \le l \le k \mid 
\xi_{l}=-\beta_{l}+c_{\beta_{l}}\delta 
\bigr\}, \label{eq:neg1} \\[1mm]
&
\Neg_{2}=
\bigl\{1 \le l \le k \mid 
\xi_{l}=\tfrac{1}{2}(-\beta_{l}+\delta)
\bigr\}, \label{eq:neg2} \\[1mm]
& 
\Neg=\Neg_{1} \cup \Neg_{2}. \label{eq:neg}
\end{align}
Note that $\Neg_{2} = \emptyset$ 
unless $\Fg$ is of type $A^{(2)}_{2\ell}$. 
%
%
%
%
\begin{rem} \label{rem:fp}
We know from Lemma~\ref{lem:fp} that 
%
%
\begin{equation} \label{eq:r-fp-1}
\begin{cases}
\supp(\fp(\nu_{l-1}))=
\supp(\fp(\nu_{l})) \cup \supp(\beta_{l}) 
& \text{if $l \in \Pos$}, \\[1.5mm]
\supp(\fp(\nu_{l}))=
\supp(\fp(\nu_{l-1})) \cup \supp(\beta_{l})
& \text{if $l \in \Neg$}, 
\end{cases}
\end{equation}
and that
%
%
\begin{equation} \label{eq:r-fp-2}
D(\nu_{l})=
\begin{cases}
D(\nu_{l-1}) & \text{if \, } l \in \Pos, \\[1mm]
D(\nu_{l-1})+\nu_{l-1}(\beta_{l}^{\vee})c_{\beta_{l}}
 & \text{if \, } l \in \Neg_{1}, \\[1mm]
D(\nu_{l-1})+\nu_{l-1}(\beta_{l}^{\vee})
 & \text{if \, } l \in \Neg_{2}.
\end{cases}
\end{equation}
\end{rem}
%
%
\begin{lem} \label{lem:pQ}
If $l \in \Pos$, then 
$\fp(\nu_{l-1})-\fp(\nu_{l}) \in p \fin{Q}_{+}$. 
\end{lem}

\begin{proof}
We see from the proof of part (1) of 
Lemma~\ref{lem:fp} that 
$\fp(\nu_{l})=\fp(r_{\xi_{l}}(\nu_{l-1}))=
\fp(\nu_{l-1})+
\nu_{l-1}(\beta_{l}^{\vee})\beta_{l}$.
Also, it follows from the definition of 
a $(q/p)$-chain that 
$\nu_{l-1}(\beta_{l}^{\vee})=
\nu_{l-1}(\xi_{l}^{\vee}) \in p\BZ_{< 0}$
(note that $l \in \Pos$). 
Therefore, we obtain that
$\fp(\nu_{l-1})-\fp(\nu_{l})=
-\nu_{l-1}(\beta_{l}^{\vee})\beta_{l} \in 
p\fin{Q}_{+}$.
This proves the lemma. 
\end{proof}
We set $I_{0}(\lambda,p):=
\bigl\{i \in I_{0} \mid m_{i} \in p\BZ\bigr\}$. 

%
\begin{lem} \label{lem:supp}
The support $\supp (\beta_{l})$ of $\beta_{l}$ is contained 
in the set $I_{0}(\lambda,p)$ for all $1 \le l \le k$. 
\end{lem}

\begin{proof}
Suppose that there exists 
$1 \le l \le k$
such that $\supp (\beta_{l}) 
\not\subset I_{0}(\lambda,p)$. We set 
\begin{equation}
\Phi:=\bigl\{ 1 \le l \le k \mid 
 \supp (\beta_{l}) \not\subset I_{0}(\lambda,p)\bigr\}
 \ne \emptyset. 
\end{equation}
%
%
\begin{claim} \label{c:supp01}
Neither $\Phi \cap \Pos$ nor
$\Phi \cap \Neg$ is empty. 
Moreover, if we set $l_{0}:=\max (\Phi \cap \Pos)$ and 
$k_{0}:=\max (\Phi \cap \Neg)$, 
then we have $k_{0} < l_{0}$. 
\end{claim}

First, we note that both of the sets $\Phi \cap \Pos$ 
and $\Phi \cap \Neg$ are not empty, 
since $\Phi \ne \emptyset$. Now, suppose that 
$\Phi \cap \Neg = \emptyset$. Then we have 
$\Phi \cap \Pos \ne \emptyset$. 
Set $l_{0}^{\prime}:=\min \bigl(\Phi \cap \Pos)$, 
and let $i \in \supp(\beta_{l_{0}^{\prime}})$ 
be such that $m_{i} \notin p\BZ$. 
Since
$\supp(\fp(\nu_{l_{0}^{\prime}-1}))
 =\supp(\fp(\nu_{l_{0}^{\prime}}))
  \cup \supp(\beta_{l_{0}^{\prime}})$
by \eqref{eq:r-fp-1}, it follows that 
$i \in \supp(\fp(\nu_{l_{0}^{\prime}-1}))$. 
Then, because 
$\supp(\fp(\nu_{0}))=
\supp(\fp(\lambda))=\supp(0)=\emptyset$,
there exists $1 \le l \le l_{0}^{\prime}-1$ 
such that $i \notin 
\supp(\fp(\nu_{l-1}))$ and $i \in 
\supp(\fp(\nu_{l}))$. 
Hence we have 
$\supp(\fp(\nu_{l-1})) \not\supset 
\supp(\fp(\nu_{l}))$. 
Therefore, we deduce from 
\eqref{eq:r-fp-1} that $l \in \Neg$, and 
$\supp(\fp(\nu_{l}))
=\supp(\fp(\nu_{l-1})) \cup \supp(\beta_{l})$.
Since $i \notin \supp(\fp(\nu_{l-1}))$ and 
$i \in \supp(\fp(\nu_{l}))$, we see that 
$i \in \supp(\beta_{l})$. 
As a result, we obtain that 
$l \in \Phi \cap \Neg$, and hence 
$\Phi \cap \Neg \ne \emptyset$, 
which is a contradiction. 
Thus we conclude that 
$\Phi \cap \Neg \ne \emptyset$. 

Set $k_{0}:=\max(\Phi \cap \Neg)$, and 
let $i \in \supp(\beta_{k_{0}})$ 
be such that $m_{i} \notin p\BZ$. 
Since 
$\supp(\fp(\nu_{k_{0}}))
=\supp(\fp(\nu_{k_{0}-1}))
\cup \supp(\beta_{k_{0}})$
by \eqref{eq:r-fp-1}, 
it follows that 
$i \in \supp(\fp(\nu_{k_{0}}))$. 
Then, because 
$\supp(\fp(\nu_{k}))=
\supp(\fp(\lambda+N\delta))=\supp(0)=\emptyset$,
there exists $k \ge l > k_{0}$ such that $i \in 
\supp(\fp(\nu_{l-1}))$ and $i \notin 
\supp(\fp(\nu_{l}))$. 
This implies that $\supp(\fp(\nu_{l-1})) \not\subset 
\supp(\fp(\nu_{l}))$. 
Hence we deduce from \eqref{eq:r-fp-1} that 
$l \in \Pos$, and that
$\supp(\fp(\nu_{l-1}))
=\supp(\fp(\nu_{l})) \cup \supp(\beta_{l})$.
Since $i \in \supp(\fp(\nu_{l-1}))$ and 
$i \not\in \supp(\fp(\nu_{l}))$, we see that 
$i \in \supp(\beta_{l})$.
As a result, we obtain that 
$l \in \Phi \cap \Pos$, and hence 
$\Phi \cap \Pos \ne \emptyset$. 
Furthermore, since $l > k_{0}$, 
we have $l_{0} > k_{0}$. 
This proves Claim~\ref{c:supp01}. 

\vsp\vsp

Since $l_{0} \in \Pos$, we see that 
$\nu_{l_{0}-1}$ and $\nu_{l_{0}}$ are 
contained in the same $\fin{W}$-orbit.
Let $\lambda^{\prime}$ be the unique 
level-zero dominant integral weight in 
$\fin{W}\nu_{l_{0}-1}=\fin{W}\nu_{l_{0}}$. 
Then, it follows that $\lambda^{\prime}=
\lambda+n^{\prime}\delta$ 
for some $n^{\prime} \in \BZ$. 
Indeed, since $W=\fin{W} \ltimes T$, 
$\nu_{l_{0}} \in W\lambda$ can be written as 
$w^{\prime}w^{\prime\prime} \lambda$ 
for some $w^{\prime} \in \fin{W}$ and 
$w^{\prime\prime} \in T$. 
Hence we have $(w^{\prime})^{-1}\nu_{l_{0}}= 
w^{\prime\prime}\lambda \in T\lambda$. 
Since $T\lambda \subset \lambda+\BZ\delta$ 
by Lemma~\ref{lem:dlam}, we have 
$w^{\prime\prime}\lambda=
\lambda+n^{\prime}\delta$
for some $n^{\prime} \in \BZ$. 
Therefore, 
$(w^{\prime})^{-1}\nu_{l_{0}}=
\lambda+n^{\prime}\delta$ is 
a level-zero dominant integral weight contained in 
$\fin{W}\nu_{l_{0}}$. 
Thus we conclude that 
$\lambda^{\prime}=(w^{\prime})^{-1}\nu_{l_{0}}=
\lambda+n^{\prime}\delta$. 

Now, let $w_{l_{0}-1} \in \fin{W}$ 
(resp., $w_{l_{0}} \in \fin{W}$) be 
the shortest element in $\fin{W}$ such that 
$\nu_{l_{0}-1}=w_{l_{0}-1}(\lambda+n^{\prime}\delta)$ 
(resp., $\nu_{l_{0}}=w_{l_{0}}
(\lambda+n^{\prime}\delta)$). 
%
%
\begin{claim} \label{c:supp02}
Assume that 
$w_{l_{0}}=r_{j_{K}}r_{j_{K-1}} \cdots r_{j_{1}}$ 
is a reduced expression of $w_{l_{0}} \in \fin{W}$, where 
$j_{1},\,\dots,\,j_{K-1},\,j_{K} \in I_{0}$.
Then, $j_{1},\,j_{2},\,\dots,\,j_{K}$ are 
all contained in $I_{0}(\lambda,p)$.
\end{claim}

We set $\mu_{0}:=\lambda+n^{\prime}\delta$, and 
$\mu_{L}:=
r_{j_{L}}r_{j_{L-1}} \cdots r_{j_{1}}
(\lambda+n^{\prime}\delta)$ 
for $1 \le L \le K$; note that $\mu_{K}=\nu_{l_{0}}$. 
Since $r_{j_{K}}r_{j_{K-1}} \cdots r_{j_{1}}$ is 
a reduced expression of 
the shortest element $w_{l_{0}} \in \fin{W}$ such that 
$\nu_{l_{0}}=w_{l_{0}}(\lambda+n^{\prime}\delta)$, 
it follows that $\mu_{L-1}(h_{j_{L}}) \ne 0$ 
for all $1 \le L \le K$. 
Suppose now that $\mu_{L-1}(h_{j_{L}}) < 0$ 
for some $1 \le L \le K$. Then we have
$r_{j_{1}} \cdots r_{j_{L-1}}(\alpha_{j_{L}}) \in -\pfr$, 
since $\lambda+n^{\prime}\delta$ is 
level-zero dominant. 
Hence, using \cite[Lemma~3.10]{Kac}, 
we deduce that 
$r_{j_{1}} \cdots r_{j_{L-1}}r_{j_{L}}$ is 
not a reduced expression, which contradicts 
our assumption that 
$r_{j_{K}}r_{j_{K-1}} \cdots r_{j_{1}}$ is 
a reduced expression of $w_{l_{0}}$. 
Thus we obtain that $\mu_{L-1}(h_{j_{L}}) > 0$ 
for all $1 \le L \le K$, and hence that 
%
%
\begin{equation} \label{eq:supp02-1}
\nu_{l_{0}}=\mu_{K}=
\lambda-\sum_{L=1}^{K}
\underbrace{\mu_{L-1}(h_{j_{L}})}_{> 0} 
\alpha_{j_{L}} + n^{\prime} \delta. 
\end{equation}
Therefore, 
if the assertion of Claim~\ref{c:supp02} is 
false, then $\supp(\fp(\nu_{l_{0}}))$ must contain
some $i \in I_{0}$ such that 
$i \notin I_{0}(\lambda,p)$. 
Then, since $\fp(\nu_{k})=\fp(\lambda+N\delta)=0$, 
there exists $k \ge l > l_{0}$ such that 
$i \in \supp(\fp(\nu_{l-1}))$ and 
$i \not\in \supp(\fp(\nu_{l}))$. 
This implies that 
$\supp(\fp(\nu_{l-1})) \not\subset 
\supp(\fp(\nu_{l}))$. 
Hence we deduce from \eqref{eq:r-fp-1} that 
$l \in \Pos$, and 
$\supp(\fp(\nu_{l-1}))
=\supp(\fp(\nu_{l})) \cup \supp(\beta_{l})$.
Since $i \in \supp(\fp(\nu_{l-1}))$ and 
$i \not\in \supp(\fp(\nu_{l}))$, we see that 
$i \in \supp(\beta_{l})$. 
As a result, we obtain that $l \in \Phi \cap \Pos$ and 
$l > l_{0}$, which contradicts the definition:
$l_{0}=\max (\Phi \cap \Pos)$. 
This proves Claim~\ref{c:supp02}. 

\vsp\vsp

Because 
$\nu_{l_{0}-1} > \nu_{l_{0}}$ and 
$\dist(\nu_{l_{0}-1},\nu_{l_{0}})=1$, 
we see from \cite[Remark~4.2]{L2} that 
$w_{l_{0}-1}$ is greater than $w_{l_{0}}$ 
with respect to the Bruhat ordering on $\fin{W}$, and 
$\ell(w_{l_{0}-1})-\ell(w_{l_{0}})=1$.
Therefore, the reduced expression of $w_{l_{0}}$ above is 
obtained from a reduced expression of $w_{l_{0}-1}$
by removing a simple reflection $r_{i_{0}}$ 
for some $i_{0} \in I_{0}$. Namely, 
$w_{l_{0}-1} \in \fin{W}$ has 
a reduced expression of the form: $w_{l_{0}-1}=
r_{j_{K}} \cdots r_{j_{L_{0}+1}}r_{i_{0}}
r_{j_{L_{0}}} \cdots r_{j_{1}}$. 
%
%
\begin{claim} \label{c:supp025}
The $i_{0} \in I_{0}$ above is 
not contained in $I_{0}(\lambda,p)$. 
\end{claim}

It immediately follows from \eqref{eq:r-fp-1} that 
$\supp(\fp(\nu_{l_{0}-1}))
 =\supp(\fp(\nu_{l_{0}})) \cup 
  \supp(\beta_{l_{0}})$.
Since $\supp(\beta_{l_{0}}) \not\subset I_{0}(\lambda,p)$ by 
the definition of $l_{0}$, we have 
$\supp(\fp(\nu_{l_{0}-1})) \not\subset I_{0}(\lambda,p)$. 
By an argument similar to that used to obtain 
\eqref{eq:supp02-1}, we deduce that 
\begin{equation*}
\nu_{l_{0}-1}=w_{l_{0}-1}(\lambda+n^{\prime}\delta) 
\in 
\lambda-\left\{
\sum_{L=1}^{L_{0}} \BZ_{> 0}\alpha_{j_{L}} + 
\BZ_{>0} \alpha_{i_{0}} + 
\sum_{L=L_{0}+1}^{K} \BZ_{> 0}\alpha_{j_{L}}
\right\}
+ n^{\prime} \delta. 
\end{equation*}
Hence we have $\supp(\fp(\nu_{l_{0}-1}))=
 \bigl\{j_{1},\,j_{2},\,\dots,\,j_{K}\bigr\} \cup \{i_{0}\}$.
But, since $j_{1},\,j_{2},\,\dots,\,j_{K} \in I_{0}(\lambda,p)$ 
by Claim~\ref{c:supp02}, and 
$\supp(\fp(\nu_{l_{0}-1})) \not\subset I_{0}(\lambda,p)$, 
we conclude that $i_{0} \notin I_{0}(\lambda,p)$. 
This proves Claim~\ref{c:supp025}. 
%
%
%
%
\begin{claim} \label{c:supp03}
The coefficient of $\alpha_{i_{0}}$ in 
$\fp(\nu_{l_{0}-1})-\fp(\nu_{l_{0}})$ 
is contained in $m_{i_{0}}+p\BZ$. 
\end{claim}

From Claims~\ref{c:supp02} and \ref{c:supp025}, 
we have $j_{L} \ne i_{0}$ 
for all $1 \le L \le K$. 
Therefore, we see from \eqref{eq:supp02-1} that 
the coefficient of $\alpha_{i_{0}}$ 
in $\fp(\nu_{l_{0}})$ is equal to $0$. 

Let us compute the coefficient of 
$\alpha_{i_{0}}$ in $\fp(\nu_{l_{0}-1})$; 
recall that 
\begin{equation*}
\nu_{l_{0}-1}=
w_{l_{0}-1}(\lambda+n^{\prime}\delta)=
r_{j_{K}} \cdots r_{j_{L_{0}+1}}r_{i_{0}}
r_{j_{L_{0}}} \cdots r_{j_{1}}
(\lambda+n^{\prime}\delta).
\end{equation*}
First we derive inductively that 
\begin{equation*}
r_{j_{L_{0}}} \cdots r_{j_{1}}
 (\lambda+n^{\prime}\delta) 
 \in \lambda-\sum_{L=1}^{L_{0}} 
(m_{j_{L_{0}}}\BZ+\cdots+m_{j_{1}}\BZ)
 \alpha_{j_{L}} + n^{\prime}\delta. 
\end{equation*}
Since 
$j_{1},\,j_{2},\,\dots,\,j_{L_{0}} \in I_{0}(\lambda,p)$ 
by Claim~\ref{c:supp02}, 
it follows that
\begin{equation*}
r_{j_{L_{0}}} \cdots r_{j_{1}}
 (\lambda+n^{\prime}\delta) 
 \in \lambda-\sum_{L=1}^{L_{0}}
 p\BZ \alpha_{j_{L}} + n^{\prime}\delta.
\end{equation*}
Furthermore, we deduce that
\begin{align*}
r_{i_{0}}r_{j_{L_{0}}} \cdots r_{j_{1}}
(\lambda+n^{\prime}\delta) 
& \in 
\lambda-m_{i_{0}}\alpha_{i_{0}} 
 -\left\{
 \sum_{L=1}^{L_{0}} 
  p\BZ \alpha_{j_{L}} -
 \sum_{L=1}^{L_{0}} 
  p\BZ \alpha_{j_{L}}(h_{i_{0}})\alpha_{i_{0}}
 \right\} + n^{\prime}\delta \\[3mm]
& \subset \lambda- (m_{i_{0}}+p\BZ)\alpha_{i_{0}} -
 \sum_{L=1}^{L_{0}} 
  p\BZ \alpha_{j_{L}} + n^{\prime}\delta. 
\end{align*}
Here, again noting that $j_{L} \ne i_{0}$ for all 
$1 \le L \le K$, we see that 
the coefficient of $\alpha_{i_{0}}$ 
in the finite part $\fp(\nu_{l_{0}-1})$ of 
$\nu_{l_{0}-1}=r_{j_{K}} \cdots r_{j_{L_{0}+1}}
r_{i_{0}}r_{j_{L_{0}}} \cdots r_{j_{1}}
(\lambda+n^{\prime}\delta)$ is equal to 
that in the finite part of 
$r_{i_{0}}r_{j_{L_{0}}} \cdots 
r_{j_{1}}(\lambda+n\delta)$. 
Hence we conclude 
that the coefficient of $\alpha_{i_{0}}$ 
in the finite part $\fp(\nu_{l_{0}-1})$ of $\nu_{l_{0}-1}$ is 
contained in $m_{i_{0}} + p\BZ$. 

Thus the coefficient of $\alpha_{i_{0}}$ in
$\fp(\nu_{l_{0}-1})-\fp(\nu_{l_{0}})$ is 
contained in $m_{i_{0}} + p\BZ$.
This proves Claim~\ref{c:supp03}. 

\vsp\vsp

We know from 
Lemma~\ref{lem:pQ} that 
$\fp(\nu_{l_{0}-1})-\fp(\nu_{l_{0}}) \in p\fin{Q}_{+}$, 
since $l_{0}=\max (\Phi \cap \Pos) \in \Pos$. 
However, because the $i_{0} \in I_{0}$ is 
not contained in $I_{0}(\lambda,p)$ 
by Claim~\ref{c:supp025}, it follows from 
Claim~\ref{c:supp03} that the coefficient of 
$\alpha_{i_{0}}$ in $\fp(\nu_{l_{0}-1})-
\fp(\nu_{l_{0}})$ is not contained in $p\BZ$. 
This implies that $\fp(\nu_{l_{0}-1})-
\fp(\nu_{l_{0}}) \notin p\fin{Q}_{+}$, 
which is a contradiction. 
Thus we have proved Lemma~\ref{lem:supp}. 
\end{proof}
It follows from \eqref{eq:r-fp-2} that
%
%
\begin{equation} \label{eq:delta}
N=D(\lambda+N\delta)=D(\nu_{k})=
 \sum_{l \in \Neg_{1}} \nu_{l-1}(\beta_{l}^{\vee})c_{\beta_{l}} + 
 \sum_{l \in \Neg_{2}} \nu_{l-1}(\beta_{l}^{\vee}). 
\end{equation}
Therefore, 
in order to prove Proposition~\ref{prop:nec} 
(i.e., that $N \in \sum_{j \in I_{0}(\lambda,p)}
  m_{j}d_{j}\BZ_{\ge 0}$), 
it suffices to establish the following lemma. 
%
%
\begin{lem} \label{lem:pair2}
{\rm (1) } 
If $l \in \Neg_{1}$, then 
$\nu_{l-1}(\beta_{l}^{\vee})c_{\beta_{l}}
\in \sum_{j \in I_{0}(\lambda,p)} 
    m_{j}d_{j}\BZ_{\ge 0}$. 

\vsp

\noindent {\rm (2) } 
If $l \in \Neg_{2}$, then 
$\nu_{l-1}(\beta_{l}^{\vee}) 
\in \sum_{j \in I_{0}(\lambda,p)} 
    m_{j}d_{j}\BZ_{\ge 0}$. 
\end{lem}

\begin{proof}
We set 
$\beta_{l}^{\prime}:=
r_{\beta_{1}}r_{\beta_{2}} \cdots 
r_{\beta_{l-1}}(\beta_{l}) \in \fr$ 
for each $1 \le l \le k$; 
note that $\beta_{1}^{\prime}=\beta_{1}$.
%
%
%
\begin{claim} \label{c:p-2-0}
We have 
$(\beta_{l}^{\prime},\beta_{l}^{\prime})=
 (\beta_{l},\beta_{l})$, and 
$(\lambda,\beta_{l}^{\prime})=
(\nu_{l-1},\beta_{l})$ for all $1 \le l \le k$.
Hence 
$\nu_{l-1}(\beta_{l}^{\vee})=
 \lambda(\beta_{l}^{\prime\,\vee})$ 
for all $1 \le l \le k$. Moreover, 
$\beta_{l}^{\prime} \in \pfr$ 
for every $l \in \Neg$. 
\end{claim} 
%
Since the bilinear form $(\cdot\,,\,\cdot)$ 
is $W$-invariant, we have 
$(\beta_{l}^{\prime},\beta_{l}^{\prime})=
 (\beta_{l},\beta_{l})$. 
In addition, 
we deduce from Lemma~\ref{lem:refl} that 
$
r_{\xi_{1}}r_{\xi_{2}} \cdots 
r_{\xi_{l-1}}(\beta_{l})=
r_{\beta_{1}}r_{\beta_{2}} \cdots 
r_{\beta_{l-1}}(\beta_{l})+
     n^{\prime\prime}\delta =
\beta_{l}^{\prime}+
     n^{\prime\prime}\delta$
for some $n^{\prime\prime} \in \BZ$, 
and hence that 
\begin{align*}
(\nu_{l-1},\beta_{l}) & =
(r_{\xi_{l-1}}r_{\xi_{l-2}} \cdots 
r_{\xi_{1}}(\nu_{0}),\beta_{l})=
(r_{\xi_{l-1}}r_{\xi_{l-2}} \cdots 
r_{\xi_{1}}(\lambda),\beta_{l}) \\
& =(\lambda, 
r_{\xi_{1}}r_{\xi_{2}} \cdots 
r_{\xi_{l-1}}(\beta_{l}))=
(\lambda, \beta_{l}^{\prime}+
 n^{\prime\prime}\delta)=
(\lambda,\beta_{l}^{\prime}). 
\end{align*}
Therefore, 
using the equations
$(\beta_{l}^{\prime},\beta_{l}^{\prime})=
 (\beta_{l},\beta_{l})$ and 
$(\nu_{l-1},\beta_{l})=
 (\lambda,\beta_{l}^{\prime})$, 
we obtain that
\begin{equation*}
\nu_{l-1}(\beta_{l}^{\vee})=
 \frac{2(\nu_{l-1},\beta_{l})}{(\beta_{l},\beta_{l})}=
 \frac{2(\lambda,\beta_{l}^{\prime})}
      {(\beta_{l}^{\prime},\beta_{l}^{\prime})}=
\lambda(\beta_{l}^{\prime\,\vee}). 
\end{equation*}
Furthermore, because $\beta_{l}$ is 
the finite part $\ba{\xi_{l}}$ of 
$\xi_{l} \in \prr$, and 
$\nu_{l-1}(\xi_{l}^{\vee}) \in p\BZ_{< 0}$ 
by the definition of $(q/p)$-chains, 
we see from Lemma~\ref{lem:pair1} that 
$\nu_{l-1}(\beta_{l}^{\vee}) > 0$ for every $l \in \Neg$. 
As a result, we obtain that
$(\lambda,\beta_{l}^{\prime})=
 (\nu_{l-1},\beta_{l}) > 0$. 
Hence we conclude that 
$\beta_{l}^{\prime} \in \pfr$, 
since $\lambda$ is a level-zero dominant integral weight.
This proves Claim~\ref{c:p-2-0}. 

\vsp

The following claim 
immediately follows from Lemma~\ref{lem:supp}. 

%
\begin{claim} \label{c:p-2-1}
The support $\supp(\beta_{l}^{\prime})$ 
of $\beta_{l}^{\prime} \in \prr$ is contained in 
$I_{0}(\lambda,p)$ for every $l \in \Neg$.  
\end{claim}

We now prove part (1) of Lemma~\ref{lem:pair2}. 
Since $l \in \Neg_{1}$, we see 
from Claim~\ref{c:p-2-0} that 
$\beta_{l}^{\prime} \in \pfr$. 
Let us write $\beta_{l}^{\prime} \in \pfr$ as: 
$\beta_{l}^{\prime}=
\sum_{j \in I_{0}} u_{lj}\alpha_{j}$, 
with $u_{lj} \in \BZ_{\ge 0}$ for $j \in I_{0}$. 

First assume that $\Fg$ is nontwisted. 
Then we know from Remark~\ref{rem:cbeta}\,(1) and 
\eqref{eq:di} that 
$c_{\beta_{l}}=1$, and 
$d_{j}=1$ for all $j \in I_{0}$. 
From Claim~\ref{c:p-2-0} and 
Lemma~\ref{lem:f-pair}, we deduce that
\begin{equation*}
\nu_{l-1}(\beta_{l}^{\vee})c_{\beta_{l}} = 
\nu_{l-1}(\beta_{l}^{\vee}) = 
\lambda({\beta_{l}^{\prime}}^{\vee})=
  \sum_{j \in I_{0}}m_{j}
  \underbrace{
  \frac{u_{lj}(\alpha_{j},\alpha_{j})}
  {(\beta_{l}^{\prime},\beta_{l}^{\prime})}
  }_{\in \BZ_{\ge 0}}.
\end{equation*}
Furthermore, since $\supp(\beta_{l}^{\prime}) \subset 
I_{0}(\lambda,p)$ by Claim~\ref{c:p-2-1}, it follows that 
$u_{lj}=0$ unless $j \in I_{0}(\lambda,p)$. 
Hence we obtain that 
$\nu_{l-1}(\beta_{l}^{\vee})c_{\beta_{l}} \in 
\sum_{j \in I_{0}(\lambda,p)}
   m_{j}d_{j}\BZ_{\ge 0} = 
 \sum_{j \in I_{0}(\lambda,p)}
  m_{j}\BZ_{\ge 0}$. 

Next assume that $\Fg$ is twisted. 
Then we know from Remark~\ref{rem:cbeta}\,(2) and 
\eqref{eq:di} that $c_{\beta_{l}}=
(\beta_{l},\beta_{l})/2$, and 
%
%
\begin{equation} \label{eq:norm}
\frac{(\alpha_{j},\alpha_{j})}{2} = 
\begin{cases}
2d_{j} & 
 \text{if $\Fg$ is of type $A^{(2)}_{2\ell}$ 
       and $j=\ell$}, \\[3mm]
d_{j} & \text{otherwise}.
\end{cases}
\end{equation}
We have 
\begin{align*}
\nu_{l-1}(\beta_{l}^{\vee})c_{\beta_{l}} & = 
\frac{2(\nu_{l-1},\beta_{l})}{(\beta_{l},\beta_{l})} 
  c_{\beta_{l}}= (\nu_{l-1},\beta_{l}) = 
  (\lambda,\beta_{l}^{\prime}) 
  \quad \text{by Claim~\ref{c:p-2-0}} \\[3mm]
& = \sum_{j \in I_{0}} u_{lj}(\lambda,\alpha_{j}) 
  = \sum_{j \in I_{0}} u_{lj} 
  \frac{2(\lambda,\alpha_{j})}{(\alpha_{j},\alpha_{j})} \cdot 
  \frac{(\alpha_{j},\alpha_{j})}{2} 
  \\[3mm]
& = \sum_{j \in I_{0}} u_{lj}
  \lambda(h_{j}) 
  \frac{(\alpha_{j},\alpha_{j})}{2}
  = \sum_{j \in I_{0}} u_{lj} m_{j}
  \frac{(\alpha_{j},\alpha_{j})}{2}. 
\end{align*}
Note that $u_{lj} \in \BZ_{\ge 0}$ 
for all $j \in I_{0}$, and that
$(\alpha_{j},\alpha_{j})/2 \in d_{j}\BZ_{> 0}$ 
for all $j \in I_{0}$ by \eqref{eq:norm}. 
Also, since $\supp(\beta_{l}^{\prime}) \subset 
I_{0}(\lambda,p)$ by Claim~\ref{c:p-2-1}, it follows that 
$u_{lj}=0$ unless $j \in I_{0}(\lambda,p)$. 
Hence we obtain that 
$\nu_{l-1}(\beta_{l}^{\vee})c_{\beta_{l}} \in 
\sum_{j \in I_{0}(\lambda,p)}
  m_{j}d_{j}\BZ_{\ge 0}$. 
This completes the proof of 
part~(1) of Lemma~\ref{lem:pair2}. 

Let us prove part (2) of Lemma~\ref{lem:pair2}; 
note that $\Fg$ is of type $A^{(2)}_{2\ell}$, 
since $\Neg_{2} \ne \emptyset$.
Therefore, we know from Remark~\ref{rem:cbeta}\,(2) and 
\eqref{eq:di} that $d_{j}=1$ for all $j \in I_{0}$. 
We see from Claim~\ref{c:p-2-0} and 
Lemma~\ref{lem:f-pair} that 
\begin{equation*}
\nu_{l-1}(\beta_{l}^{\vee}) = 
\lambda({\beta_{l}^{\prime}}^{\vee})= 
  \sum_{j \in I_{0}}m_{j}
   \underbrace{
   \frac{u_{lj}(\alpha_{j},\alpha_{j})}
   {(\beta_{l}^{\prime},\beta_{l}^{\prime})}
   }_{\in \BZ_{\ge 0}}.
\end{equation*}
Furthermore, since $\supp(\beta_{l}^{\prime}) \subset 
I_{0}(\lambda,p)$ by Claim~\ref{c:p-2-1}, 
it follows that $u_{lj}=0$ unless 
$j \in  I_{0}(\lambda,p)$. 
Hence we obtain that 
$\nu_{l-1}(\beta_{l}^{\vee}) \in 
 \sum_{j \in I_{0}(\lambda,p)}
  m_{j}d_{j} \BZ_{\ge 0} = 
 \sum_{j \in I_{0}(\lambda,p)}
  m_{j}\BZ_{\ge 0}$. 
This completes 
the proof of part~(2) of Lemma~\ref{lem:pair2}. 
\end{proof}

As mentioned just above Lemma~\ref{lem:pair2}, 
Proposition~\ref{prop:nec} 
immediately follows from Lemma~\ref{lem:pair2} and 
\eqref{eq:delta}. Thus we have completed the proof of 
Proposition~\ref{prop:nec}.

%
\subsection{Sufficient condition for the existence of a $(q/p)$-chain.}
\label{subsec:suf}
Let $\lambda$ be the (fixed) level-zero 
dominant integral weight of the form \eqref{eq:lambda}. 
This subsection is devoted to proving 
the following proposition.
%
%
%
%
\begin{prop} \label{prop:suf}
Let $1 \le q < p$ be coprime integers, 
and let $N \in \BZ_{\ge 0}$. 
Then, there exists a $(q/p)$-chain for 
$(\lambda,\,\lambda+N\delta)$ if 
$N \in \sum_{
 j \in I_{0}(\lambda,p)} 
  m_{j}d_{j} \BZ_{\ge 0}$.
\end{prop}
From Proposition~\ref{prop:suf} along with 
Lemma~\ref{lem:d-shift}, we obtain 
the following corollary. 
%
%
\begin{cor} \label{cor:suf}
Let $1 \le q < p$ be coprime integers, 
and let $N^{\prime},\,N^{\prime\prime} \in \BZ$. 
Then, there exists a $(q/p)$-chain for 
$(\lambda-N^{\prime}\delta, \, 
  \lambda-N^{\prime\prime}\delta)$ if 
$N^{\prime}-N^{\prime\prime} \in 
\sum_{
 j \in I_{0}(\lambda,p)} 
 m_{j}d_{j} \BZ_{\ge 0}$.
\end{cor}
In order to prove Proposition~\ref{prop:suf}, 
we need some lemmas. From now till
the end of the proof of Proposition~\ref{prop:mi-ch} below, 
we fix $i \in \Supp_{\ge 2}(\lambda)$ arbitrarily. 
We set 
$\Supp(\lambda):=
\bigl\{i \in I_{0} \mid m_{i} \ge 1\bigr\} \subset I_{0}$,
and $S(\lambda,i):=
\bigl(I_{0} \setminus \Supp(\lambda)\bigr) \cup \{i\}$. 
%
%
\begin{dfn} \label{dfn:prr-lam} 
Let $\prr(\lambda,i)$ denote
the set of positive real roots $\zeta \in \prr$, 
with $\gamma:=\ba{\zeta} \in \pfr$ its finite part, 
satisfying the the following conditions: 
\\[1.5mm]
(a) $\zeta$ is either of the form 
$\zeta=-\gamma+c_{\gamma}\delta$, or of the form 
$\zeta=\frac{1}{2}(-\gamma+\delta)$, 
\\[1.5mm]
(b) the support $\supp(\gamma)$ of $\gamma$ is contained in 
$S(\lambda,i)$, and contains the (fixed) 
$i \in \Supp_{\ge 2}(\lambda)$, 
\\[1.5mm]
(c) $r_{\zeta}(\lambda)=
    \lambda-\alpha+m_{i}d_{i}\delta$ 
    for some $\alpha \in \fin{Q}_{+}$. 
\end{dfn}
%
%
%
%
\begin{lem} \label{lem:prr-lam}
The set $\prr(\lambda,i)$ is not empty. 
\end{lem}

\begin{proof}
First assume that 
$\Fg$ is of type $A^{(2)}_{2\ell}$, and that 
the simple root corresponding to 
the $i \in \Supp_{\ge 2}(\lambda)$ is 
the long simple root (i.e., $i=\ell$). 
We set $\zeta:=\frac{1}{2}(-\alpha_{i}+\delta)$; 
note that $\zeta$ is indeed a positive real root 
by Proposition~\ref{prop:real}\,(2). 
It is obvious that 
this $\zeta \in \prr$ 
satisfies conditions~(a) and (b) of 
Definition~\ref{dfn:prr-lam}. 
Also, we see from Lemma~\ref{lem:pair1}\,(2) that 
$r_{\zeta}(\lambda)
 =\lambda-\lambda(h_{i})\alpha_{i}+\lambda(h_{i})\delta
 =\lambda-m_{i}\alpha_{i}+m_{i}\delta$.
Since $d_{i}=1$ by \eqref{eq:di}, 
it follows that $\zeta$ 
satisfies condition~(c) of Definition~\ref{dfn:prr-lam}. 
Thus we obtain that $\zeta \in \prr(\lambda,i)$, and hence 
$\prr(\lambda,i) \ne \emptyset$. 

Next assume that $\Fg$ is not 
of type $A^{(2)}_{2\ell}$, or 
$\Fg$ is of type $A^{(2)}_{2\ell}$ and 
the simple root corresponding to 
the $i \in \Supp_{\ge 2}(\lambda)$ is not 
the long simple root (i.e., $i \ne \ell$). 
We set $\zeta:=-\alpha_{i}+c_{\alpha_{i}}\delta$; 
note that $\zeta$ is indeed a positive real root 
by Proposition~\ref{prop:real}. It is obvious that 
this $\zeta \in \prr$ satisfies conditions~(a) and (b) of 
Definition~\ref{dfn:prr-lam}. 
Also, 
we see from Lemma~\ref{lem:pair1}\,(1) that 
$r_{\zeta}(\lambda)
 =\lambda-\lambda(h_{i})\alpha_{i}+
  \lambda(h_{i})c_{\alpha_{i}}\delta
 =\lambda-m_{i}\alpha_{i}+m_{i}c_{\alpha_{i}}\delta$.
Since $d_{i}=c_{\alpha_{i}}$ by \eqref{eq:di}, 
it follows that $\zeta$ satisfies condition~(c) of 
Definition~\ref{dfn:prr-lam}.
Thus we obtain that 
$\zeta \in \prr(\lambda,i)$, 
and hence $\prr(\lambda,i) \ne \emptyset$. 
This proves the lemma. 
\end{proof}
%
%
\begin{lem} \label{lem:pair3}
For every $\zeta \in \prr(\lambda,i)$, we have 
$\lambda(\zeta^{\vee}) \in m_{i}\BZ_{< 0}$. 
\end{lem}

\begin{proof}
Let $\gamma \in \pfr$ 
be the finite part $\ba{\zeta}$ of $\zeta \in \prr$. 
We deduce from condition~(a) of 
Definition~\ref{dfn:prr-lam}, 
using Lemma~\ref{lem:pair1}, 
that $\lambda(\zeta^{\vee})=
 -\lambda(\gamma^{\vee})$, or 
$-2\lambda(\gamma^{\vee})$.
Therefore, it suffices to show that 
$\lambda(\gamma^{\vee}) \in m_{i}\BZ_{> 0}$. 
Write $\gamma \in \pfr$ 
as $\sum_{j \in I_{0}} u_{j} \alpha_{j}$, 
with $u_{j} \in \BZ_{\ge 0}$ for $j \in I_{0}$. 
It follows from condition~(b) of 
Definition~\ref{dfn:prr-lam} that 
$u_{i} \in \BZ_{> 0}$ for 
the $i \in \Supp_{\ge 2}(\lambda)$, 
and $u_{j}=0$ for $j \notin S(\lambda,i)$. 
Hence, using Lemma~\ref{lem:f-pair}, we deduce that
\begin{equation*}
\lambda(\gamma^{\vee})
  =\sum_{j \in I_{0}} m_{j}
   \frac{u_{j}(\alpha_{j},\alpha_{j})}{(\gamma,\gamma)}
  = \sum_{j \in S(\lambda,i)} m_{j}
   \frac{u_{j}(\alpha_{j},\alpha_{j})}{(\gamma,\gamma)}.
\end{equation*}
Note that, by the definition of $S(\lambda,i)$, 
we have $m_{j}=0$ for all 
$j \in S(\lambda,i)$ with $j \ne i$. 
Consequently,
\begin{equation*}
\lambda(\gamma^{\vee})
  = \sum_{j \in S(\lambda,i)} m_{j}
   \frac{u_{j}(\alpha_{j},\alpha_{j})}{(\gamma,\gamma)}
  =m_{i}
   \frac{u_{i}(\alpha_{i},\alpha_{i})}{(\gamma,\gamma)}.
\end{equation*}
Since $u_{i} \in \BZ_{> 0}$ as seen above, 
we obtain from Lemma~\ref{lem:f-pair} that 
$u_{i}(\alpha_{i},\alpha_{i})/(\gamma,\gamma)
\in \BZ_{> 0}$. Thus we conclude that 
$\lambda(\gamma^{\vee}) \in m_{i}\BZ_{> 0}$, 
as desired. This proves the lemma. 
\end{proof}
%
%
\begin{lem} \label{lem:supp2}
Let $\zeta \in \prr(\lambda,i)$, 
with $\gamma:=\ba{\zeta} \in \pfr$ its finite part. 
Then, the finite part $\fp(r_{\zeta}(\lambda))$ of 
$r_{\zeta}(\lambda) \in W\lambda$ equals 
$\lambda(\gamma^{\vee})\gamma$. 
Hence the support of 
$\fp(r_{\zeta}(\lambda))$ is contained in 
$S(\lambda,i)$, and contains 
the $i \in \Supp_{\ge 2}(\lambda)$.
\end{lem}

\begin{proof}
By condition~(a) of Definition~\ref{dfn:prr-lam}, 
$\zeta$ is either of the form 
$\zeta=-\gamma+c_{\gamma}\delta$, or of the form 
$\zeta=\frac{1}{2}(-\gamma+\delta)$. 
It immediately follows from 
Lemma~\ref{lem:pair1} that 
\begin{equation*}
r_{\zeta}(\lambda)=
\begin{cases}
\lambda-\lambda(\gamma^{\vee})\gamma+
\lambda(\gamma^{\vee})c_{\gamma}\delta
  & \text{if $\zeta=-\gamma+c_{\gamma}\delta$}, \\[3mm]
\lambda-\lambda(\gamma^{\vee})\gamma+
\lambda(\gamma^{\vee})\delta
  & \text{if $\zeta=
    \frac{1}{2}(-\gamma+\delta)$}.
\end{cases}
\end{equation*}
Recall from the proof of Lemma~\ref{lem:pair3} 
that $\lambda(\gamma^{\vee}) \in m_{i}\BZ_{> 0}$. 
Hence the finite part $\fp(r_{\zeta}(\lambda))$ of 
$r_{\zeta}(\lambda) \in W\lambda$ equals 
$\lambda(\gamma^{\vee})\gamma$. 
Furthermore, 
it follows from condition~(b) of 
Definition~\ref{dfn:prr-lam} that the support of 
$\fp(r_{\zeta}(\lambda))=
\lambda(\gamma^{\vee})\gamma$ is 
contained in $S(\lambda,i)$, and 
contains the $i \in \Supp_{\ge 2}(\lambda)$.
This proves the lemma. 
\end{proof}
Now we take an element $\zeta_{i} \in \prr(\lambda,i)$
such that the finite part of $\zeta_{i}$ is maximal 
with respect to the (usual) partial ordering
defined by $\fin{Q}_{+}$. 
Namely, let $\zeta_{i}$ be an element of 
$\prr(\lambda,i)$ such that 
there exists no $\zeta \in \prr(\lambda,i)$ 
with $\fp(r_{\zeta}(\lambda)) - 
\fp(r_{\zeta_{i}}(\lambda)) \in \fin{Q}_{+} \setminus \{0\}$.
Let us denote by $\gamma_{i} \in \pfr$ 
the finite part $\ba{\zeta_{i}}$ of $\zeta_{i}$. 
%
%
\begin{lem} \label{lem:dist}
We have $\lambda(\zeta_{i}^{\vee}) \in 
m_{i}\BZ_{< 0}$, and hence 
$\lambda > r_{\zeta_{i}}(\lambda)$. 
Moreover, we have 
$\dist(\lambda,\,r_{\zeta_{i}}(\lambda))=1$. 
\end{lem}

\begin{proof}
It is obvious from Lemma~\ref{lem:pair3} that 
$\lambda(\zeta_{i}^{\vee}) \in 
m_{i}\BZ_{< 0}$, and hence 
$\lambda > r_{\zeta_{i}}(\lambda)$. 
We will prove that 
$\dist(\lambda,\,r_{\zeta_{i}}(\lambda))=1$. 
Assume that 
$\dist(\lambda,\,r_{\zeta_{i}}(\lambda))=k \ge 1$, 
and let $\lambda=\nu_{0} > \nu_{1} > \cdots > \nu_{k}=
r_{\zeta_{i}}(\lambda)$ be a longest chain for 
$(\lambda,\,r_{\zeta_{i}}(\lambda))$, with 
$\xi_{1},\,\xi_{2},\,\dots,\,\xi_{k}$
the corresponding positive real roots. 
Then we have $\nu_{l-1} > 
r_{\xi_{l}}(\nu_{l-1})=\nu_{l}$ and
$\dist (\nu_{l-1},\nu_{l})=1$ 
for each $1 \le l \le k$. 
Let $\beta_{l}$ denote
the finite part $\ba{\xi_{l}}$ of $\xi_{l} \in \prr$ 
for $1 \le l \le k$.
Since $\dist (\nu_{l-1},\nu_{l})=1$, 
we see from Lemma~\ref{lem:dist=1} that 
$\xi_{l}$ is one of the following forms: 
$\xi_{l}=\beta_{l}$, 
$\xi_{l}=-\beta_{l}+c_{\beta_{l}}\delta$, or 
$\xi_{l}=\frac{1}{2}(-\beta_{l}+\delta)$. 
We define subsets $\Pos$, $\Neg_{1}$, $\Neg_{2}$, and 
$\Neg$ of $\bigl\{1,\,2,\,\dots,\,k\bigr\}$ 
as in \eqref{eq:pos}, \eqref{eq:neg1}, 
\eqref{eq:neg2}, and \eqref{eq:neg}, respectively. 
It follows from Lemma~\ref{lem:fp} 
(see also Remark~\ref{rem:fp}) that 
%
%
\begin{equation} \label{eq:dist00}
\begin{cases}
\supp(\fp(\nu_{l-1}))=
\supp(\fp(\nu_{l})) \cup \supp(\beta_{l}) 
& \text{if $l \in \Pos$}, \\[1.5mm]
\supp(\fp(\nu_{l}))=
\supp(\fp(\nu_{l-1})) \cup \supp(\beta_{l})
& \text{if $l \in \Neg$}, 
\end{cases}
\end{equation}
and that
%
%
\begin{equation} \label{eq:dist001}
D(\nu_{l})=
\begin{cases}
D(\nu_{l-1}) & \text{if \, } l \in \Pos, \\[1mm]
D(\nu_{l-1})+\nu_{l-1}(\beta_{l}^{\vee})c_{\beta_{l}}
 & \text{if \, } l \in \Neg_{1}, \\[1mm]
D(\nu_{l-1})+\nu_{l-1}(\beta_{l}^{\vee})
 & \text{if \, } l \in \Neg_{2}.
\end{cases}
\end{equation}

Because $D(\nu_{k})=
D(r_{\zeta_{i}}(\lambda))=m_{i}d_{i}$ by 
condition (c) of Definition~\ref{dfn:prr-lam}, 
we have
%
%
\begin{equation} \label{eq:dist005}
m_{i}d_{i}=D(\nu_{k})
\stackrel{\text{by \eqref{eq:dist001}}}{=}
\sum_{l \in \Neg_{1}}
   \nu_{l-1}(\beta_{l}^{\vee})c_{\beta_{l}} + 
\sum_{l \in \Neg_{2}}
   \nu_{l-1}(\beta_{l}^{\vee}).
\end{equation}
For each $l \in \Neg$, we set 
%
%
\begin{equation} \label{eq:Dl}
D_{l}:=
\begin{cases}
\nu_{l-1}(\beta_{l}^{\vee})c_{\beta_{l}} 
  & \text{if $l \in \Neg_{1}$}, \\[1.5mm]
\nu_{l-1}(\beta_{l}^{\vee})
  & \text{if $l \in \Neg_{2}$}.
\end{cases}
\end{equation}
We see from Lemma~\ref{lem:pair1} that 
$D_{l} \in \BZ_{> 0}$ for all $l \in \Neg$, 
since $\nu_{l-1}(\xi_{l}^{\vee}) \in \BZ_{< 0}$ and 
$c_{\beta_{l}} \in \BZ_{> 0}$. Also, 
we see from \eqref{eq:dist005} that 
%
%
\begin{equation} \label{eq:dist01}
m_{i}d_{i}=\sum_{l \in \Neg} D_{l}. 
\end{equation}
%
%
\begin{claim} \label{c:dist02}
There exists $l \in \Neg$ such that 
the {\rm(}fixed\,{\rm)} 
$i \in \Supp_{\ge 2}(\lambda)$ is contained in $\supp(\beta_{l})$. 
\end{claim}

We know from Lemma~\ref{lem:supp2} that 
the support of the finite part $\fp(\nu_{k})$ of 
$\nu_{k}=r_{\zeta_{i}}(\lambda) \in W\lambda$ 
contains the $i \in \Supp_{\ge 2}(\lambda)$. 
Because $\supp(\fp(\nu_{0}))=\supp(\fp(\lambda))= 
\supp(0)=\emptyset$, it follows that 
there exists $1 \le l \le k$ 
such that $i \notin \supp(\fp(\nu_{l-1}))$ and 
$i \in \supp(\fp(\nu_{l}))$. 
Therefore, we see from 
\eqref{eq:dist00} that $l \in \Neg$, and 
$\supp(\fp(\nu_{l}))=
\supp(\fp(\nu_{l-1})) \cup \supp(\beta_{l})$. 
But, since $i \notin \supp(\fp(\nu_{l-1}))$ and 
$i \in \supp(\fp(\nu_{l}))$, we have
$i \in \supp(\beta_{l})$. This proves 
Claim~\ref{c:dist02}. 

\vsp\vsp

We set $l_{0}:=\min \bigl\{l \in \Neg \mid 
i \in \supp(\beta_{l})\bigr\}$.
%
%
\begin{claim} \label{c:dist03}
If $l \in \Pos$ and $l < l_{0}$, then 
$\supp(\beta_{l})$ does not contain 
the {\rm(}fixed\,{\rm)} $i \in \Supp_{\ge 2}(\lambda)$. 
\end{claim}

Suppose that there exists $l \in \Pos$ such that 
$l < l_{0}$ and $i \in \supp(\beta_{l})$. 
It follows from \eqref{eq:dist00} that 
$\supp(\fp(\nu_{l-1}))=
\supp(\fp(\nu_{l})) \cup \supp(\beta_{l})$, 
and hence $i \in \supp(\fp(\nu_{l-1}))$. 
Because $\supp(\fp(\nu_{0}))=\supp(\fp(\lambda))= 
\supp(0)=\emptyset$, 
there exists $1 \le l^{\prime} \le l-1$ 
such that 
$i \notin \supp(\fp(\nu_{l^{\prime}-1}))$ and 
$i \in \supp(\fp(\nu_{l^{\prime}}))$. 
Therefore, we see from 
\eqref{eq:dist00} that $l^{\prime} \in \Neg$, 
and 
$\supp(\fp(\nu_{l^{\prime}}))=
\supp(\fp(\nu_{l^{\prime}-1})) \cup 
\supp(\beta_{l^{\prime}})$. 
But, since $i \notin 
\supp(\fp(\nu_{l^{\prime}-1}))$ and 
$i \in \supp(\fp(\nu_{l^{\prime}}))$, 
we have $i \in \supp(\beta_{l^{\prime}})$, 
which contradicts the definition:
$l_{0}=\min \bigl\{l \in \Neg \mid 
i \in \supp(\beta_{l})\bigr\}$. 
This proves Claim~\ref{c:dist03}.

\vsp\vsp

We set 
$\beta_{l_{0}}^{\prime}:=
r_{\beta_{1}}r_{\beta_{2}} \cdots 
r_{\beta_{l_{0}-1}}(\beta_{l_{0}})$. 
%
%
\begin{claim} \label{c:dist04}
The support $\supp(\beta_{l_{0}}^{\prime})$ of 
$\beta_{l_{0}}^{\prime}$ contains the {\rm(}fixed\,{\rm)} 
$i \in \Supp_{\ge 2}(\lambda)$. 
\end{claim}

We set 
$\beta_{l_{0},l^{\prime}}:=
r_{\beta_{l^{\prime}}}
r_{\beta_{l^{\prime}+1}} \cdots 
r_{\beta_{l_{0}-1}}(\beta_{l_{0}})$
for $1 \le l^{\prime} \le l_{0}-1$, and 
$\beta_{l_{0},l_{0}}=\beta_{l_{0}}$; 
note that 
$\beta_{l_{0},1}=\beta_{l_{0}}^{\prime}$.
It suffices to show that $i \in 
\supp(\beta_{l_{0},l^{\prime}})$ 
for all $1 \le l^{\prime} \le l_{0}$ 
(put $l^{\prime}=1$). 
We show this 
by descending induction on $l^{\prime}$. 
If $l^{\prime}=l_{0}$, 
then it follows from the definition of $l_{0}$ that 
$i \in \supp(\beta_{l_{0},l_{0}})=
 \supp(\beta_{l_{0}})$. Assume that 
$l^{\prime} < l_{0}$. Then we have 
$\beta_{l_{0},l^{\prime}}=
r_{\beta_{l^{\prime}}}
 (\beta_{l_{0},l^{\prime}+1})=
\beta_{l_{0},l^{\prime}+1}-
\beta_{l_{0},l^{\prime}+1}(\beta_{l^{\prime}}^{\vee})
\beta_{l^{\prime}}$.
Because $i \notin \supp(\beta_{l^{\prime}})$ 
by Claim~\ref{c:dist03} and the definition of $l_{0}$, 
we deduce that the coefficient of $\alpha_{i}$ in 
$\beta_{l_{0},l^{\prime}}$ equals 
that in $\beta_{l_{0},l^{\prime}+1}$. 
This implies that 
$i \in \supp(\beta_{l_{0},l^{\prime}})$ 
since $i \in \supp(\beta_{l_{0},l^{\prime}+1})$ 
by the inductive assumption. 
Thus we have shown that $i \in 
\supp(\beta_{l_{0},l^{\prime}})$ 
for all $1 \le l^{\prime} \le l_{0}$, 
as desired. This proves Claim~\ref{c:dist04}.

\vsp\vsp

By an argument similar to that used to prove
Claim~\ref{c:p-2-0} in the proof of 
Lemma~\ref{lem:pair2}, we can prove 
the following claim.
%
%
%
\begin{claim} \label{c:dist035}
We have 
$(\beta_{l_{0}}^{\prime},\beta_{l_{0}}^{\prime})=
 (\beta_{l_{0}},\beta_{l_{0}})$, and 
$(\lambda,\beta_{l_{0}}^{\prime})=
(\nu_{l_{0}-1},\beta_{l_{0}})$. 
Moreover, 
$\nu_{l_{0}-1}(\beta_{l_{0}}^{\vee})=
 \lambda(\beta_{l_{0}}^{\prime\,\vee})$, and 
 $\beta_{l_{0}}^{\prime} \in \pfr$.
\end{claim}
%
%
%
%
\begin{claim} \label{c:dist045}
The following inequality holds\,{\rm:}
$D_{l_{0}} \ge 
 \sum_{j \in \supp(\beta_{l_{0}}^{\prime})} 
    m_{j}d_{j}$.
\end{claim}
%
Write $\beta_{l_{0}}^{\prime} \in \pfr$ as 
$\sum_{j \in I_{0}}
u_{j}\alpha_{j}$, with $u_{j} \in \BZ_{\ge 0}$ for 
$j \in I_{0}$. 

First assume that $\Fg$ is nontwisted 
(note that $l_{0} \in \Neg_{1}$, 
since $\Neg_{2}=\emptyset$).
Then we know from Remark~\ref{rem:cbeta}\,(1) and 
\eqref{eq:di} that $c_{\beta_{l_{0}}}=1$, and 
$d_{j}=1$ for all $j \in I_{0}$. 
Using Claim~\ref{c:dist035} and 
Lemma~\ref{lem:f-pair}, we deduce that 
\begin{align*}
D_{l_{0}} & =
\nu_{l_{0}-1}(\beta_{l_{0}}^{\vee})c_{\beta_{l_{0}}}= 
\nu_{l_{0}-1}(\beta_{l_{0}}^{\vee}) = 
\lambda(\beta_{l_{0}}^{\prime\,\vee})=
\sum_{j \in I_{0}}
 m_{j} 
 \frac{u_{j}(\alpha_{j},\alpha_{j})}
 {(\beta_{l_{0}}^{\prime},\beta_{l_{0}}^{\prime})} 
\\[3mm]
& =
\sum_{j \in  \supp(\beta_{l_{0}}^{\prime})}
 m_{j} 
 \underbrace{
 \frac{u_{j}(\alpha_{j},\alpha_{j})}
 {(\beta_{l_{0}}^{\prime},\beta_{l_{0}}^{\prime})}
 }_{\in \BZ_{> 0}}
\ge 
 \sum_{j \in \supp(\beta_{l_{0}}^{\prime})} m_{j}
=
 \sum_{j \in \supp(\beta_{l_{0}}^{\prime})} 
 m_{j}d_{j}. 
\end{align*}

Next assume that $\Fg$ is twisted, and 
$l_{0} \in \Neg_{1}$. 
Then we know from 
Remark~\ref{rem:cbeta}\,(2) and 
\eqref{eq:di} that 
$c_{\beta_{l_{0}}}=
(\beta_{l_{0}},\beta_{l_{0}})/2$, and that
\begin{equation*}
\frac{(\alpha_{j},\alpha_{j})}{2} = 
\begin{cases}
2d_{j} & 
 \text{if $\Fg$ is of type $A^{(2)}_{2\ell}$ 
       and $j=\ell$}, \\[3mm]
d_{j} & \text{otherwise}.
\end{cases}
\end{equation*}
In the same manner as
in the proof of Lemma~\ref{lem:pair2}, 
we deduce that
\begin{align*}
D_{l_{0}}
& =
 \nu_{l_{0}-1}(\beta_{l_{0}}^{\vee})c_{\beta_{l_{0}}}
=
 \sum_{j \in I_{0}} u_{j} m_{j}
 \frac{(\alpha_{j},\alpha_{j})}{2}
\\[3mm]
& =
  \sum_{j \in  \supp(\beta_{l_{0}}^{\prime})}
  u_{j} m_{j}
  \underbrace{\frac{(\alpha_{j},\alpha_{j})}{2}}_{%
  \in d_{j}\BZ_{> 0}} 
\ge 
  \sum_{j \in \supp(\beta_{l_{0}}^{\prime})} 
  m_{j}d_{j}. 
\end{align*}

Finally, assume that 
$\Fg$ is of type $A^{(2)}_{2\ell}$, and 
$l_{0} \in \Neg_{2}$.
Then we know from Remark~\ref{rem:cbeta}\,(2) and 
\eqref{eq:di} that $d_{j}=1$ for all $j \in I_{0}$. 
Using Claim~\ref{c:dist035} and 
Lemma~\ref{lem:f-pair}, we deduce that 
\begin{align*}
D_{l_{0}}
& = \nu_{l_{0}-1}(\beta_{l_{0}}^{\vee}) 
  = \lambda(\beta_{l_{0}}^{\prime\,\vee}) 
  = \sum_{j \in I_{0}}m_{j} 
   \frac{u_{j}(\alpha_{j},\alpha_{j})}
   {(\beta_{l_{0}}^{\prime},\beta_{l_{0}}^{\prime})}
\\[3mm]
& = \sum_{j \in \supp(\beta_{l_{0}}^{\prime})}m_{j} 
   \underbrace{
   \frac{u_{j}(\alpha_{j},\alpha_{j})}
   {(\beta_{l_{0}}^{\prime},\beta_{l_{0}}^{\prime})}
   }_{\in \BZ_{> 0}}
 \ge 
 \sum_{j \in \supp(\beta_{l_{0}}^{\prime})} m_{j}
 = \sum_{j \in \supp(\beta_{l_{0}}^{\prime})} 
 m_{j}d_{j}. 
\end{align*}
This proves Claim~\ref{c:dist045}.

\vsp\vsp

From \eqref{eq:dist01} and Claim~\ref{c:dist045}, 
we have 
\begin{align*}
m_{i}d_{i} & = 
 \sum_{l \in \Neg} D_{l}
   \ge D_{l_{0}} 
 \qquad \text{since $D_{l} \in \BZ_{> 0}$
 for all $l \in \Neg$} \\[3mm]
 & \ge 
 \sum_{j \in \supp(\beta_{l_{0}}^{\prime})} 
 m_{j}d_{j} \ge m_{i}d_{i}
 \qquad \text{since 
 $i \in \supp(\beta_{l_{0}}^{\prime})$ by 
 Claim~\ref{c:dist04}}.
\end{align*}
Therefore, we conclude that
\begin{equation*}
m_{i}d_{i} = 
 \sum_{l \in \Neg} D_{l} = 
 D_{l_{0}} =
 \sum_{j \in \supp(\beta_{l_{0}}^{\prime})} 
 m_{j}d_{j}, 
\end{equation*}
and hence that
%
%
\begin{equation} \label{eq:45-5}
\Neg=\bigl\{l_{0}\bigr\}, \quad 
D_{l_{0}}=m_{i}d_{i}, 
\quad \text{and} \quad
\supp(\beta_{l_{0}}^{\prime}) \cap 
\Supp(\lambda)=\bigl\{i\bigr\}.
\end{equation}
%
%
\begin{claim} \label{c:dist05}
We have $l_{0}=1$. Moreover, the support 
$\supp(\beta_{l_{0}})$ of 
$\beta_{l_{0}}=\beta_{1}$ is contained in 
$S(\lambda,i)$, and 
contains the {\rm(}fixed\,{\rm)} 
$i \in \Supp_{\ge 2}(\lambda)$.
\end{claim}

Suppose that $1 \in \Pos$. Then we have 
$\lambda(\xi_{1}^{\vee})=\lambda(\beta_{1}^{\vee})$. 
But, since $\lambda$ is a level-zero dominant integral weight, 
we have $\lambda(\beta_{1}^{\vee}) > 0$, and hence 
$\lambda(\xi_{1}^{\vee}) > 0$, which is a contradiction. 
Thus we have $1 \in \Neg$. 
Since $\Neg=\bigl\{l_{0}\bigr\}$, 
we conclude that $l_{0}=1$.
Note that $\beta_{l_{0}}^{\prime}=
\beta_{1}^{\prime}=\beta_{1}$ by the definition of 
$\beta_{l_{0}}^{\prime}$.
Hence it follows from \eqref{eq:45-5} that 
$\supp(\beta_{1}) \, \cap \, 
 \Supp(\lambda)=\bigl\{i\bigr\}$, i.e., that 
$\supp(\beta_{1})$ is contained in $S(\lambda,i)=
\bigl(I_{0} \setminus \Supp(\lambda)\bigr) 
\cup \{i\}$. 
This proves Claim~\ref{c:dist05}. 
%
%
\begin{claim} \label{c:dist06}
We have $\xi_{l_{0}}=\xi_{1} \in \prr(\lambda,i)$. 
\end{claim}
%
It is obvious that $\xi_{1}$ satisfies condition~(a) 
of Definition~\ref{dfn:prr-lam}, since $1 \in \Neg$. 
Also, it follows from Claim~\ref{c:dist05} that 
$\xi_{1}$ satisfies 
condition~(b) of Definition~\ref{dfn:prr-lam}.
Furthermore, from Lemma~\ref{lem:pair1} and 
\eqref{eq:Dl}, we see that 
$r_{\xi_{1}}(\lambda)=
 \lambda-\lambda(\beta_{1}^{\vee})
 \beta_{1} + D_{1} \delta$, with
 $\lambda(\beta_{1}^{\vee}) \in \BZ_{> 0}$.
But, it follows from \eqref{eq:45-5} and 
Claim~\ref{c:dist05} that $D_{1}=m_{i}d_{i}$. 
Thus, $\xi_{1}$ satisfies condition~(c) of 
Definition~\ref{dfn:prr-lam}. 
This proves Claim~\ref{c:dist06}.

\vsp\vsp

From Claim~\ref{c:dist05}, we have 
$\Neg=\bigl\{1\bigr\}$ and 
$\Pos=\bigl\{2,\,3,\,\dots,\,k\bigr\}$. 
Hence we deduce that
\begin{align*}
\nu_{k} & = 
 r_{\xi_{k}}r_{\xi_{k-1}} \cdots r_{\xi_{1}}(\lambda)
 = 
 \lambda-\sum_{l=1}^{k}
 \nu_{l-1}(\xi_{l}^{\vee}) \xi_{l}
 = 
 \lambda-\nu_{0}(\xi_{1}^{\vee})\xi_{1} -
     \sum_{l=2}^{k} 
     \nu_{l-1}(\xi_{l}^{\vee})\xi_{l} 
 \\[1.5mm]
 & =
 \lambda-\lambda(\xi_{1}^{\vee})\xi_{1} -
     \sum_{l=2}^{k} 
     \nu_{l-1}(\beta_{l}^{\vee})\beta_{l} 
 \quad \text{since $\nu_{0}=\lambda$, and $\Pos=
       \bigl\{2,\,3,\,\dots,\,k\bigr\}$,}
 \\[1.5mm]
 & = r_{\xi_{1}}(\lambda) - 
     \sum_{l=2}^{k}
     \nu_{l-1}(\beta_{l}^{\vee}) \beta_{l}. 
\end{align*}
Consequently, we have
%
%
\begin{equation} \label{eq:dist-fp}
\fp(r_{\zeta_{i}}(\lambda)) = 
\fp(\nu_{k})=\fp(r_{\xi_{1}}(\lambda)) + 
\sum_{l=2}^{k}
\nu_{l-1}(\beta_{l}^{\vee}) \beta_{l}.
\end{equation}
Since $\nu_{l-1}(\beta_{l}^{\vee})=
\nu_{l-1}(\xi_{l}^{\vee}) \in \BZ_{< 0}$ 
for all $l \in \Pos=\bigl\{2,\,3,\,\dots,\,k\bigr\}$, 
we see from \eqref{eq:dist-fp} that 
$\fp(r_{\xi_{1}}(\lambda)) - \fp(r_{\zeta_{i}}(\lambda)) = 
   -\sum_{l=2}^{k}
   \nu_{l-1}(\beta_{l}^{\vee}) \beta_{l} 
 \in \fin{Q}_{+}$.
Therefore, the definition of 
$\zeta_{i} \in \prr(\lambda,i)$ implies that 
$\fp(r_{\xi_{1}}(\lambda)) - \fp(r_{\zeta_{i}}(\lambda)) = 0$, 
and hence that $\Pos=\emptyset$, $k=1$. 
Thus we have proved $\dist(\lambda, r_{\zeta_{i}}(\lambda))=1$. 
This completes the proof of Lemma~\ref{lem:dist}.
\end{proof}

Now we set $\mu:=r_{\zeta_{i}}(\lambda)$. 
Then we see from condition~(c) of 
Definition~\ref{dfn:prr-lam}, 
using Lemma~\ref{lem:pair1}, that 
$\mu=r_{\gamma_{i}}(\lambda+m_{i}d_{i}\delta)$; note that 
$r_{\gamma_{i}} \in \fin{W}$. 
Let $w \in \fin{W}$ be the shortest element such that 
$w(\lambda+m_{i}d_{i}\delta)=\mu$, and let 
$w=r_{j_{1}}r_{j_{2}} \cdots r_{j_{K}}$ be 
a reduced expression, 
where $j_{1},\,j_{2},\,\dots,\,j_{K} \in I_{0}$. 
We set $\mu_{0}=\mu$, and 
$\mu_{L}:=r_{j_{L}}r_{j_{L-1}} \cdots r_{j_{1}}(\mu)$
for $1 \le L \le K$; 
note that $\mu_{K}=w^{-1}(\mu)=
\lambda+m_{i}d_{i}\delta$. 
Then, using \cite[Lemma~3.10]{Kac}, 
we can show that 
$\mu_{L-1}(h_{j_{L}}) \in \BZ_{< 0}$ 
for all $1 \le L \le K$. 
Hence we have 
$\mu=\mu_{0} > \mu_{1} > \mu_{2} 
> \cdots > \mu_{K}=\lambda+m_{i}d_{i}\delta$. 

Moreover, it follows that 
$\mu_{K}=w^{-1}(\mu)=
 r_{j_{K}} \cdots r_{j_{2}}r_{j_{1}}(\mu)=
 \mu-\sum_{L=1}^{K} 
 \mu_{L-1}(h_{j_{L}}) \alpha_{j_{L}}$. 
Also, we see from 
condition~(c) of Definition~\ref{dfn:prr-lam} that 
$\mu=r_{\zeta_{i}}(\lambda)=\lambda-
\fp(r_{\zeta_{i}}(\lambda))+m_{i}d_{i}\delta$. 
Combining these equations, we obtain that
\begin{equation*}
\mu_{K}=
\lambda-\biggl(
  \fp(r_{\zeta_{i}}\lambda)+\sum_{L=1}^{K} 
 \underbrace{\mu_{L-1}(h_{j_{L}})}_{
 \in \BZ_{< 0}} \alpha_{j_{L}}\biggr) 
 + m_{i}d_{i}\delta.
\end{equation*}
Since $\mu_{K}=\lambda+m_{i}d_{i}\delta$, 
we have
\begin{equation*}
 \fp(r_{\zeta_{i}}\lambda)+\sum_{L=1}^{K} 
 \underbrace{\mu_{L-1}(h_{j_{L}})}_{
 \in \BZ_{< 0}} \alpha_{j_{L}}=0.
\end{equation*}
Recall from Lemma~\ref{lem:supp2} that 
the support of the finite part 
$\fp(r_{\zeta_{i}}(\lambda))$ of 
$r_{\zeta_{i}}(\lambda) \in W\lambda$ is 
contained in $S(\lambda,i)$. 
Therefore, we conclude that 
$j_{1},\,j_{2},\,\dots,\,j_{K} \in 
S(\lambda,i)$. 
%
%
\begin{lem} \label{lem:mi-ch}
{\rm (1) }
We have $\dist(\mu_{L-1},\mu_{L})=1$ 
for all $1 \le L \le K$. 

\vsp

\noindent
{\rm (2) } 
We have $\mu_{L-1}(h_{j_{L}}) \in m_{i}\BZ_{< 0}$ 
for all $1 \le L \le K$.
\end{lem}

\begin{proof}
(1) Since the positive real root corresponding to 
$(\mu_{L-1},\mu_{L})$ is 
the simple root $\alpha_{j_{L}} \in \fin{\Pi}$, 
the assertion immediately follows 
from Remark~\ref{rem:Q+}\,(2). 

\vsp 

\noindent
(2) Because $\mu_{L-1}(h_{j_{L}}) \in \BZ_{< 0}$
for all $1 \le L \le K$ as shown 
in the argument above this lemma, 
we need only show that 
$\mu_{L-1}(h_{j_{L}}) \in m_{i}\BZ$ 
for all $1 \le L \le K$. 
We will show this by induction on $L$.

Assume first that $L=1$. 
From condition~(c) of Definition~\ref{dfn:prr-lam} 
and Lemma~\ref{lem:supp2}, we deduce that
$\mu_{0}=\mu=r_{\zeta_{i}}(\lambda)=\lambda
  -\lambda(\gamma_{i}^{\vee})\gamma_{i}
  +m_{i}d_{i}\delta$.
Since $j_{1} \in S(\lambda,i)=
\bigl(I_{0} \setminus \Supp(\lambda)\bigr) \cup \{i\}$ 
as shown just above this lemma, 
it follows that $\lambda(h_{j_{1}}) = m_{i}$ if $j_{1} = i$, 
and $\lambda(h_{j_{1}}) = 0$ if $j_{1} \ne i$. 
In addition, we see from the proof of 
Lemma~\ref{lem:pair3} that 
$\lambda(\gamma_{i}^{\vee}) \in 
m_{i}\BZ_{> 0}$. Therefore, we obtain that
\begin{equation*}
\mu_{0}(h_{j_{1}})= 
\underbrace{\lambda(h_{j_{1}})}_{0 \text{ or } m_{i}}-
\underbrace{\lambda(\gamma_{i}^{\vee})}_{\in m_{i}\BZ_{> 0}}
\underbrace{\gamma_{i}(h_{j_{1}})}_{\in \BZ}+ 
m_{i}d_{i}\underbrace{\delta(h_{j_{1}})}_{=0} 
\in m_{i}\BZ.
\end{equation*}

Assume next that $L > 1$. Because
$\mu_{L-1}=
  \mu-
  \sum_{L^{\prime}=1}^{L-1} 
     \mu_{L^{\prime}-1}(h_{j_{L^{\prime}}})
     \alpha_{j_{L^{\prime}}}$, 
it follows that 
$\mu_{L-1}(h_{j_{L}})=
  \mu(h_{j_{L}})-
  \sum_{L^{\prime}=1}^{L-1} 
     \mu_{L^{\prime}-1}(h_{j_{L^{\prime}}})
     \alpha_{j_{L^{\prime}}}(h_{j_{L}})$.
Here, in the same manner as above, 
we see that
\begin{equation*}
\mu(h_{j_{L}})= 
\underbrace{\lambda(h_{j_{L}})}_{0 \text{ or } m_{i}}-
\underbrace{\lambda(\gamma_{i}^{\vee})}_{\in m_{i}\BZ_{> 0}}
\underbrace{\gamma_{i}(h_{j_{L}})}_{\in \BZ}+ 
m_{i}d_{i}\underbrace{\delta(h_{j_{L}})}_{=0} 
\in m_{i}\BZ. 
\end{equation*}
Also, it follows from the inductive assumption that 
$\mu_{L^{\prime}-1}(h_{j_{L^{\prime}}}) \in m_{i}\BZ$ 
for all $1 \le L^{\prime} \le L-1$. 
Since $\alpha_{j_{L^{\prime}}}(h_{j_{L}}) \in \BZ$ 
for all $1 \le L^{\prime} \le L-1$, we conclude that 
\begin{equation*}
\mu_{L-1}(h_{j_{L}})=
  \underbrace{\mu(h_{j_{L}})}_{\in m_{i}\BZ}-
  \sum_{L^{\prime}=1}^{L-1} 
   \underbrace{\mu_{L^{\prime}-1}(h_{j_{L^{\prime}}})}_{\in m_{i}\BZ}
   \underbrace{\alpha_{j_{L^{\prime}}}(h_{j_{L}})}_{\in \BZ}
   \in m_{i}\BZ.
\end{equation*}
Thus we have proved the lemma.
\end{proof}
%
%
\begin{prop} \label{prop:mi-ch}
For each $1 \le q^{\prime} \le m_{i}-1$, 
there exists a $(q^{\prime}/m_{i})$-chain 
for $(\lambda, \lambda+m_{i}d_{i}\delta)$. 
\end{prop}

\begin{proof}
Combining Lemmas~\ref{lem:dist} and \ref{lem:mi-ch}, 
we see that 
$\lambda > r_{\zeta_{i}}(\lambda)=\mu_{0} > 
  \mu_{1} > \mu_{2} > \cdots > \mu_{K}=
  \lambda+m_{i}d_{i}\delta$
is a $(q^{\prime}/m_{i})$-chain for 
$(\lambda,\,\lambda+m_{i}d_{i}\delta)$, with 
$\zeta_{i},\,
 \alpha_{j_{1}},\,\alpha_{j_{2}},\,\dots,\,\alpha_{j_{K}}$
the corresponding positive real roots.
This proves the proposition. 
\end{proof}

Now we are ready to prove Proposition~\ref{prop:suf}.

\begin{proof}[Proof of Proposition~\ref{prop:suf}]
We will prove this proposition by induction on $N$. 
If $N=0$, then the assertion is obvious. 
Assume that $N > 0$. We write $N$ as 
$\sum_{j \in I_{0}(\lambda,p)}m_{j}d_{j}k_{j}$ 
with $k_{j} \in \BZ_{\ge 0}$ for $j \in I_{0}(\lambda,p)$, and 
take $i \in I_{0}(\lambda,p)$ such that $k_{i} > 0$, 
$m_{i} > 0$ 
(and hence $m_{i} \ge 2$, i.e., $i \in \Supp_{\ge 2}(\lambda)$).
Note that we have $q/p = q^{\prime}/m_{i}$ 
for some $1 \le q^{\prime} \le m_{i}-1$, since 
$m_{i} \in p\BZ$. We know from 
Proposition~\ref{prop:mi-ch} that 
there exists a $(q/p)$-chain 
for $(\lambda, \lambda+m_{i}d_{i}\delta)$; 
we write it as: 
$\lambda=\nu_{0} > \nu_{1} > \dots > \nu_{K}=
\lambda+m_{i}d_{i}\delta$, 
with $\xi_{1},\,\xi_{2},\,\dots,\,\xi_{K}$ 
the corresponding positive real roots.
Since $m_{i} > 0$ and $k_{i} > 0$, it follows that 
$N > N-m_{i}d_{i} \in \sum_{j \in I_{0}(\lambda,p)} 
 m_{j}d_{j}\BZ_{\ge 0}$. 
Therefore, 
by the inductive assumption, 
there exists a $(q/p)$-chain for 
$(\lambda, \lambda+(N-m_{i}d_{i})\delta)$. 
Hence we deduce from Lemma~\ref{lem:d-shift} that 
there exists a $(q/p)$-chain for 
$(\lambda+m_{i}d_{i}\delta, \, \lambda+N\delta)$; 
we write it as: $\lambda+m_{i}d_{i}\delta=
\nu_{0}^{\prime} > \nu_{1}^{\prime} > \dots > 
\nu_{K^{\prime}}^{\prime}=\lambda+N\delta$, with 
$\xi_{1}^{\prime},\,\xi_{2}^{\prime},\,\dots,\,
 \xi_{K^{\prime}}^{\prime}$ 
the corresponding positive real roots.
From these two $(q/p)$-chains, we can construct 
a $(q/p)$-chain for $(\lambda,\,\lambda+N\delta)$
as follows: 
\begin{equation*}
\lambda=\nu_{0} > \nu_{1} > \dots > \nu_{K}=
\lambda+m_{i}d_{i}\delta=
\nu_{0}^{\prime} > \nu_{1}^{\prime} > \dots > 
\nu_{K^{\prime}}^{\prime}=\lambda+N\delta, 
\end{equation*}
where $\xi_{1},\,\xi_{2},\,\dots,\,
 \xi_{K},\,
 \xi_{1}^{\prime},\,\xi_{2}^{\prime},\,\dots,\,
 \xi_{K^{\prime}}^{\prime}$ are the corresponding positive real roots.
Thus we have completed the proof of 
Proposition~\ref{prop:suf}. 
\end{proof}

%
\subsection{Proof of Theorem~\ref{thm:comps}.}
\label{subsec:proof}

It follows from 
Corollary~\ref{cor:suf} 
that a path $\pi \in \BP$ 
having an expression of the form \eqref{eq:comps} is 
indeed an LS path of shape $\lambda$ 
(see Definition~\ref{dfn:LS} of an LS path).
Furthermore, we see from Lemma~\ref{lem:d-shift-ext} that 
such a path $\pi \in \BP$ is extremal. 
This proves part (1) of Theorem~\ref{thm:comps}. 

Let us prove part~(2) of Theorem~\ref{thm:comps}. 
Let $\pi \in \BB(\lambda)$ be an LS path of shape $\lambda$. 
We know from 
Proposition~\ref{thm:NSz}\,(1) 
that the $P_{\cl}$-crystal $\BB(\lambda)_{\cl}$ is connected. 
Therefore, there exists some monomial $X$ 
in the root operators 
$e_{j}$, $f_{j}$ for $j \in I$ such that 
$X \cl(\pi)=\cl(\pi_{\lambda})$. 
Hence it follows from \eqref{eq:cl-ro} that 
$\cl(X\pi)=X\cl(\pi)=\cl(\pi_{\lambda})$. 
From the fact that $\cl(X\pi)=\cl(\pi_{\lambda})$, 
we deduce, using Remark~\ref{rem:exp-cl}, 
that the path $\pi^{\prime}:=X\pi \in \BP$ 
has an expression of the form:
\begin{equation*}
\pi^{\prime} 
=
(\lambda-N_{1}^{\prime}\delta,\,\dots,\,
 \lambda-N_{s^{\prime}-1}^{\prime}\delta,\ 
 \lambda-N_{s^{\prime}}^{\prime}\delta \ ; \ 
 \sigma_{0},\,\sigma_{1},\,\dots,\,\sigma_{s^{\prime}-1},\,
 \sigma_{s^{\prime}}),
\end{equation*}
with $N_{1}^{\prime},\,N_{2}^{\prime},\,\dots,\,
N_{s^{\prime}}^{\prime} \in \BQ$ and 
$0=\sigma_{0} < \sigma_{1} < \cdots < 
 \sigma_{s^{\prime}-1} < \sigma_{s^{\prime}}=1$. 
Here we note that the path $\pi^{\prime}=X\pi \in \BP$ is 
an LS path of shape $\lambda$. 
Therefore, it follows 
from Remark~\ref{rem:LS}\,(1) that 
$\lambda-N_{u^{\prime}}^{\prime}\delta \in W\lambda$ 
for all $1 \le u^{\prime} \le s^{\prime}$, 
which along with Lemma~\ref{lem:Worbit} 
implies that $N_{u^{\prime}}^{\prime} \in \BZ$
for all $1 \le u^{\prime} \le s^{\prime}$. 
Now, if we define a path $\psi \in \BP$ by:
$\psi=\pi_{\lambda}-\pi^{\prime}$, then 
$\psi \in \BP$ has an expression of the form:
\begin{equation*}
\psi=
(N_{1}^{\prime}\delta,\,\dots,\,
 N_{s^{\prime}-1}^{\prime}\delta,\ 
 N_{s^{\prime}}^{\prime}\delta \ ; \ 
 \sigma_{0},\,
 \sigma_{1},\,\dots,\,
 \sigma_{s^{\prime}-1},\,
 \sigma_{s^{\prime}}).
\end{equation*}
Since $\lambda-N_{s^{\prime}}^{\prime}\delta \in W\lambda$ 
by Remark~\ref{rem:LS}\,(1), it follows that 
there exists some $w \in W$ such that 
$w\lambda=\lambda+N_{s^{\prime}}^{\prime}\delta$. 
Then, from \eqref{eq:weyl}, we have 
$S_{w}\pi_{\lambda}=\pi_{w\lambda}=
\pi_{\lambda+N_{s^{\prime}}^{\prime}\delta}=
\pi_{\lambda}+\pi_{N_{s^{\prime}}^{\prime}\delta}$.
Therefore, by using Lemma~\ref{lem:d-shift3}, 
we deduce that 
\begin{equation*}
S_{w}\pi^{\prime}=S_{w}(\pi_{\lambda}-\psi)=
S_{w}\pi_{\lambda}-\psi=
\pi_{\lambda}+\pi_{N_{s^{\prime}}^{\prime}\delta} - \psi = 
\pi^{\prime}+\pi_{N_{s^{\prime}}^{\prime}\delta}, 
\end{equation*}
and hence that the path 
$\pi^{\prime\prime}:=
S_{w}\pi^{\prime} \in \BB(\lambda)$ 
has an expression of the form:
%
%
\begin{equation} \label{eq:comps1}
\pi^{\prime\prime}=(
 \lambda-N_{1}^{\prime\prime}\delta,\,\dots,\,
 \lambda-N_{s^{\prime}-1}^{\prime\prime}\delta,\ 
 \lambda \ ; \ 
 \sigma_{0},\,
 \sigma_{1},\,\dots,\,
 \sigma_{s^{\prime}-1},\,
 \sigma_{s^{\prime}}),
\end{equation}
where we set $N_{u^{\prime}}^{\prime\prime}:=
N_{u^{\prime}}^{\prime}-N_{s^{\prime}}^{\prime}$ 
for $1 \le u^{\prime} \le s^{\prime}$; 
note that $N_{s^{\prime}}^{\prime\prime}=0$. 
We will show that 
$\pi^{\prime\prime} \in \BB(\lambda)$ has an expression of 
the form \eqref{eq:comps}.

Assume now that there exists 
some $1 \le u^{\prime} \le s^{\prime}-1$ such that 
$\sigma_{u^{\prime}} \notin \turn(\lambda)$, 
and write $\sigma_{u^{\prime}}=q/p$, 
where $1 \le q < p$ are coprime integers. 
Then it can easily be seen 
from the definition \eqref{eq:turn} of 
$\turn(\lambda)$ that for each $i \in I_{0}$, 
we have $m_{i} \notin p\BZ$ unless $m_{i}=0$. 
This means that there exists no $j \in I_{0}(\lambda,p)$
such that $m_{j}> 0$. Hence it follows from 
Corollary~\ref{cor:nec} that $N_{u^{\prime}}^{\prime\prime}-
N_{u^{\prime}+1}^{\prime\prime} \in \sum_{j \in I_{0}(\lambda,p)}
m_{j}d_{j}\BZ_{\ge 0}=\bigl\{0\bigr\}$, i.e., that 
$N_{u^{\prime}}^{\prime\prime}=N_{u^{\prime}+1}^{\prime\prime}$. 
As a result, we have 
$\lambda-N_{u^{\prime}}^{\prime\prime}\delta=
 \lambda-N_{u^{\prime}+1}^{\prime\prime}\delta$. 
Thus we can ``omit'' 
$\lambda-N_{u^{\prime}}^{\prime\prime}\delta$ and 
$\sigma_{u^{\prime}}$ from the expression 
\eqref{eq:comps1} (see Remark~\ref{rem:expr}\,(1)). 
Namely, the path $\pi^{\prime\prime} \in \BB(\lambda)$ has 
an expression of the form:
%
%
\begin{equation} \label{eq:comps2}
\pi^{\prime\prime}=
(\lambda-N_{s_{1}^{\prime}}^{\prime\prime}\delta,\,
 \dots,\,
 \lambda-N_{s_{L-1}^{\prime}}^{\prime\prime}\delta,\ 
 \lambda \ ; \ 
 \sigma_{s_{0}^{\prime}},\,
 \sigma_{s_{1}^{\prime}},\,\dots,\,
 \sigma_{s_{L-1}^{\prime}},\,
 \sigma_{s_{L}^{\prime}}),
\end{equation}
where $0=s_{0}^{\prime} < s_{1}^{\prime} 
< \cdots < s_{L}^{\prime}=s^{\prime}$, and 
$\sigma_{s_{l}^{\prime}} \in \turn(\lambda)=
\bigl\{\tau_{1} < \tau_{2} < \cdots < \tau_{s-1}\bigr\}$ 
for every $1 \le l \le L-1$. 
We define an increasing sequence 
$0=s_{0} < s_{1} < \cdots < s_{L}=s$ so that 
$\tau_{s_{l}}=\sigma_{s_{l}^{\prime}}$ 
for every $0 \le l \le L$, and set 
$N_{s_{l}}:=N_{s_{l}^{\prime}}^{\prime\prime}$ 
for $1 \le l \le L-1$. 
Then the expression 
\eqref{eq:comps2} of 
$\pi^{\prime\prime} \in \BB(\lambda)$ 
can be written as: 
%
%
\begin{equation} \label{eq:comps3}
\pi^{\prime\prime}=
(\lambda-N_{s_{1}}\delta,\,\dots,\,
 \lambda-N_{s_{L-1}}\delta,\ 
 \lambda \ ; \ 
 \tau_{s_{0}},\,
 \tau_{s_{1}},\,\dots,\,
 \tau_{s_{L-1}},\,
 \tau_{s_{L}}).
\end{equation}
Assume next that there exists some $1 \le u \le s-1$ 
such that 
$\tau_{u} \in \turn(\lambda) \setminus 
\bigl\{\tau_{s_{1}},\,\dots,\,\tau_{s_{L-1}}\bigr\}$, 
and let $1 \le l_{0} \le L$ be such that 
$\tau_{s_{l_{0}-1}} < \tau_{u} < \tau_{s_{l_{0}}}$. 
Then we can ``insert'' $\tau_{u}$ 
between $\tau_{s_{l_{0}-1}}$ and $\tau_{s_{l_{0}}}$ 
in the expression \eqref{eq:comps3} of 
$\pi^{\prime\prime}$ 
(see Remark~\ref{rem:expr}\,(2)). 
Namely, the path $\pi^{\prime\prime} \in \BB(\lambda)$ 
has an expression of the form:
\begin{align*}
& \pi^{\prime\prime}=
(\lambda-N_{s_{1}}\delta,\,\dots,\,
 \lambda-N_{s_{l_{0}-1}}\delta, \ 
 \underbrace{\lambda-N_{u}\delta}_{\text{inserted}}, \ 
 \lambda-N_{s_{l_{0}}}\delta,\,\dots,\,
 \lambda-N_{s_{L-1}}\delta,\ 
 \lambda \ ;
\\
& \hspace{30mm} \tau_{s_{0}},\,\tau_{s_{1}},\,\dots,\,
 \tau_{s_{l_{0}-1}},\,
 \underbrace{\tau_{u}}_{\text{inserted}},\,
 \tau_{s_{l_{0}}},\,\dots,\,\tau_{s_{L-1}},\,\tau_{s_{L}}),
\end{align*}
where we set $N_{u}:=N_{s_{l_{0}}}$. 
Repeating these two procedures if necessary, 
we finally obtain an expression of 
$\pi^{\prime\prime} \in \BP$ of the form:
%
%
\begin{equation}
\pi^{\prime\prime}=
(\lambda-N_{1}\delta,\,\dots,\,
 \lambda-N_{s-1}\delta,\ 
 \lambda \ ; \ 
 \tau_{0},\,\tau_{1},\,\dots,\,\tau_{s-1},\,\tau_{s}),
\end{equation}
where $N_{1},\,N_{2},\,\dots,\,N_{s-1} \in \BZ$. 
Since $\pi^{\prime\prime} \in \BB(\lambda)$ is 
an LS path of shape $\lambda$, 
there exists a $\tau_{u}$-chain for $(\lambda-N_{u}\delta,\,
\lambda-N_{u+1}\delta)$ for each $1 \le u \le s-1$, where we set 
$N_{s}:=0$. Hence it follows from Corollary~\ref{cor:nec} that 
\begin{equation*}
  N_{u}-N_{u+1} \in 
   \sum_{j \in I_{0}(\lambda,p_{u})} 
   m_{j}d_{j} \BZ_{\ge 0}
  \quad \text{for each $1 \le u \le s-1$}. 
\end{equation*}
Thus we have proved that every $\pi \in \BB(\lambda)$ is 
connected to an LS path of shape $\lambda$ 
having an expression 
of the form \eqref{eq:comps}, i.e., that 
each connected component of $\BB(\lambda)$ 
contains an LS path of shape $\lambda$ 
having an expression of the form \eqref{eq:comps}.

Now, let $\pi_{1} \in \BB(\lambda)$ and 
$\pi_{2} \in \BB(\lambda)$ be 
LS paths having an expression of 
the form \eqref{eq:comps}:
\begin{align*}
& \pi_{1}=
(\lambda-N_{1,1}\delta,\,\dots,\,
 \lambda-N_{1,s-1}\delta,\ 
 \lambda \ ; \ 
 \tau_{0},\,\tau_{1},\,\dots,\,
 \tau_{s-1},\,\tau_{s}), \\[1.5mm]
& \pi_{2}=
(\lambda-N_{2,1}\delta,\,
 \dots,\,
 \lambda-N_{2,s-1}\delta,\ 
 \lambda \,;\, 
 \tau_{0},\,\tau_{1},\,\dots,\,
 \tau_{s-1},\,\tau_{s}), 
\end{align*}
and assume that the paths 
$\pi_{1}$ and $\pi_{2}$ are in the same 
connected component of $\BB(\lambda)$. 
Then, there exists some monomial $X^{\prime}$ 
in the root operators $e_{j}$, $f_{j}$ for $j \in I$ such that 
$X^{\prime} \pi_{1} = \pi_{2}$. 
We define a path $\psi_{1} \in \BP$ by:
\begin{equation*}
\psi_{1}=
(N_{1,1}\delta,\,\dots,\,
 N_{1,s-1}\delta,\ 0 \,;\, 
 \tau_{0},\,\tau_{1},\,\dots,\,\tau_{s-1},\,\tau_{s}). 
\end{equation*}
Note that $\pi_{1}=\pi_{\lambda}-\psi_{1}$, and hence 
that $\pi_{\lambda}=\pi_{1}+\psi_{1}$. 
We deduce from Lemma~\ref{lem:d-shift2} that 
\begin{align*}
X^{\prime} \pi_{\lambda} & = 
X^{\prime} (\pi_{1}+\psi_{1}) = 
 \pi_{2}+\psi_{1}
 \\[1.5mm]
& =
(\lambda-(N_{2,1}-N_{1,1})\delta,\,\dots,\,
 \lambda-(N_{2,s-1}-N_{1,s-1})\delta,\ \lambda \ ; \ 
 \tau_{0},\,\tau_{1},\,\dots,\,\tau_{s-1},\,\tau_{s}).
\end{align*}
Because $X^{\prime}\pi_{\lambda} \in \BB(\lambda)$ is 
an LS path of shape $\lambda$, we see from 
Remark~\ref{rem:LS}\,(1) that 
\begin{equation*}
\lambda-(N_{2,1}-N_{1,1})\delta \ge 
\dots \ge 
\lambda-(N_{2,s-1}-N_{1,s-1})\delta \ge 
\lambda.
\end{equation*}
Hence it follows from Remark~\ref{rem:Q+}\,(1) that 
$N_{2,u}-N_{1,u} \ge 0$ for all $1 \le u \le s$. 
By interchanging the roles of $\pi_{1}$ and $\pi_{2}$, 
we also deduce that 
$N_{1,u}-N_{2,u} \ge 0$ for all $1 \le u \le s$. 
Thus we conclude that $N_{1,u}=N_{2,u}$ 
for all $1 \le u \le s$, and hence $\pi_{1}=\pi_{2}$. 
This proves part (2) of Theorem~\ref{thm:comps}. 

Finally, let us prove part~(3) of Theorem~\ref{thm:comps}. 
Let $\pi \in \BB(\lambda)$ be 
an extremal LS path of shape $\lambda$. 
Then we know from Remark~\ref{rem:cl-ext} that 
$\cl(\pi) \in \BB(\lambda)_{\cl}$ is extremal. 
Hence it follows from Theorem~\ref{thm:NSz}\,(3) that 
there exists some $w \in W$ such that 
$\cl(S_{w}\pi)=S_{w}\cl(\pi)=
\cl(\pi_{\lambda})$. 
Here we recall that 
in the proof of part (2) above, we proved that 
an LS path $\pi^{\prime} \in \BB(\lambda)$ such that 
$\cl(\pi^{\prime})=\cl(\pi_{\lambda})$ is 
contained in the $W$-orbit of an extremal LS path 
having an expression of the form \eqref{eq:comps}. 
Therefore, we conclude that 
$S_{w}\pi$ (and hence $\pi$) is contained in 
the $W$-orbit of an extremal LS path 
having an expression of the form \eqref{eq:comps}.
This proves part (3) of Theorem~\ref{thm:comps}. 
Thus we have completed the proof of Theorem~\ref{thm:comps}. 
%
%
%

\section{Crystal structure of 
the connected component $\BB_{0}(\lambda)$ of $\BB(\lambda)$.}
\label{sec:aff}
%
%
%
\subsection{Affinization of the 
$P_{\cl}$-crystal $\BB(\lambda)_{\cl}$.}
\label{subsec:aff}

We take (and fix) a $\BQ$-linear embedding 
$\af: \Fh^{\ast}/\BQ\delta \hookrightarrow \Fh^{\ast}$ 
such that $\cl \circ \af = \id$, 
$\af(\cl(\alpha_{j}))=\alpha_{j}$ 
for all $j \in I_{0}$, and such that 
$\af(P_{\cl}) \subset P$ 
(note that $\Fh^{\ast} \simeq \BQ \otimes_{\BZ} P$ and 
$\Fh^{\ast}/\BQ\delta \simeq \BQ \otimes_{\BZ} P_{\cl}$). 
Since $\delta-a_{0}\alpha_{0}=\theta \in \fin{Q}_{+}$, 
we see that 
$\af(\cl(\alpha_{0}))=
\af(\cl(a_{0}^{-1}(\delta-\theta)))=
-a_{0}^{-1}\theta$.
In addition, since the level-zero fundamental weight $\vpi_{i}$ 
is contained in $\bigoplus_{j \in I_{0}} \BQ \alpha_{j}$ 
for all $i \in I_{0}$, it follows that
$\af(\cl(\vpi_{i}))=\vpi_{i}$ for all $i \in I_{0}$.

Let $\lambda \in P$ 
be the (fixed) level-zero dominant integral weight 
of the form \eqref{eq:lambda}. We define a positive integer 
$d_{\lambda} \in \BZ_{> 0}$ by: 
\begin{equation}
W\lambda \cap (\lambda+\BZ\delta) = 
 \lambda+d_{\lambda}\BZ\delta.
\end{equation}
%
%
\begin{rem} \label{rem:dlam}
It follows from Lemma~\ref{lem:dlam} 
(see the proof of Lemma~\ref{lem:Worbit}) that 
$d_{\lambda}$ is equal to 
the greatest common divisor of the integers
  $\bigl\{m_{i}d_{i}\bigr\}_{i \in I_{0}}$. 
\end{rem}

Define the affinization $\ha{\BB(\lambda)_{\cl}}$ of 
the $P_{\cl}$-crystal $\BB(\lambda)_{\cl}$ to be 
the direct product set 
\begin{equation}
\ha{\BB(\lambda)_{\cl}}:=
 \BB(\lambda)_{\cl} \times (a_{0}^{-1}\BZ)
\end{equation}
of $\BB(\lambda)_{\cl}$ and $a_{0}^{-1}\BZ$; 
let us denote an element 
$(\eta,n) \in \ha{\BB(\lambda)_{\cl}}$ by 
$\eta \otimes z^{n}$.
Now we give a $P$-crystal structure to 
$\ha{\BB(\lambda)_{\cl}}$ as follows. 
We define the Kashiwara operators 
$e_{j}$, $f_{j}$  on $\ha{\BB(\lambda)_{\cl}}$ 
for $j \in I$ by:
\begin{align}
& e_{j}(\eta \otimes z^{n})=
  (e_{j}\eta) \otimes z^{n}, 
  \quad 
  f_{j}(\eta \otimes z^{n})=
  (f_{j}\eta) \otimes z^{n} 
  \quad
  \text{for $j \in I_{0}$}, 
\\[3mm]
& e_{0}(\eta \otimes z^{n})=
  (e_{0}\eta) \otimes z^{n+a_{0}^{-1}}, 
  \quad 
  f_{0}(\eta \otimes z^{n})=
  (f_{0}\eta) \otimes z^{n-a_{0}^{-1}}, 
\end{align}
for $\eta \otimes z^{n} \in \ha{\BB(\lambda)_{\cl}}$. 
Here $\bzero \otimes z^{n}$ is understood to be 
$\bzero$ for each $n \in a_{0}^{-1}\BZ$. 
We define the weight map 
$\wt:\ha{\BB(\lambda)_{\cl}} \rightarrow P$ by:
%
%
\begin{equation} \label{eq:aff-wt}
\wt(\eta \otimes z^{n})=\af(\eta(1)) + n \delta \quad 
 \text{for \, } 
 \eta \otimes z^{n} \in \ha{\BB(\lambda)_{\cl}}, 
\end{equation}
and define the maps 
$\ve_{j},\,\vp_{j}:
 \ha{\BB(\lambda)_{\cl}} \rightarrow \BZ_{\ge 0}$ 
 for $j \in I$ by:
%
%
\begin{align}
& \ve_{j}(\eta \otimes z^{n})=
  \max\bigl\{l \ge 0 \mid 
  e_{j}^{l}(\eta \otimes z^{n}) \ne \bzero \bigr\}
  \quad
  \text{for \, } 
  \eta \otimes z^{n} \in \ha{\BB(\lambda)_{\cl}},  
  \label{eq:aff-ve}
\\[1.5mm]
& \vp_{j}(\eta \otimes z^{n})=
  \max\bigl\{l \ge 0 \mid 
  f_{j}^{l}(\eta \otimes z^{n}) \ne \bzero \bigr\}
  \quad
  \text{for \, } 
  \eta \otimes z^{n} \in \ha{\BB(\lambda)_{\cl}}.
  \label{eq:aff-vp}
\end{align}
The proof of the following 
proposition is straightforward.
\begin{prop}
The set $\ha{\BB(\lambda)_{\cl}}$ equipped with 
the Kashiwara operators $e_{j}$, $f_{j}$ for $j \in I$, and 
the maps $\wt$, $\ve_{j}$, $\vp_{j}$ for $j \in I$, is 
a $P$-crystal. 
\end{prop}

%
\subsection{Isomorphism theorem.}
\label{subsec:isom}

Let $\lambda \in P$ be the (fixed) 
level-zero dominant integral weight 
of the form \eqref{eq:lambda}. 
For each $M \in a_{0}^{-1}\BZ$, 
let $\BB_{0}(\lambda+M\delta)$ denote
the connected component of the crystal 
$\BB(\lambda+M\delta)$ of 
all LS paths of shape $\lambda+M\delta$ 
containing the straight line $\pi_{\lambda+M\delta}$. 
%
%
\begin{lem} \label{lem:con}
Let $M_{1},\,M_{2} \in a_{0}^{-1}\BZ$. Then, 
$\BB_{0}(\lambda+M_{1}\delta)=
 \BB_{0}(\lambda+M_{2}\delta)$ if and only if 
$M_{1}-M_{2} \in d_{\lambda}\BZ$.
\end{lem}

\begin{proof}
Assume that 
$\BB_{0}(\lambda+M_{1}\delta)=
 \BB_{0}(\lambda+M_{2}\delta)$. Then 
the straight line $\pi_{\lambda+M_{1}\delta}$ is 
an LS path of shape $\lambda+M_{2}\delta$. 
Therefore, it follows from Remark~\ref{rem:LS}\,(2)
that $\lambda+M_{1}\delta \in W(\lambda+M_{2}\delta)$. 
Since $W(\lambda+M_{2}\delta)=W\lambda+M_{2}\delta$, 
we have $\lambda + (M_{1}-M_{2})\delta \in W\lambda$. 
Thus we obtain that $M_{1}-M_{2} \in d_{\lambda}\BZ$ 
from the definition of $d_{\lambda}$. 

Conversely, if $M_{1}-M_{2} \in d_{\lambda}\BZ$, then 
we deduce from the definition of $d_{\lambda}$ that 
$\lambda + (M_{1}-M_{2})\delta \in W\lambda$, and hence 
$\lambda+M_{1}\delta \in W(\lambda+M_{2}\delta)$. 
Therefore, there exists some $w \in W$ such that 
$\lambda+M_{1}\delta=w(\lambda+M_{2}\delta)$.
As a result, it follows from \eqref{eq:weyl} that 
$S_{w}\pi_{\lambda+M_{2}\delta}=
\pi_{\lambda+M_{1}\delta}$. 
Thus we obtain that 
$\pi_{\lambda+M_{2}\delta} \in 
\BB_{0}(\lambda+M_{1}\delta)$, 
and hence $\BB_{0}(\lambda+M_{1}\delta)=
 \BB_{0}(\lambda+M_{2}\delta)$.
This proves the lemma. 
\end{proof}
The following theorem (along with Corollary~\ref{cor:con}) 
is a generalization of \cite[Proposition~5.9]{GL}, 
in which $\Fg$ is nontwisted, and 
$\lambda \in P$ is a positive integer multiple of 
a level-zero fundamental weight. 
%
%
\begin{thm} \label{thm:isom}
Let $\lambda \in P$ be 
the {\rm(}fixed\,{\rm)} 
level-zero dominant integral weight 
of the form \eqref{eq:lambda}.
Then, there exists an isomorphism
$\Theta$ of $P$-crystals\,{\rm:}
\begin{equation}
\Theta:\ha{\BB(\lambda)_{\cl}} 
 \stackrel{\sim}{\rightarrow}
\bigsqcup_{
 \begin{subarray}{c}
 M \in a_{0}^{-1}\BZ \\[1mm]
 0 \le M < d_{\lambda}
 \end{subarray}}
\BB_{0}(\lambda+M\delta). 
\end{equation}
\end{thm}

In order to define the map $\Theta$, 
we need some lemmas. 
%
%
%
%
\begin{lem} \label{lem:fiber}
Let $\eta \in \BB(\lambda)_{\cl}$. 

\vsp

\noindent {\rm (1) }
The set 
$\cl^{-1}(\eta) \cap \BB_{0}(\lambda)$ is not empty. 

\vsp

\noindent {\rm (2) }
Take an arbitrary $\pi \in \cl^{-1}(\eta) \cap 
\BB_{0}(\lambda)$. Then, 
$\cl^{-1}(\eta) \cap \BB_{0}(\lambda) = 
 \bigl\{\pi+\pi_{kd_{\lambda}\delta} \mid k \in \BZ\bigr\}$.
\end{lem}

\begin{proof}
(1) It follows from Theorem~\ref{thm:NSz}\,(1)
that there exists some monomial $X$ in the root operators 
$e_{j}$, $f_{j}$ for $j \in I$ such that 
$\eta=X \cl(\pi_{\lambda})$. 
We deduce from \eqref{eq:cl-ro} that 
$\eta=
 \cl(X \pi_{\lambda})$, 
and hence that 
$X \pi_{\lambda} \in \cl^{-1}(\eta)$.
Since $X \pi_{\lambda} \in \BB_{0}(\lambda)$, 
we see that $X\pi_{\lambda} \in \cl^{-1}(\eta) \cap 
\BB_{0}(\lambda)$. 
This proves part (1) of the lemma.

\vsp

\noindent 
(2) First we prove that 
$\cl^{-1}(\eta) \cap \BB_{0}(\lambda) \supset
 \bigl\{\pi + \pi_{kd_{\lambda}\delta} \mid k \in \BZ\bigr\}$.
Take an arbitrary $k \in \BZ$, and set 
$\pi^{\prime}:=\pi+\pi_{kd_{\lambda}\delta}$. 
Since it is obvious that $\pi^{\prime} \in \cl^{-1}(\eta)$, 
it remains to show that 
$\pi^{\prime} \in \BB_{0}(\lambda)$. 
Because $\pi \in \BB_{0}(\lambda)$ by assumption, 
there exists some monomial $X^{\prime}$ in the root operators 
$e_{j}$, $f_{j}$ for $j \in I$ such that 
$X^{\prime} \pi = \pi_{\lambda}$. 
Therefore, using Lemma~\ref{lem:d-shift2}, 
we deduce that 
\begin{equation*}
X^{\prime} \pi^{\prime} =
X^{\prime} (\pi + \pi_{kd_{\lambda}\delta})=
X^{\prime} \pi + \pi_{kd_{\lambda}\delta} =
\pi_{\lambda}+\pi_{kd_{\lambda}\delta}=
\pi_{\lambda+kd_{\lambda}\delta}.
\end{equation*}
Note that 
$\lambda+kd_{\lambda}\delta \in W\lambda$ 
by the definition of $d_{\lambda}$. 
Hence there exists some $w \in W$ such that 
$\lambda=w(\lambda+kd_{\lambda}\delta)$.
Therefore, we see from \eqref{eq:weyl} that 
$S_{w}\pi_{\lambda+kd_{\lambda}\delta}=
\pi_{\lambda}$. Thus we obtain that 
$S_{w}X^{\prime}\pi^{\prime}
 =S_{w}\pi_{\lambda+kd_{\lambda}\delta}
 =\pi_{\lambda}$, 
and hence $\pi^{\prime} \in \BB_{0}(\lambda)$.
This proves the desired inclusion $\supset$.

Next we prove that 
$\cl^{-1}(\eta) \cap \BB_{0}(\lambda) \subset
 \bigl\{\pi + \pi_{kd_{\lambda}\delta} \mid k \in \BZ\bigr\}$.
Let $\pi^{\prime} \in 
\cl^{-1}(\eta) \cap \BB_{0}(\lambda)$. 
Since $\cl(\pi)=\cl(\pi^{\prime})=\eta$, 
there exists a piecewise linear, continuous function 
$F(t):[0,1] \rightarrow \BQ$ such that 
$F(0)=0$, $F(1) \in a_{0}^{-1}\BZ$, and 
$\pi^{\prime}(t)=\pi(t)+F(t)\delta$
for all $t \in [0,1]$.
Define a path $\psi \in \BP$ by: 
$\psi(t)=F(t)\delta$, $t \in [0,1]$; 
note that $\pi^{\prime}=\pi+\psi$. 
We will show that $\psi=\pi_{kd_{\lambda}\delta}$ 
for some $k \in \BZ$. 
Since $X^{\prime} \pi = \pi_{\lambda}$ 
for the above monomial $X^{\prime}$ 
in the root operators, 
we deduce from Lemma~\ref{lem:d-shift2} that 
$X^{\prime}\pi^{\prime}=
X^{\prime}(\pi+\psi)=
X^{\prime}\pi+\psi=
\pi_{\lambda}+\psi$.
If we write the path $\psi \in \BP$ as:
$\psi=
(N_{1}^{\prime}\delta,\,
 N_{2}^{\prime}\delta,\,
 \dots,\,
 N_{s^{\prime}}^{\prime}\delta \, ; \, 
 \sigma_{0},\,\sigma_{1},\,\dots,\,\sigma_{s^{\prime}})$,
with $N_{1}^{\prime},\,N_{2}^{\prime},\,\dots,\,
N_{s^{\prime}}^{\prime} \in \BZ$ and 
$0=\sigma_{0} < \sigma_{1} < \cdots < \sigma_{s^{\prime}}=1$, 
then we have 
\begin{equation*}
\BB_{0}(\lambda) \ni 
X^{\prime}\pi^{\prime} =
\pi_{\lambda}+\psi =
(\lambda+N_{1}^{\prime}\delta,\,
 \lambda+N_{2}^{\prime}\delta,\,
 \dots,\,
 \lambda+N_{s^{\prime}}^{\prime}\delta \, ; \, 
 \sigma_{0},\,\sigma_{1},\,\dots,\,
 \sigma_{s^{\prime}}). 
\end{equation*}
Since $X^{\prime}\pi^{\prime} \in \BB(\lambda)$, 
it follows from Remark~\ref{rem:LS}\,(1) that 
$\lambda+N_{1}^{\prime}\delta \ge 
 \lambda+N_{2}^{\prime}\delta \ge \dots \ge 
 \lambda+N_{s^{\prime}}^{\prime}\delta$,
and hence that 
$N_{1}^{\prime} \le N_{2}^{\prime} \le \dots \le 
 N_{s^{\prime}}^{\prime}$ by Remark~\ref{rem:Q+}\,(1).
Also, because $X^{\prime}\pi^{\prime} \in \BB_{0}(\lambda)$, 
there exists some monomial $X^{\prime\prime}$ 
in the root operators $e_{j}$, $f_{j}$ for $j \in I$ such that 
$X^{\prime\prime}X^{\prime}\pi^{\prime}=\pi_{\lambda}$. 
Therefore, using Lemma~\ref{lem:d-shift2}, we deduce that 
\begin{align*}
\BB_{0}(\lambda) \ni 
X^{\prime\prime}\pi_{\lambda} & = 
X^{\prime\prime}(X^{\prime}\pi^{\prime}-\psi) = 
X^{\prime\prime}X^{\prime}\pi^{\prime}-\psi = 
\pi_{\lambda}-\psi \\[1.5mm]
& = 
(\lambda-N_{1}^{\prime}\delta,\,
 \lambda-N_{2}^{\prime}\delta,\,
 \dots,\,
 \lambda-N_{s^{\prime}}^{\prime}\delta \, ; \, 
 \sigma_{0},\,\sigma_{1},\,\dots,\,
 \sigma_{s^{\prime}}). 
\end{align*}
Again, since $X^{\prime\prime}\pi_{\lambda} \in \BB(\lambda)$, 
it follows from Remark~\ref{rem:LS}\,(1) that 
$\lambda-N_{1}^{\prime}\delta \ge 
 \lambda-N_{2}^{\prime}\delta \ge \dots \ge 
 \lambda-N_{s^{\prime}}^{\prime}\delta$,
and hence that 
$N_{1}^{\prime} \ge N_{2}^{\prime} \ge \dots \ge 
 N_{s^{\prime}}^{\prime}$ by Remark~\ref{rem:Q+}\,(1). 
Thus we obtain that 
$N_{1}^{\prime} = N_{2}^{\prime} = \dots =
 N_{s^{\prime}}^{\prime}$.
But, since $\lambda+N_{1}^{\prime}\delta \in W\lambda$ 
by Remark~\ref{rem:LS}\,(1), 
we see from the definition of $d_{\lambda}$ 
that $N_{1}^{\prime}=kd_{\lambda}$ 
for some $k \in \BZ$.
Therefore, we conclude that 
$N_{1}^{\prime} = N_{2}^{\prime} = \dots =
 N_{s^{\prime}}^{\prime}=kd_{\lambda}$, and 
hence that
\begin{equation*}
\psi=
(\underbrace{
 kd_{\lambda}\delta,\ 
 kd_{\lambda}\delta,\, 
 \dots,\,
 kd_{\lambda}\delta}_{
 \text{$s^{\prime}$ times}} \ ; \ 
 \sigma_{0},\,\sigma_{1},\,\dots,\,\sigma_{s^{\prime}})
=\pi_{kd_{\lambda}\delta},
\end{equation*}
as desired. This proves part (2) of the lemma. 
\end{proof}

For each $\eta \in \BB(\lambda)_{\cl}$, we define 
$\pi_{\eta}^{0} \in \BP$ as follows: 
Take an arbitrary 
$\pi \in \cl^{-1}(\eta) \cap \BB_{0}(\lambda)$. 
It follows from Lemma~\ref{lem:wt} that 
$\pi(1) \in P$ can be written as 
$\lambda-\alpha+n^{\prime}\delta$, 
with $\alpha \in a_{0}^{-1}\fin{Q}_{+}$ and 
$n^{\prime} \in a_{0}^{-1}\BZ$. 
Then we set
\begin{equation}
\pi_{\eta}^{0}:=\pi-\pi_{n^{\prime}\delta}. 
\end{equation}
%
%
\begin{rem} \label{rem:fiber}
We see from Lemma~\ref{lem:fiber}\,(2) that 
$\pi_{\eta}^{0}$ does not depend on the choice of 
$\pi \in \cl^{-1}(\eta) \cap \BB_{0}(\lambda)$. 
\end{rem}
%
%
\begin{lem} \label{lem:zero}
Let $\eta \in \BB(\lambda)_{\cl}$. 
Then, for each $j \in I$, we have\,{\rm:}

\vsp 

\noindent {\rm (1) }
$\pi_{e_{j}\eta}^{0}=
  e_{j}\pi_{\eta}^{0} - 
  \pi_{\delta_{j,0}a_{0}^{-1}\delta}$
  if $e_{j}\eta \ne \bzero$, 

\vsp 

\noindent {\rm (2) }
$\pi_{f_{j}\eta}^{0}=
  f_{j}\pi_{\eta}^{0} + 
  \pi_{\delta_{j,0}a_{0}^{-1}\delta}$
  if $f_{j}\eta \ne \bzero$. 
\end{lem}

\begin{proof}
We give only the proof of part~(1); 
the proof of part~(2) is similar.
Let $\pi \in \cl^{-1}(\eta) \cap 
\BB_{0}(\lambda)$, and assume that 
$\pi(1)=\lambda-\alpha+ n^{\prime}\delta$, 
with $\alpha \in a_{0}^{-1}\fin{Q}_{+}$ and  
$n^{\prime} \in a_{0}^{-1}\BZ$. 
Note that $e_{j}\pi \in \cl^{-1}(e_{j}\eta) \cap 
\BB_{0}(\lambda)$ by \eqref{eq:cl-ro}. 
Because $\alpha_{0}=a_{0}^{-1}(\delta-\theta)$, 
we have 
\begin{align*}
(e_{j}\pi)(1) & = \pi(1)+\alpha_{j}=
\lambda-\alpha+ n^{\prime}\delta + \alpha_{j} \\[1.5mm]
& = 
\begin{cases}
\lambda-(\alpha-\alpha_{j})+n^{\prime}\delta 
 & \text{if $j \ne 0$}, \\[3mm]
\lambda-(\alpha+a_{0}^{-1}\theta)+
(n^{\prime}+a_{0}^{-1})\delta 
 & \text{if $j = 0$}.
\end{cases}
\end{align*}
Therefore, it follows from the definition of 
$\pi_{e_{j}\eta}^{0}$ that 
%
%
\begin{equation} \label{eq:zero-1}
\pi_{e_{j}\eta}^{0} = 
 e_{j}\pi-
 \pi_{(n^{\prime}+\delta_{j,0}a_{0}^{-1})\delta}. 
\end{equation}
Also, because 
$\pi_{\eta}^{0} = \pi-\pi_{n^{\prime}\delta}$ 
by the definition of $\pi_{\eta}^{0}$, 
we have $e_{j}\pi_{\eta}^{0} 
=e_{j}\pi - \pi_{n^{\prime}\delta}$ 
by Lemma~\ref{lem:d-shift2}. 
Combining this equation with \eqref{eq:zero-1}, 
we obtain that 
$\pi_{e_{j}\eta}^{0}=
  e_{j}\pi_{\eta}^{0} - 
  \pi_{\delta_{j,0}a_{0}^{-1}\delta}$. 
This proves part~(1) of the lemma. 
\end{proof}

Now, let us define a map 
$\Theta:\ha{\BB(\lambda)_{\cl}} 
\rightarrow \BP$ by:
%
%
\begin{equation} \label{eq:img0}
\Theta(\eta \otimes z^{n})=
 \pi_{\eta}^{0} + \pi_{n\delta}
\quad \text{for \, } 
\eta \otimes z^{n} \in \ha{\BB(\lambda)_{\cl}}. 
\end{equation}
%
%
\begin{lem} \label{lem:img1}
Let $\eta \otimes z^{n} \in \ha{\BB(\lambda)_{\cl}}$, 
and let $\pi \in \cl^{-1}(\eta) \cap \BB_{0}(\lambda)$. 
Assume that $\pi(1)=\lambda-\alpha+n^{\prime}\delta$, 
with $\alpha \in a_{0}^{-1}\fin{Q}_{+}$ 
and $n^{\prime} \in a_{0}^{-1}\BZ$. 
Then, $\Theta(\eta \otimes z^{n}) \in 
\BB(\lambda+(n-n^{\prime})\delta)$. 
\end{lem}
\begin{proof}
It follows from the definition of the map $\Theta$ that
%
%
\begin{equation} \label{eq:img1-1}
\Theta(\eta \otimes z^{n}) =
\pi_{\eta}^{0} + \pi_{n\delta} =
(\pi-\pi_{n^{\prime}\delta})+\pi_{n\delta} =
\pi + \pi_{(n-n^{\prime})\delta}.
\end{equation}
Therefore, we see from Lemma~\ref{lem:d-shift4} that 
$\Theta(\eta \otimes z^{n}) \in 
\BB(\lambda+(n-n^{\prime})\delta)$. 
This proves the lemma.
\end{proof}
\begin{proof}[Proof of Theorem~\ref{thm:isom}]
We  will establish this theorem by proving Claims~\ref{c:wt}, 
\ref{c:kas}, and \ref{c:inj-img} below.
%
%
\begin{claim} \label{c:wt}
Let $\eta \otimes z^{n} \in \ha{\BB(\lambda)_{\cl}}$. 
Then, $\wt(\eta \otimes z^{n})=
\wt(\Theta(\eta \otimes z^{n}))$. 
\end{claim}

Let $\pi \in \cl^{-1}(\eta) \cap \, \BB_{0}(\lambda)$, and 
assume that $\pi(1)=\lambda-\alpha+n^{\prime}\delta$, 
with $\alpha \in a_{0}^{-1}\fin{Q}_{+}$ and 
$n^{\prime} \in a_{0}^{-1}\BZ$. Then we see from 
the definition of the $\BQ$-linear embedding 
$\af:\Fh^{\ast}/\BQ\delta \hookrightarrow \Fh^{\ast}$ 
that 
$\af(\eta(1))
  =\af(\cl(\pi(1)))
  =\af(\cl(\lambda))-\af(\cl(\alpha))
  =\lambda-\alpha$.
Therefore, by \eqref{eq:aff-wt}, we have
\begin{equation*}
\wt(\eta \otimes z^{n}) 
 = \af(\eta(1))+n\delta 
 = \lambda-\alpha+n\delta.
\end{equation*}
Also, it follows from the definition of 
$\pi_{\eta}^{0}$ that $\pi_{\eta}^{0}(1)=
\pi(1)-n^{\prime}\delta= \lambda-\alpha$. 
Therefore, by the definition of the map $\Theta$, 
we have
\begin{equation*}
\wt(\Theta(\eta \otimes z^{n})) = 
\wt(\pi_{\eta}^{0}+\pi_{n\delta})=
\pi_{\eta}^{0}(1)+\pi_{n\delta}(1)=
\lambda-\alpha+n\delta.
\end{equation*}
This proves Claim~\ref{c:wt}. 
%
%
%
\begin{claim} \label{c:kas}
Let $\eta \otimes z^{n} \in \ha{\BB(\lambda)_{\cl}}$. 
Then, for all $j \in I$, 
\begin{align}
& 
\Theta\bigl(e_{j}(\eta \otimes z^{n})\bigr)=
  e_{j}\Theta(\eta \otimes z^{n})
\quad \text{and} \quad
\Theta\bigl(f_{j}(\eta \otimes z^{n})\bigr)=
  f_{j}\Theta(\eta \otimes z^{n}), \\[1.5mm]
& 
\ve_{j}(\eta \otimes z^{n})=
  \ve_{j}(\Theta(\eta \otimes z^{n}))
\quad \text{and} \quad
\vp_{j}(\eta \otimes z^{n})=
  \vp_{j}(\Theta(\eta \otimes z^{n})).
\end{align}
Here $\Theta(\bzero)$ is understood to be $\bzero$. 
\end{claim}

We give only the proof of the equations 
$\Theta\bigl(e_{j}(\eta \otimes z^{n})\bigr)=
 e_{j}\Theta(\eta \otimes z^{n})$ and 
$\ve_{j}(\eta \otimes z^{n})=
  \ve_{j}(\Theta(\eta \otimes z^{n}))$, since 
the proof of the corresponding equations for $f_{j}$ 
is similar.
First, let us show 
that $e_{j}(\eta \otimes z^{n}) = \bzero$ if and only if 
$e_{j}\Theta(\eta \otimes z^{n})=\bzero$. 
Indeed, we have:
\begin{align*}
e_{j}(\eta \otimes z^{n}) = \bzero \quad 
& \Leftrightarrow \quad 
e_{j}\eta = \bzero
\quad \text{by the definition of $e_{j}$ for 
            $\ha{\BB(\lambda)_{\cl}}$}, \\
& \Leftrightarrow \quad 
e_{j}\pi = \bzero 
\quad \text{by \eqref{eq:cl-ro}}, \\
& \Leftrightarrow \quad 
e_{j}\pi_{\eta}^{0} = \bzero 
\quad \text{by Lemma~\ref{lem:d-shift2} and 
      the definition of $\pi_{\eta}^{0}$}, \\
& \Leftrightarrow \quad 
e_{j}(\pi_{\eta}^{0}+\pi_{n\delta}) = \bzero 
\quad \text{by Lemma~\ref{lem:d-shift2}}, \\
& \Leftrightarrow \quad 
e_{j}\Theta(\eta \otimes z^{n}) = \bzero 
\quad \text{by the definition of the map $\Theta$}.
\end{align*}

Now assume that 
$e_{j}\pi_{\eta}^{0} \ne \bzero$ 
(and hence $e_{j}\eta \ne \bzero$).
From the definitions of 
the Kashiwara operator $e_{j}$ for 
$\ha{\BB(\lambda)_{\cl}}$ and the map $\Theta$, 
we deduce that 
$\Theta\bigl(e_{j}(\eta \otimes z^{n})\bigr) = 
\Theta(e_{j}\eta \otimes z^{n+\delta_{j,0}a_{0}^{-1}})
=\pi_{e_{j}\eta}^{0} + \pi_{(n+\delta_{j,0}a_{0}^{-1})\delta}$.
Furthermore, this equals 
$e_{j}\pi_{\eta}^{0}+ \pi_{n\delta}$ 
by Lemma~\ref{lem:zero}\,(1).
Also, we see from Lemma~\ref{lem:d-shift2} that 
$e_{j}\Theta(\eta \otimes z^{n}) = 
e_{j}(\pi_{\eta}^{0}+\pi_{n\delta}) = 
e_{j}\pi_{\eta}^{0}+\pi_{n\delta}$.
Thus we have proved the equation 
$\Theta\bigl(e_{j}(\eta \otimes z^{n})\bigr)=
 e_{j}\Theta(\eta \otimes z^{n})$.
This equation along with \eqref{eq:aff-ve} 
immediately implies the equation
$\ve_{j}(\eta \otimes z^{n})=
 \ve_{j}(\Theta(\eta \otimes z^{n}))$. 
This proves Claim~\ref{c:kas}. 

\vsp\vsp

It follows from 
Claims~\ref{c:wt} and \ref{c:kas}, 
in view of Lemma~\ref{lem:img1}, that 
the image $\img \Theta$ of $\Theta$ is a subcrystal of 
$\bigcup_{M \in a_{0}^{-1}\BZ} \BB(\lambda+M\delta)$.
%
%
%
\begin{claim} \label{c:inj-img}
The map $\Theta:\ha{\BB(\lambda)_{\cl}} \rightarrow \BP$ is 
injective, and its image $\img \Theta$ is equal to
%
%
\begin{equation} \label{eq:inj-img}
\bigsqcup_{
 \begin{subarray}{c}
 M \in a_{0}^{-1}\BZ \\[1mm]
 0 \le M < d_{\lambda}
 \end{subarray}}
\BB_{0}(\lambda+M\delta).
\end{equation}
\end{claim}
%
First we show that 
the map $\Theta$ is injective. Assume that 
$\Theta(\eta_{1} \otimes z^{n_{1}})=
 \Theta(\eta_{2} \otimes z^{n_{2}})$ for some 
$\eta_{1} \otimes z^{n_{1}},\,
 \eta_{2} \otimes z^{n_{2}} \in \ha{\BB(\lambda)_{\cl}}$. 
It immediately follows from the definition 
of the map $\Theta$ that 
$\pi_{\eta_{1}}^{0} + \pi_{n_{1}\delta} = 
\pi_{\eta_{2}}^{0} + \pi_{n_{2}\delta}$. 
Also, we have
$\cl(\pi_{\eta_{1}}^{0})=\eta_{1}$ and 
$\cl(\pi_{\eta_{2}}^{0})=\eta_{2}$ 
by the definitions of 
$\pi_{\eta_{1}}^{0}$ and 
$\pi_{\eta_{2}}^{0}$.
Hence we obtain that 
\begin{equation*}
\eta_{1}=\cl(\pi_{\eta_{1}}^{0})=
\cl(\pi_{\eta_{1}}^{0}+\pi_{n_{1}\delta})=
\cl(\pi_{\eta_{2}}^{0}+\pi_{n_{2}\delta})=
\cl(\pi_{\eta_{2}}^{0})=\eta_{2}.
\end{equation*}
This implies that 
$\pi_{n_{1}\delta}=\pi_{n_{2}\delta}$, 
and hence $n_{1}=n_{2}$. 
Thus we have shown that 
$\eta_{1} \otimes z^{n_{1}}=\eta_{2} \otimes z^{n_{2}}$, 
and hence that the map $\Theta$ is injective. 

Next we show that the image $\img \Theta$ of $\Theta$ 
is equal to the subset \eqref{eq:inj-img} of $\BP$; 
note that
\begin{equation*}
\bigcup_{M \in a_{0}^{-1}\BZ} 
\BB_{0}(\lambda+M\delta)=
\bigsqcup_{
 \begin{subarray}{c}
 M \in a_{0}^{-1}\BZ \\[1mm]
 0 \le M < d_{\lambda}
 \end{subarray}}
\BB_{0}(\lambda+M\delta)
\end{equation*}
as seen from Lemma~\ref{lem:con}.
Since $\pi^{0}_{\cl(\pi_{\lambda})}$ is clearly equal to 
$\pi_{\lambda}$ by the definition of $\pi^{0}_{\cl(\pi_{\lambda})}$, 
it follows from the definition of $\Theta$ that 
$\Theta(\cl(\pi_{\lambda}) \otimes t^{M})=
 \pi_{\lambda}+\pi_{M\delta}=
 \pi_{\lambda+M\delta}$ 
for every $M \in a_{0}^{-1}\BZ$. 
Consequently, 
we have $\pi_{\lambda+M\delta} \in \img \Theta$ 
for all $M \in a_{0}^{-1}\BZ$. 
This implies that
%
%
\begin{equation} \label{eq:img-sub}
\bigcup_{M \in a_{0}^{-1}\BZ} 
\BB_{0}(\lambda+M\delta) \subset \img \Theta,
\end{equation}
since $\img \Theta$ is a subcrystal of 
$\bigcup_{M \in a_{0}^{-1}\BZ} \BB(\lambda+M\delta)$.

Now, let us show the reverse inclusion. 
As seen in the proof of Lemma~\ref{lem:img1} 
(see \eqref{eq:img1-1}), every element of $\img \Theta$ is 
of the form $\pi+\pi_{n^{\prime\prime}\delta}$, 
with $\pi \in \BB_{0}(\lambda)$ and 
$n^{\prime\prime} \in a_{0}^{-1}\BZ$.
Since $\pi \in \BB_{0}(\lambda)$, 
there exists some monomial $X$ 
in the root operators 
$e_{j}$, $f_{j}$ for $j \in I$ such that 
$\pi=X\pi_{\lambda}$. 
Hence it follows from 
Lemma~\ref{lem:d-shift2} that 
$X(\pi_{\lambda+n^{\prime\prime}\delta})=
X(\pi_{\lambda}+\pi_{n^{\prime\prime}\delta})=
\pi+\pi_{n^{\prime\prime}\delta}$.
Therefore, we obtain that 
%
%
\begin{equation} \label{eq:img-sup}
\bigcup_{M \in a_{0}^{-1}\BZ} 
\BB_{0}(\lambda+M\delta) \supset \img \Theta.
\end{equation}
Combining \eqref{eq:img-sub} and \eqref{eq:img-sup}, 
we conclude that 
\begin{equation*}
\img \Theta
=
\bigcup_{M \in a_{0}^{-1}\BZ} 
\BB_{0}(\lambda+M\delta)
=
\bigsqcup_{
 \begin{subarray}{c}
 M \in a_{0}^{-1}\BZ \\[1mm]
 0 \le M < d_{\lambda}
 \end{subarray}}
\BB_{0}(\lambda+M\delta).
\end{equation*}
This proves Claim~\ref{c:inj-img}.

\vsp\vsp

Thus the proof of 
Theorem~\ref{thm:isom} is completed. 
\end{proof}
%
%
\begin{cor} \label{cor:con}
The connected component $\BB_{0}(\lambda)$ of $\BB(\lambda)$ 
is isomorphic as a $P$-crystal to the subcrystal of 
$\ha{\BB(\lambda)_{\cl}}$ consisting of 
all elements $\eta \otimes z^{n} \in \ha{\BB(\lambda)_{\cl}}$ 
that satisfy the following condition {\rm (C):}

\vsp

\noindent 
{\rm (C) } 
Let $\pi \in \cl^{-1}(\eta) \cap \BB_{0}(\lambda)$, and 
assume that $\pi(1)=\lambda-\alpha+n^{\prime}\delta$ 
with $\alpha \in a_{0}^{-1}\fin{Q}_{+}$ and 
$n^{\prime} \in a_{0}^{-1}\BZ$. 
Then it holds that $n^{\prime}-n \in d_{\lambda}\BZ$.
\end{cor}

\begin{rem}
As is easily seen from Lemma~\ref{lem:fiber}\,(2), 
it does not depend on the choice of 
$\pi \in \cl^{-1}(\eta) \cap \BB_{0}(\lambda)$
whether or not an element 
$\eta \otimes z^{n} \in \ha{\BB(\lambda)_{\cl}}$ 
satisfies condition~(C).
\end{rem}

\begin{proof}[Proof of Corollary~\ref{cor:con}]
Let $\eta \otimes z^{n} \in \ha{\BB(\lambda)_{\cl}}$, 
and let $\pi \in \cl^{-1}(\eta) \cap \BB_{0}(\lambda)$. 
Assume that $\pi(1)=\lambda-\alpha+n^{\prime}\delta$, 
with $\alpha \in a_{0}^{-1}\fin{Q}_{+}$ and 
$n^{\prime} \in a_{0}^{-1}\BZ$. 
Then, as in \eqref{eq:img1-1}, 
we have $\Theta(\eta \otimes z^{n})=
\pi+\pi_{(n-n^{\prime})\delta}$. 
By the argument used to show \eqref{eq:img-sup} above, 
we can show that
$\Theta(\eta \otimes z^{n})=
\pi+\pi_{(n-n^{\prime})\delta} \in 
\BB_{0}(\lambda+(n-n^{\prime})\delta)$.
Therefore, if the element 
$\eta \otimes z^{n} \in \ha{\BB(\lambda)_{\cl}}$ 
satisfies condition (C) 
(i.e., that $n^{\prime}-n \in d_{\lambda}\BZ$), then 
it follows from Lemma~\ref{lem:con} that
$\Theta(\eta \otimes z^{n}) \in 
\BB_{0}(\lambda+(n-n^{\prime})\delta) = 
\BB_{0}(\lambda)$.
Conversely, if 
$\Theta(\eta \otimes z^{n}) \in \BB_{0}(\lambda)$, then 
we obtain that $\BB_{0}(\lambda+(n-n^{\prime})\delta)=
\BB_{0}(\lambda)$, and hence that $n-n^{\prime} \in 
d_{\lambda}\BZ$ by Lemma~\ref{lem:con}. 
Therefore, 
the element $\eta \otimes z^{n} \in 
\ha{\BB(\lambda)_{\cl}}$ 
satisfies condition (C). 
This completes the proof of the corollary. 
\end{proof}

\appendix

%
\section{Appendix.}
\label{sec:appendix}

%
\subsection{Relation with
the crystal base of an extremal weight module.}
\label{subsec:rel-ext}

For an integral weight $\lambda \in P$, 
we denote by $\CB(\lambda)$ the crystal base of 
the extremal weight $U_{q}(\Fg)$-module 
$V(\lambda)$ of extremal weight $\lambda$ 
(see \cite[Sections~3 and 5]{Kaslz}), 
where $U_{q}(\Fg)$ is 
the quantized universal 
enveloping algebra of $\Fg$, and 
denote by $u_{\lambda} \in \CB(\lambda)$
the extremal element of weight $\lambda$ 
corresponding to the canonical generator of $V(\lambda)$. 
The aim in the Appendix is to prove the following proposition.
%
%
%
\begin{prop} \label{prop:ext}
Let $\lambda=\sum_{i \in I_{0}}m_{i}\vpi_{i}$, 
with $m_{i} \in \BZ_{\ge 0}$ for $i \in I_{0}$, be 
a level-zero dominant integral weight. 
Then, the crystal base $\CB(\lambda)$ 
is isomorphic as a $P$-crystal to the crystal
$\BB(\lambda)$  of all LS paths of shape $\lambda$ 
if and only if $\lambda$ is of the form $m_{i}\vpi_{i}$ 
for some $i \in I_{0}$. 
\end{prop}

The main ingredient in our proof of 
Proposition~\ref{prop:ext} is 
the following theorem 
(see \cite[Remark~4.17]{BN}, and also 
\cite[Conjecture~13.1\,(iii)]{Kaslz}). 
%
%
\begin{thm} \label{thm:BN}
Let $\lambda=\sum_{i \in I_{0}}m_{i}\vpi_{i}$, 
with $m_{i} \in \BZ_{\ge 0}$ for $i \in I_{0}$, be 
a level-zero dominant integral weight. 
Then, there exists an isomorphism of $P$-crystals 
$\CB(\lambda) \stackrel{\sim}{\rightarrow} 
 \bigotimes_{i \in I_{0}} \CB(m_{i}\vpi_{i})$ 
that maps $u_{\lambda} \in \CB(\lambda)$ to 
$\ti{u}_{\lambda}:=
 \bigotimes_{i \in I_{0}} u_{m_{i}\vpi_{i}} 
 \in \bigotimes_{i \in I_{0}} \CB(m_{i}\vpi_{i})$. 
\end{thm}
Let $\CB$ be the crystal base of 
an extremal weight $U_{q}(\Fg)$-module, 
or a tensor product of 
the crystal bases of 
extremal weight $U_{q}(\Fg)$-modules. 
Then we have an action, 
denoted by $S_{w}$, $w \in W$, of the Weyl group $W$ 
on the set $\CB$ such that $\wt(S_{w}b)=w\wt(b)$ 
for every $w \in W$ and $b \in \CB$ 
(see \cite[Section~7]{Kasmod}). 
The proof of the following lemma is standard 
(see the proof of \cite[Lemma~1.6\,(1)]{AK}).
%
%
\begin{lem} \label{lem:w-ext}
Let $\lambda=\sum_{i \in I_{0}}m_{i}\vpi_{i}$, 
 with $m_{i} \in \BZ_{\ge 0}$ for $i \in I_{0}$, 
be a level-zero dominant integral weight. 
For every $w \in W$, we have 
$S_{w} \ti{u}_{\lambda}=
 \bigotimes_{i \in I_{0}} S_{w}u_{m_{i}\vpi_{i}}$.
\end{lem}
%
%
%
\begin{lem} \label{lem:wt-lam}
Let $\lambda=\sum_{i \in I_{0}}m_{i}\vpi_{i}$, 
with $m_{i} \in \BZ_{\ge 0}$ for $i \in I_{0}$, 
be a level-zero dominant integral weight. 
Then each connected component of $\BB(\lambda)$ 
contains at most one element of weight $\lambda$.
\end{lem}

\begin{proof}
Take an arbitrary connected component
of $\BB(\lambda)$, and denote it by $\BB_{1}(\lambda)$. 
We know from Theorem~\ref{thm:comps}\,(2)
that $\BB_{1}(\lambda)$ contains 
a unique (extremal) LS path $\pi$ 
of shape $\lambda$ 
having an expression of the form \eqref{eq:comps}: 
\begin{equation*}
\pi=(
  \lambda-N_{1}\delta,\,\dots,\,
  \lambda-N_{s-1}\delta,\ \lambda \ ; \ 
  \tau_{0},\,\tau_{1},\,\dots,\,\tau_{s-1},\,\tau_{s})
  \in \BB_{1}(\lambda). 
\end{equation*}
Define a path $\psi \in \BP$ by:
\begin{equation*}
\psi=(N_{1}\delta,\,\dots,\,
  N_{s-1}\delta,\ 0 \ ; \ 
  \tau_{0},\,\tau_{1},\,\dots,\,\tau_{s-1},\,\tau_{s}), 
\end{equation*}
and $n_{1} \in a_{0}^{-1}\BZ$ by: $\psi(1)=n_{1}\delta$. 
Then we deduce from Lemma~\ref{lem:d-shift2} that 
$\BB_{1}(\lambda)=\BB_{0}(\lambda)-\psi$. 
Therefore, the number of elements of 
weight $\lambda$ in $\BB_{1}(\lambda)$ is equal to 
the number of elements of weight 
$\lambda+n_{1}\delta$ in $\BB_{0}(\lambda)$. 
But, we see from Theorem~\ref{thm:NSz}\,(3) 
and the definition of the $P$-crystal 
$\ha{\BB(\lambda)_{\cl}}$ 
that the number of elements of weight 
$\lambda+n_{1}\delta$ in $\ha{\BB(\lambda)_{\cl}}$ 
is exactly $1$. 
Hence it follows from Theorem~\ref{thm:isom} that 
the number of elements of weight 
$\lambda+n_{1}\delta$ in $\BB_{0}(\lambda)$ 
is at most $1$.
This proves the lemma. 
\end{proof}

\begin{proof}[Proof of Proposition~\ref{prop:ext}]
We proved in \cite[Corollary~3.8.1]{NSp2} that 
for each $m \in \BZ_{\ge 0}$ and $i \in I_{0}$, 
the crystal base $\CB(m\vpi_{i})$ is 
isomorphic to $\BB(m\vpi_{i})$ as a $P$-crystal. 
Therefore, it suffices to show the ``only if'' part. 
Assume that the support $\Supp(\lambda)$ of $\lambda$ 
contains at least two elements, i.e., that 
$\# \Supp(\lambda)  \ge 2$. 
Let us denote by $\CB_{0}(\lambda)$ the connected component of 
$\CB(\lambda)$ containing $u_{\lambda}$. 
Note that under the isomorphism of $P$-crystals
in Theorem~\ref{thm:BN}, $\CB_{0}(\lambda)$ is isomorphic 
to the connected component of 
$\bigotimes_{i \in I_{0}} \CB(m_{i}\vpi_{i})$ 
containing $\ti{u}_{\lambda}:=
\bigotimes_{i \in I_{0}} u_{m_{i}\vpi_{i}}$, 
which we denote by $\ti{\CB}_{0}(\lambda)$. 
\begin{aclaim}
The number of elements of weight $\lambda$ in 
$\CB_{0}(\lambda)$ is greater than or 
equal to $2$. 
\end{aclaim}

Because $\# \Supp(\lambda)  \ge 2$, 
we can take $\beta_{1},\beta_{2} \in \Gamma$ 
such that $(\lambda,\beta_{1})=
(\lambda,\beta_{2})=0$ and such that
$(\vpi_{i_{0}},\beta_{1}) \ne (\vpi_{i_{0}},\beta_{2})$ 
for some $i_{0} \in I_{0}$ with $m_{i_{0}} > 0$. 
We see from Lemma~\ref{lem:w-ext} 
that 
\begin{equation*}
S_{t_{\beta_{1}}}\ti{u}_{\lambda} = 
 \bigotimes_{i \in I_{0}}
 S_{t_{\beta_{1}}}u_{m_{i}\vpi_{i}}, 
\qquad
S_{t_{\beta_{2}}}\ti{u}_{\lambda} = 
 \bigotimes_{i \in I_{0}} 
 S_{t_{\beta_{2}}}u_{m_{i}\vpi_{i}}. 
\end{equation*}
Since $(\lambda,\beta_{1})=
(\lambda,\beta_{2})=0$, it immediately follows 
that both of the elements 
$S_{t_{\beta_{1}}}\ti{u}_{\lambda}$ and 
$S_{t_{\beta_{2}}}\ti{u}_{\lambda}$ are 
of weight $\lambda$. 
Note that the weight of the $i_{0}$-th factor 
$S_{t_{\beta_{1}}}u_{m_{i_{0}}\vpi_{i_{0}}}$ of 
$S_{t_{\beta_{1}}}\ti{u}_{\lambda}$ equals
$m_{i_{0}}\vpi_{i_{0}}-
 m_{i_{0}}(\vpi_{i_{0}},\beta_{1})\delta$, and 
 the weight of the $i_{0}$-th factor 
$S_{t_{\beta_{2}}}u_{m_{i_{0}}\vpi_{i_{0}}}$ of 
$S_{t_{\beta_{2}}}\ti{u}_{\lambda}$
equals
$m_{i_{0}}\vpi_{i_{0}}-
 m_{i_{0}}(\vpi_{i_{0}},\beta_{2})\delta$. 
But, since 
$(\vpi_{i_{0}},\beta_{1}) \ne (\vpi_{i_{0}},\beta_{2})$, 
the weight of 
$S_{t_{\beta_{1}}}u_{m_{i_{0}}\vpi_{i_{0}}}$ is not equal to 
that of $S_{t_{\beta_{2}}}u_{m_{i_{0}}\vpi_{i_{0}}}$. 
This implies that 
$S_{t_{\beta_{1}}}u_{m_{i_{0}}\vpi_{i_{0}}} \ne 
 S_{t_{\beta_{2}}}u_{m_{i_{0}}\vpi_{i_{0}}}$, and hence 
that 
$S_{t_{\beta_{1}}}\ti{u}_{\lambda} \ne 
 S_{t_{\beta_{2}}}\ti{u}_{\lambda}$. 
Because both of the elements 
$S_{t_{\beta_{1}}}\ti{u}_{\lambda}$ and 
$S_{t_{\beta_{2}}}\ti{u}_{\lambda}$ lie 
in the connected component $\ti{\CB}_{0}(\lambda)$ 
of $\bigotimes_{i \in I_{0}} \CB(m_{i}\vpi_{i})$ containing 
$\ti{u}_{\lambda}$, we conclude that the number of elements of 
weight $\lambda$ in $\ti{\CB}_{0}(\lambda)$ is greater than or 
equal to $2$. 
Hence the same is true for the connected component 
$\CB_{0}(\lambda)$ of $\CB(\lambda)$, 
since $\CB_{0}(\lambda)$ is isomorphic to 
$\ti{\CB}_{0}(\lambda)$ as a $P$-crystal.
This proves the claim.

\vsp\vsp

Now suppose that 
the crystal base $\CB(\lambda)$ is isomorphic to 
$\BB(\lambda)$ as a $P$-crystal. 
Then the connected component $\CB_{0}(\lambda)$ of $\CB(\lambda)$ 
is isomorphic as a $P$-crystal to 
some connected component of $\BB(\lambda)$; 
denote it by $\BB_{1}(\lambda)$ as in the proof of 
Lemma~\ref{lem:wt-lam}. 
We know from Lemma~\ref{lem:wt-lam} that 
the number of elements of weight $\lambda$ in 
$\BB_{1}(\lambda)$ is at most $1$. 
Hence we deduce that the number of 
elements of weight $\lambda$ in $\CB_{0}(\lambda)$ 
is at most $1$, which contradicts the claim above. 
This completes the proof of the proposition. 
\end{proof}

\begin{eg}
Assume that $\Fg$ is nontwisted, and 
$\lambda \in P$ is of the form 
$\sum_{i \in J} \vpi_{i}$, 
with $J \subset I_{0}$ and $\#J \ge 2$.
Note that in this case, $d_{i}=1$ for all $i \in I_{0}$ 
by Remark~\ref{rem:cbeta}\,(1) and \eqref{eq:di}. 
We know from \cite[Proposition~5.4\,(ii)]{Kaslz} that 
the crystal base $\CB(\lambda)$ is connected. 
Also, we see from Theorem~\ref{thm:comps} that 
the crystal $\BB(\lambda)$ is connected, 
since $\Supp_{\ge 2}(\lambda) = \emptyset$ and hence 
$\turn(\lambda)=\emptyset$. However, the crystals
$\CB(\lambda)$ and $\BB(\lambda)$ are not isomorphic, 
as can be seen from Proposition~\ref{prop:ext}. 
Let us study this case more precisely.

We know from Theorem~\ref{thm:BN} that 
$\CB(\lambda) \cong \bigotimes_{i \in J} \CB(\vpi_{i})$ 
as $P$-crystals.
In addition, we know from \cite[Corollary~2.2.1]{NSp2} that 
for each $i \in I_{0}$, the crystal base $\CB(\vpi_{i})$ 
is isomorphic as a $P$-crystal to $\BB(\vpi_{i})=\BB_{0}(\vpi_{i})$, 
which in turn is isomorphic to 
the affinization $\ha{\BB(\vpi_{i})_{\cl}}$ of 
the $P_{\cl}$-crystal $\BB(\vpi_{i})_{\cl}$ 
by Theorem~\ref{thm:isom} (note that $d_{\vpi_{i}}=1$). 
Hence we have
%
%
\begin{equation} \label{eq:eg-ext01}
\CB(\lambda) 
\cong 
\bigotimes_{i \in J} 
\ha{\BB(\vpi_{i})_{\cl}}
=
\bigotimes_{i \in J} 
\bigl(\CB(\vpi_{i})_{\cl} \times \BZ \bigr)
 \quad \text{as $P$-crystals}.
\end{equation}
Now, in our case, $d_{\lambda}=1$ 
since by Remark~\ref{rem:dlam} it is equal to 
the greatest common divisor of 
$\bigl\{m_{i}d_{i}\bigr\}_{i \in I_{0}}=
 \bigl\{m_{i}\bigr\}_{i \in I_{0}}$ 
 with $m_{i}=0,\,1$ for $i \in I_{0}$. 
Hence it follows from Theorem~\ref{thm:isom} and 
the connectedness of $\BB(\lambda)$, 
along with Theorem~\ref{thm:NSz}, that 
%
%
\begin{equation} \label{eq:eg-ext02}
\BB(\lambda) = \BB_{0}(\lambda) \cong 
\ha{\BB(\lambda)_{\cl}} =
 \BB(\lambda)_{\cl} \times \BZ \cong 
 \left(\bigotimes_{i \in J} 
 \BB(\vpi_{i})_{\cl}\right) \times \BZ
 \quad \text{as $P$-crystals}.
\end{equation}
The observations \eqref{eq:eg-ext01} and 
\eqref{eq:eg-ext02} above illustrate how 
the $P$-crystals $\CB(\lambda)$ and $\BB(\lambda)$ differ. 
\end{eg}


{\small
\setlength{\baselineskip}{13pt}
\renewcommand{\refname}{References}

}


\begin{thebibliography}{XXXX}

\bibitem[AK]{AK}
T. Akasaka and M. Kashiwara, 
Finite-dimensional representations of 
quantum affine algebras, 
{\it Publ. Res. Inst. Math. Sci.} {\bf 33} (1997), 839--867. 

\bibitem[BN]{BN}
J. Beck and H. Nakajima, 
Crystal bases and two-sided cells of quantum affine algebras, 
{\it Duke Math. J.} {\bf 123} (2004), 335--402.

\bibitem[GL]{GL}
J. Greenstein and P. Lamprou, Path model for quantum loop
modules of fundamental type, {\it Int. Math. Res. Not.} 
{\bf 2004}, no.14, 675--711. 

\bibitem[HK]{HK}
J. Hong and S.-J. Kang, 
``Introduction to quantum groups and crystal bases'', 
Graduate Studies in Mathematics Vol.~42, 
Amer. Math. Soc., Providence, RI, 2002.

\bibitem[J]{J}
A. Joseph,
``Quantum Groups and Their Primitive Ideals'',
Ergebnisse der Mathematik und ihrer Grenzgebiete Vol.~29, 
Springer--Verlag, Berlin, 1995.

\bibitem[Kac]{Kac}
V. G. Kac, 
``Infinite Dimensional Lie Algebras'', 3rd Edition, 
Cambridge University Press, Cambridge, UK, 1990.

\bibitem[Kas1]{Kasmod}
M. Kashiwara, Crystal bases of modified quantized 
enveloping algebra, 
{\it Duke Math. J.} {\bf 73} (1994), 383--413.

\bibitem[Kas2]{Kassim}
M. Kashiwara, Similarity of crystal bases, 
{\it in} ``Lie Algebras and Their Representations'' 
(S.-J. Kang et al., Eds.), Contemp. Math. Vol.~194, 
pp. 177--186, Amer. Math. Soc., Providence, RI, 1996.

\bibitem[Kas3]{Kaslz}
M. Kashiwara, On level-zero representations of 
quantized affine algebras, 
{\it Duke Math. J.} {\bf 112} (2002), 117--175.

\bibitem[L1]{L1}
P. Littelmann, 
A Littlewood--Richardson rule for symmetrizable 
Kac--Moody algebras, 
{\it Invent. Math.} {\bf 116} (1994), 329--346.

\bibitem[L2]{L2}
P. Littelmann, 
Paths and root operators in representation theory, 
{\it Ann. of Math.} (2) {\bf 142} (1995), 499--525.

\bibitem[NS1]{NSp1}
S. Naito and D. Sagaki, 
Path model for a level-zero extremal weight module over a quantum 
affine algebra, {\it Int. Math. Res. Not.} 
{\bf 2003}, no.32, 1731--1754.

\bibitem[NS2]{NSp2}
S. Naito and D. Sagaki, 
Path model for a level-zero extremal weight module over a quantum 
affine algebra. I\hspace{-1pt}I, to appear in {\it Adv. Math.}

\bibitem[NS3]{NSz}
S. Naito and D. Sagaki, 
Crystal of Lakshmibai-Seshadri paths
associated to an integral weight of level zero
for an affine Lie algebra,
{\it Int. Math. Res. Not.} {\bf 2005}, no.14, 815--840.

\end{thebibliography}
\end{document}